\documentclass[reqno,11pt]{amsart}
\usepackage{amssymb}
\usepackage{verbatim}
\usepackage{geometry}
\newif\ifpdf
\ifpdf
  \usepackage[pdftex]{graphicx}
  \usepackage[pdftex]{hyperref}
\else
  \usepackage{graphicx}
\fi
\numberwithin{equation}{section}   
\theoremstyle{plain}
\newtheorem{Thm}{Theorem}[section]
\newtheorem{Prop}[Thm]{Proposition}
\newtheorem{Lemma}[Thm]{Lemma}
\newtheorem{Cor}[Thm]{Corollary}
\newtheorem{Prop-def}[Thm]{Proposition-Definition}

\newtheorem*{ThmA}{Theorem A}
\newtheorem*{ThmB}{Theorem B}
\newtheorem*{ThmC}{Theorem C}
\newtheorem*{ThmAp}{Theorem A'}
\newtheorem*{ThmBp}{Theorem B'}

\theoremstyle{definition}
\newtheorem{Example}[Thm]{Example}
\newtheorem{Def}[Thm]{Definition}
\newtheorem{Remark}[Thm]{Remark}

\newtheorem*{Remark*}{Remark}
\newtheorem*{Ackn}{Acknowledgment}

\newcommand{\C}{{\mathbf{C}}}

\newcommand{\N}{{\mathbf{N}}}
\renewcommand{\P}{{\mathbf{P}}}
\newcommand{\Q}{{\mathbf{Q}}}
\newcommand{\R}{{\mathbf{R}}}
\newcommand{\Z}{{\mathbf{Z}}}

\newcommand{\Rbar}{\overline{\R}_+}

\newcommand{\fm}{{\mathfrak{m}}}
\newcommand{\fp}{{\mathfrak{p}}}
\newcommand{\fC}{{\mathfrak{C}}}

\newcommand{\abar}{{\bar{a}}}
\newcommand{\bbar}{{\bar{b}}}
\newcommand{\cbar}{{\bar{c}}}
\newcommand{\dbar}{{\bar{d}}}

\newcommand{\cC}{{\mathcal{C}}}
\newcommand{\cD}{{\mathcal{D}}}
\newcommand{\cE}{{\mathcal{E}}}
\newcommand{\cF}{{\mathcal{F}}}

\newcommand{\cL}{{\mathcal{L}}}

\newcommand{\cO}{{\mathcal{O}}}

\newcommand{\cS}{{\mathcal{S}}}
\newcommand{\cT}{{\mathcal{T}}}
\newcommand{\cV}{{\mathcal{V}}}

\newcommand{\hF}{{\hat{F}}}
\newcommand{\hf}{{\hat{f}}}

\newcommand{\tm}{{\tilde{m}}}

\newcommand{\talpha}{{\tilde{\alpha}}}
\newcommand{\tbeta}{{\tilde{\beta}}}

\newcommand{\tnu}{{\tilde{\nu}}}

\newcommand{\tA}{{\tilde{A}}}
\newcommand{\tC}{{\tilde{C}}}

\newcommand{\tP}{{\tilde{P}}}

\newcommand{\tU}{{\tilde{U}}}

\newcommand{\tcV}{{\tilde{\mathcal{V}}}}

\renewcommand{\=}{:=}
\renewcommand{\a}{\alpha}

\newcommand{\e}{\varepsilon}

\newcommand{\id}{\operatorname{id}}

\newcommand{\ratrk}{\operatorname{rat.rk}}

\newcommand{\trdeg}{\operatorname{tr.deg}}

\renewcommand{\div}{\operatorname{div}}
\newcommand{\notto}{\nrightarrow}

\newcommand{\ie}{i.e.\ }
\newcommand{\eg}{e.g.\ }
\newcommand{\qand}{{\quad\text{and}\quad}}

\newcommand{\vv}{{\vec{v}}}
\newcommand{\ww}{{\vec{w}}}

\newcommand{\cVqm}{{\mathcal{V}_{\mathrm{qm}}}}
\newcommand{\cVdiv}{{\mathcal{V}_{\mathrm{div}}}}
\newcommand{\Gast}{\Gamma^\ast}

\newcommand{\affroot}{{-\deg}}

\begin{document}
%
%
\setcounter{tocdepth}{1}

\title{Eigenvaluations}
\date{\today}
\author{Charles Favre \and Mattias Jonsson}
\address{CNRS\\
         Institut de Math{\'e}matiques\\
         Equipe G{\'e}om{\'e}trie et Dynamique\\
         F-75251 Paris Cedex 05\\
         France}
\address{IMPA\\
Estrada Dona Castorina 110\\
Rio de Janeiro\\ Brazil 22460-320}
\email{favre@math.jussieu.fr}
\address{Department of Mathematics\\ 
         University of Michigan\\
         Ann Arbor, MI 48109\\
         USA}
\address{Department of Mathematics\\ 
         KTH\\
	 SE-100 44 Stockholm\\
	 Sweden}
\email{mattiasj@umich.edu}
\thanks{Second author supported by the NSF, the Swedish Research
  Council and the Gustafsson foundation.}

%
%

\begin{abstract}
  We study the dynamics in $\C^2$ of superattracting fixed point germs
  and of polynomial maps near infinity.  In both cases we show that
  the asymptotic attraction rate is a quadratic integer, and construct
  a plurisubharmonic function with the adequate
  invariance property.  This is done by finding an infinitely 
  near point
  at which the map becomes rigid: the critical set is contained in a
  totally invariant set with normal crossings.  We locate this
  infinitely near point through the induced action of the dynamics on
  a space of valuations.  This space carries an $\R$-tree structure
  and conveniently encodes local data: an infinitely near point
  corresponds to a open subset of the tree.  The action respects the
  tree structure and admits a fixed point---or eigenvaluation---which
  is attracting in a certain sense.  A suitable basin of attraction
  corresponds to the desired infinitely near point.

\bigskip

\noindent{\sc R\'esum\'e}. Nous \'etudions la dynamique dans $\C^2$
des germes d'applications holomorphes superattractives et des applications
polynomiales pr\`es de l'infini. Dans les deux cas, nous montrons que
le taux asymptotique d'attraction (vers l'origine ou vers l'infini
respectivement) est un entier quadratique, et nous construisons une
fonction plurisousharmonique v\'erifiant l'\'equation d'invariance
correspondante. Pour cel\`a, nous exhibons un point infiniment proche
en lequel l'application devient rigide au sens o\`u son ensemble
critique est contenu dans une courbe \`a croisements normaux et
totalement invariante. Nous localisons ce point infiniment proche en
analysant l'action induite par l'application sur un espace ad\'equat
de valuations. Cet espace est un arbre r\'eel et code de mani\`ere
efficace des donn\'ees locales li\'ees aux singularit\'es de courbes:
un point infiniment proche correspond ainsi \`a un ouvert dans
cet arbre. L'action respecte la structure d'arbre et admet un point
fixe (appel\'e ``eigenvaluation'') qui est attirant en un certain
sens. Le bassin d'attraction de ce point correspond au point
infiniment proche d\'esir\'e.
\end{abstract}
%
%
\maketitle
%
%
%
%
%
%
%
%
\section*{Introduction}
Our objective in this paper is to study dynamics in $\C^2$ in two
situations: holomorphic fixed point germs $f$ at the origin and
polynomial maps $F$ at infinity.  We shall refer to these two cases as
\emph{local} and \emph{affine},\footnote{The affine case refers to
$\C^2$ with a fixed algebraic embedding into $\P^2$. The local case
applies to a (classical) neighborhood of any smooth point on a
surface.}  and we are interested in how fast orbits approach the
origin and infinity, respectively.

In the affine case we can measure the rate of attraction to
infinity by the algebraic degree $\deg(F)$.  Indeed, for ``most''
points $p$ near infinity we have $\|F(p)\|\simeq\|p\|^{\deg(F)}$.
It is easy to see that $d_\infty:=\lim_{n\to\infty}\deg(F^n)^{1/n}$
exists.\footnote{This number is also known as 
  the (first) dynamical degree of $F$, see~\cite{RS}.}  
However, it is not obvious in what sense
that orbits are actually attracted at the rate~$d_\infty$.

Similarly, in the local case, we define $c(f)\ge1$ as the
smallest degree of any term in the Taylor expansion of $f$ in
local coordinates.
Thus $\|f(p)\|\simeq\|p\|^{c(f)}$ for ``most'' $p$ near the
origin.  Again, the existence of the limit
$c_\infty=\lim_{n\to\infty}c(f^n)^{1/n}$ is easy to establish, 
but its dynamical significance is less clear.

To avoid trivialities, all maps considered will be 
\emph{dominant}, \ie their Jacobian determinants 
do not vanish identically.
We are also mainly interested in superattracting behavior, so
we shall assume throughout
that $d_\infty>1$ and $c_\infty>1$, respectively.

Our aim is then to answer two questions. 
First, what numbers $d_\infty$ and $c_\infty$ may appear,
and what can be said about the convergence of 
$\deg(F^n)^{1/n}$ and $c(f^n)^{1/n}$ to $d_\infty$ and $c_\infty$,
respectively? 
Second, in what sense are orbits actually attracted to
infinity (the origin) at the rate $d_\infty$ ($c_\infty$)?

The next two theorems provide answers in the local case.
\begin{ThmA}
  For any dominant fixed point germ $f:(\C^2,0)\to(\C^2,0)$,
  the asymptotic attraction rate $c_\infty=c_\infty(f)$ 
  is a quadratic integer,
  \ie there exist integers $a$, $b$ such that
  $c_\infty^2+ac_\infty+b=0$.
  
  Moreover, there exists $\delta\in(0,1]$ such that
  $\delta c_\infty^n\le c(f^n)\le c_\infty^n$ for all $n\ge1$.
\end{ThmA}
Our proof also goes through for formal
fixed points germs (defined by formal power series).
The bound on $c(f^n)$ allows us to improve upon Theorem~A'
from~\cite{brolin} (see the remark after that theorem).
\begin{ThmB}
  If $f:(\C^2,0)\to(\C^2,0)$ is a dominant
  fixed point germ with $c_\infty(f)>1$, 
  then there exists a plurisubharmonic
  function $u$ defined near the origin and a constant $C>0$
  such that $C^{-1}\log\|\cdot\|\le u\le C\log\|\cdot\|$
  and such that $c_\infty^{-n}u\circ f^n$ decreases to
  a plurisubharmonic function $u_\infty\not\equiv-\infty$
  with $u_\infty\circ f=c_\infty u_\infty$.
\end{ThmB}
In particular, for any point $p$ outside a pluripolar set, we have
$\delta_1^{c_\infty^n}\le\|f^n(p)\|\le \delta_2^{c_\infty^n}$ for
all $n\ge1$ and some constants $\delta_i= \delta_i(p) >0$.
Of course, not every orbit must converge to the origin
at the rate $c_\infty$ and $f$ may even contract some curves.

Hubbard and Papadopol~\cite{HP1} also studied superattracting fixed
point germs, but their objective was different from ours: for the
mappings they considered, $c(f^n)=c^n$ and the
existence of $u_\infty$ is easier to establish.

\smallskip
The corresponding results in the affine case are as follows.
\begin{ThmAp}
  For any dominant polynomial map $F:\C^2\to\C^2$,
  the asymptotic attraction rate $d_\infty=d_\infty(F)$ 
  is a quadratic integer.
  
  Moreover, if $d_\infty>1$, then 
  either there exists $D\in[1,\infty)$ such that
  $d_\infty^n\le \deg(F^n)\le Dd_\infty^n$ for all $n\ge1$; or there
  exists a change of coordinates by a polynomial automorphism of
  $\C^2$ such that $F$ becomes a skew product, 
  $F(X,Y)=(P(X),Q(X,Y))$,
  with $d_\infty(F)=\deg(P)=\deg_Y(Q)>\deg_Y Q(X_0,Y)$ 
  for some $X_0\in\C$. 
\end{ThmAp}
\begin{ThmBp}
  Suppose $F:\C^2\to\C^2$ is a dominant polynomial map with
  $d_\infty(F)>1$ and not conjugate to a skew product.
  Then there exist $C>0$ and a plurisubharmonic function $U\ge0$
  on $\C^2$,
  such that $C^{-1}\log\|\cdot\|\le U\le C\log\|\cdot\|$
  near infinity, and such that $d_\infty^{-n}U\circ F^n$ decreases to
  a plurisubharmonic function $U_\infty\not\equiv0$ with
  $U_\infty\circ F=d_\infty U_\infty$. 
  The positive closed current $T\=dd^cU_\infty$ 
  satisfies $F^*T=d_\infty T$.
\end{ThmBp}
Skew products were considered in~\cite[Section~6]{FG}.
The construction of an invariant current for rational maps is
due to Sibony~\cite{Sibony} in the
``algebraically stable'' case,
when $\deg(F^n) = \deg(F)^n$ for all $n\ge1$.
This was later extended to birational surface maps
in~\cite{DF}, and to special non-invertible rational maps in~\cite{Ng}. 
All these constructions are based on a precise understanding 
on the growth of degrees; the same is true here.

The set $\{U_\infty=0\}$ above is not pluripolar in general
and points in $\{U_\infty=0\}$ may converge to infinity 
at a speed different from $d_\infty$: 
see~\cite{DDS,Vigny} for interesting (algebraically stable) examples.

Let us now explain our approach to the main results above.  
We first focus on the local case where the ideas are less involved.

A common approach to study dynamic behavior
is to find coordinates in 
which the mapping has a particularly simple 
expression---a \emph{normal form}. 
In one dimension, a very simple normal form is provided
by Boettcher's Theorem (see~\cite{CG}): there exists
a local coordinate near the origin such that the map
becomes $\zeta\mapsto\zeta^c$ for some $c\ge2$. 

As noted by Hubbard and Papadopol, the naive
generalization of Boet\-tcher's Theorem fails in higher dimensions.
For instance, $f(x,y)=(x^2,y^2-x^3)$ cannot be conjugated
to $f_0(x,y)=(x^2,y^2)$. 
Indeed, the critical set $C_{f_0}=\{xy=0\}$ is totally invariant
for $f_0$, but $C_f=\{xy=0\}$ is not totally invariant for $f$.
For this example, the analysis by Hubbard and Papadopol 
(see also~\cite{BJ} for the affine case) provides a 
conjugacy between $f$ and $f_0$ on a subset of $\C^2$, 
but that relies on much stronger (hyperbolicity)
assumptions on $f$ than we wish to impose.

Instead, we take inspiration from singularity theory and allow
ourselves to work with birational changes of coordinates
and with slightly more general normal forms.

Let us call a fixed point germ \emph{rigid} if the critical set is
contained in a totally invariant set with normal crossings
singularities. Thus $f_0$ above is rigid while $f$ is not.  Rigid
(contracting) germs were introduced and classified
in~\cite{F-rigid}.  There are seven classes of them, each containing
a simple normal form.

The concept of a rigid germ is a dynamic version of
simple normal crossings singularities for curves. 
The next theorem, which we shall refer to as \emph{rigidification},
can thus be viewed as an analogue of embedded resolution 
of plane curve singularities. 
By a \emph{modification} we shall mean a composition
of point blowups above the origin. A (closed) point on the exceptional
divisor is called an \emph{infinitely near point}.
\begin{ThmC}
  For any dominant fixed point germ $f:(\C^2,0)\to(\C^2,0)$
  with $c_\infty(f)>1$,  there exists a
  modification $\pi:X\to(\C^2,0)$ and an infinitely near 
  point $p\in\pi^{-1}(0)$ such that the lift
  $\hf:=\pi^{-1}\circ f\circ\pi$ is
  holomorphic at $p$, $\hf(p)=p$, 
  and the fixed point germ $\hf:(X,p)\to(X,p)$ is rigid. 
\end{ThmC}
If we allow ourselves to use \emph{different} modifications at the
source and target spaces, we can make any map 
holomorphic and locally \emph{monomial} 
in suitable coordinates, see~\cite{Cut} 
and Proposition~\ref{P-monomap} below.  Although this fact is not a
direct consequence of Hironaka's desingularization theorem, it is
considerably simpler to establish than Theorem~C, 
and serves as an important ingredient in the proof of the latter.  

We now indicate how Theorem~C can be used to prove Theorems~A and~B.
One can show that $c_\infty(f)$ is controlled by the speed at which
orbits at $p$ under the lift $\hf$ approach the exceptional divisor
$\pi^{-1}(0)$. 
In view of the classification of rigid germs, this speed is
either an integer or the spectral radius of a
$2\times 2$ matrix with integer coefficients, hence a quadratic
integer.

The construction of the function $u$ 
with $u\circ f\le c_\infty u$ in Theorem~B
requires a separate
argument, but once we have found $u$, the fact that 
$c_\infty^{-n}u\circ f^n\not\to-\infty$ can
be established by showing that at least one orbit of 
$\hf$ cannot approach $\pi^{-1}(0)$ too fast.

\smallskip 
So how, then, could one prove rigidification
as in Theorem~C?  When resolving curve singularities, the
modification $\pi:X\to(\C^2,0)$ can be found inductively by
successively blowing up points where we do not have normal
crossings.  The main problem then becomes proving that the inductive
procedure terminates: this is done by showing that the curves involved
always become less singular in the sense that a certain discrete
numerical invariant decreases.

To prove rigidification, such an approach seems difficult to
carry out successfully.
Instead we take a different tack: we consider a space
parameterizing \emph{all} possible modifications at the same time,
study the dynamics induced by $f$ on this space, 
and deduce dynamical information about $f$ itself.

To explain this strategy, note that any irreducible component $E$
of the exceptional divisor of a modification 
defines a \emph{divisorial valuation} on the ring $R$ 
of holomorphic germs at the origin: the value on a germ 
is the order of vanishing along $E$ of its pullback.
It is possible to define a natural induced action of $f$
on the space $\cVdiv$ of divisorial valuations.
However, this space has a lack of completeness which makes it hard to
properly analyze the induced action.

Instead, we consider a larger space $\cV$ consisting of \emph{all}
valuations on $R$ (centered at the origin, and suitably
normalized). This space is automatically compact for the topology of
the weak convergence, and we showed in~\cite{valtree} that $\cV$
admits a natural \emph{tree structure}. Roughly speaking, $\cV$ is a
collection of real line segments welded together in such a way that no
cycles appear.  More precisely, $\cV$ is an $\R$-tree for a natural
metric and possesses strong self-similarity properties. The line
segments are made up of monomial valuations in suitable coordinates at
infinitely near points above the origin, and are related to the Farey
blowups studied in~\cite{HP2}. The set $\cVdiv$ above is a
subset of $\cV$ in essentially the same way as $\Q$ is a subset of the
(compactified) real line.  In the sequel, we refer to $\cV$ as the
\emph{valuative tree}.

The valuative tree can be used to
efficiently encode local information near the origin in $\C^2$.  
We used it to study singularities of curves and ideals
in~\cite{valtree}, of plurisubharmonic functions in~\cite{pshsing},
and multiplier ideals in~\cite{valmul}.

Here we prove that any dominant fixed point germ
$f$ induces a natural selfmap $f_\bullet:\cV\to\cV$, and
that this selfmap respects the tree structure in a strong sense.
For topological reasons, $f_\bullet$ must then admit a
fixed point, or \emph{eigenvaluation}. 
We show that this
eigenvaluation admits a basin of attraction $U$, which can 
be chosen to be of a special type $U(p)$ 
corresponding to an infinitely near
point $p$ as follows.  Fix a modification $\pi:X\to(\C^2,0)$ such that
$p\in\pi^{-1}(0)\subset X$.
Then $U(p)\subset\cV$ is the set of 
all valuations whose center on $X$ is $p$.
The invariance of $U(p)$ under $f_\bullet$ translates into
the lift $\hf$ of $f$ being holomorphic at $p$ with $\hf(p)=p$.
By choosing the basin carefully, we can even make $\hf$ rigid at $p$,
leading to Theorem~C. 

The induced selfmap $f_\bullet$ of $\cV$ seems to have quite
strong global attracting features, similar to 
holomorphic selfmaps of the unit disk in $\C$.
There are also close connections to $p$-adic dynamics, more precisely
to the study of rational maps on $\C_p$ (see~\cite{Rivera1}).

\smallskip
In the affine setting, \ie when proving Theorems~A' and~B', the main
idea is the same as in the local case, and consists of studying
the action of $F$ on the set $\cV_0$ of (normalized) valuations on
$\C[X,Y]$ centered at infinity. However, in the same way that the
study of plane affine curves has a quite different flavor than the
theory of germs of analytic curves, the affine case proves to be more
delicate than the local case. Let us indicate two
difficulties that we need to overcome.

First, a dominant polynomial map of $\C^2$ need not be proper, so the
line at infinity is not invariant.  In particular, the pushforward of
a valuation centered at infinity may no longer be centered at
infinity.  As a consequence, the induced selfmap of $\cV_0$ is not
everywhere defined.  Second (and more importantly), a valuation is by
its very nature a local object and measures behavior at an infinitely
near point.  But for Theorem~A' we need to control degrees, which has
a more global flavor.

As it turns out, we can address both difficulties 
by working with a suitable subtree $\cV_1$ of $\cV_0$.  On this
tree, the relation between local and global information is well
behaved.  Moreover, any dominant polynomial
map acts on this tree in a way that preserves the tree structure.  In
particular, there is again an eigenvaluation, and this provides us
with a starting point for the proofs of Theorems~A' and~B'.  The case
when $F$ can be conjugated to a skew product arises exactly when the
eigenvaluation is an end of the tree of a special type: it is
associated to a pencil of rational curves having one place at
infinity. The Line Embedding Theorem~\cite{abh-moh,Suzuki} then
provides the conjugacy.

The space $\cV_1$ is the key new idea in the affine case. 
Roughly speaking, a valuation is in $\cV_1$ if it takes negative
values on any polynomial as well as on the form $dX\wedge dY$ 
on $\C^2$. 
To our knowledge, such valuations have not 
previously been systematically studied. 
A crucial property that they enjoy is that they are
dominated by affine curves with one place at infinity.
While such curves have received much attention in the literature
(see \eg~\cite{abh-moh,Suzuki,CPR,FujimotoSuzuki}),
our focus on valuations made it difficult 
to extract the information that we needed.
We therefore explore the structure on $\cV_1$ in an appendix, using
key polynomials as in~\cite{valtree}, a technique originating with
MacLane~\cite{mac} and subsequently developed by 
Abhyankar-Moh~\cite{abh-moh} under the name approximate roots.

\smallskip
The paper is divided into eight sections. In the first one, we
briefly recall the structure of the valuative tree as described
in~\cite{valtree}.  In the next five sections, we only deal with the
local case, leaving affine problems for Section~\ref{sec-affine}. 
We thus fix a fixed point germ $f$.  
In Section~\ref{sec-actsingle},
we describe how $f$ maps one valuation to another. 
We then investigate the action
of $f$ (without iteration) on the valuative tree in
Section~\ref{sec-acttree}, and show that it 
preserves---in a strong sense---the tree structure. 
Section~\ref{sec-fixed} is devoted to
the existence of an eigenvaluation and the second part of Theorem~A. 
In Section~\ref{sec-normal}, we construct basins of attraction,
deduce a stronger version of Theorem~C, 
and exhibit normal forms of the lift $f$ at the point where
it is rigid. This enables us to complete the proof of
Theorem~A. Finally, Section~\ref{sec-green} is devoted to the proof of
Theorem~B, which relies heavily on Theorem~C. 
In Section~\ref{sec-affine}, we turn to the affine case. 
The paper ends with appendix detailing 
the structure of the valuation spaces $\cV_0$ and $\cV_1$ 
in the same spirit as the local case in~\cite{valtree}. 
%
%
%
%
\begin{Ackn}
  This work was partially carried out when the authors were
  visiting RIMS, Kyoto. We thank this institute for its
  hospitality, and Prof.~Shishikura for his exceptional kindness.
  We also thank R.~A.~Mollin for suggesting the proof of
  Lemma~\ref{L-quadra} and the referee for useful remarks.
\end{Ackn}
%
%
%
%
\section{Background}\label{sec-background}
In this section we briefly review the approach to valuation
theory worked out in~\cite{valtree}.
Except in Section~\ref{sec-affine} and Appendix~\ref{sec-app}
we let $R$ denote the ring of holomorphic germs at the origin in $\C^2$. 
This is a local ring. 
Its maximal ideal $\fm$ is the set of germs vanishing
at the origin, and its residue field is $\C$.
We write $(\hat{R},\hat{\fm})$ for the completion of $R$.
It is the ring of formal power series in two complex variables.
%
%
\subsection{Valuations}
We consider the space $\cV$ of centered, normalized 
\emph{valuations} on $R$,
\ie the set of functions $\nu:R\to[0,\infty]$ satisfying:
\begin{enumerate}
\item[(i)]\label{p1} 
  $\nu(\psi\psi')=\nu(\psi)+\nu(\psi')$ for all $\psi,\psi'$; 
\item[(ii)]\label{p2}
  $\nu(\psi+\psi')\ge\min\{\nu(\psi),\nu(\psi')\}$ for all
  $\psi,\psi'$; 
\item[(iii)]\label{p3} 
  $\nu(0)=\infty$, 
  $\nu|_{\C^*}=0$, 
  $\nu(\fm):=\min\{\nu(\psi)\ ;\ \psi\in\fm\}=1$.
\end{enumerate} 
Then $\cV$ is equipped with a natural \emph{partial ordering}: 
$\nu\le\mu$ iff $\nu(\psi)\le\mu(\psi)$ for all $\psi\in\fm$. 
The \emph{multiplicity valuation} $\nu_\fm$ defined by
$\nu_\fm(\psi)=m(\psi)=\max\{k\ ;\ \psi\in\fm^k\}$ 
is the unique minimal element of $\cV$. 
%
%
\subsection{Curve valuations}\label{sec-curve}
Some natural maximal elements (for the partial ordering above)
are the \emph{curve valuations} defined as follows. To each
irreducible curve $C$ we associate $\nu_C\in\cV$ defined by
$\nu_C(\psi)=C\cdot\psi^{-1}\{0\}/m(C)$, where ``$\cdot$'' denotes
intersection multiplicity and $m$ multiplicity.  We call $\nu_C$
\emph{analytic} when $C$ is an analytic curve; otherwise $\nu_C$ is
\emph{formal}.  If $C=\phi^{-1}(0)$ for $\phi\in\hat\fm$, then we also
write $\nu_C=\nu_\phi$.  Given a parameterization
$h:(\C,0)\to(\C^2,0)$ of $C$ we have
$\nu_C(\psi)=m(C)^{-1}\div_t(\psi\circ h(t))$, where $\div_t$ denotes
the order of vanishing at $t=0$.  Note that $\nu_\phi(\phi)=\infty$.

The set $\cC$ of local irreducible curves carries a natural
(ultra)metric in which $\cC$ has diameter 1.
It is given by $d_\cC(C,D)=m(C)m(D)/C\cdot D$.
%
%
\subsection{Quasimonomial valuations and tree structure}\label{sec-qm}
Arguably the most important valuations in $\cV$ are the
\emph{quasimonomial} ones.
They are of the form
$\nu_{C,t}$, where $C\in\cC, t\in[1,\infty)$, and satisfy
$\nu_{C,t}(\psi)=\min\{\nu_D(\psi)\ ;\ d_\cC(C,D)\le t^{-1}\}$.  We
have $\nu_{C,s}=\nu_{D,t}$ iff $s=t\le d_\cC(C,D)^{-1}$; and
$\nu_{C,t} \ge \nu_{C,t'}$ iff $t \ge t'$.  Thus $\cVqm$, the set of
all quasimonomial valuations, is naturally a quotient of
$\cC\times[1,\infty)$, and has a natural tree structure: if
$\nu,\nu'\in\cVqm$ and $\nu<\nu'$, then the \emph{segment}
$[\nu,\nu']=\{\mu\in\cVqm\ ;\ \nu\le\mu\le\nu'\}$ is isomorphic to a
compact real interval. 
We set $\nu_{\phi,t}:=\nu_{C,t}$ when
$C=\phi^{-1}(0)$.  
Quasimonomial valuations are of two types: \emph{divisorial}
and \emph{irrational}, depending on whether the parameter $t$
is rational or irrational.

The full space $\cV$ is the completion of $\cVqm$ in the sense
that every element in $\cV$ is the limit of an increasing sequence
in $\cVqm$. It is hence also naturally a tree, called the 
\emph{valuative tree}. 
The \emph{ends} (\ie maximal elements) of $\cV$ 
are exactly the elements
of $\cV\setminus\cVqm$ and are either curve valuations or
\emph{infinitely singular} valuations.\footnote{The latter are
represented by infinite Puiseux series whose exponents are
rational numbers with unbounded denominators.}

A \emph{subtree} of $\cV$ is a subset $\cT$ such that $\nu\in\cT$
and $\mu\le\nu$ implies $\mu\in\cT$. A subtree is \emph{finite}
if it has finitely many ends.
%
%
\subsection{Numerical invariants}\label{sec-numerical}
A curve valuation $\nu=\nu_\phi$ is characterized by the fact that
the prime ideal $\fp_\nu:=\{\nu=+\infty\}$ is nontrivial, generated 
by $\phi$.

If $\nu\in\cV$ is not a curve valuation, then its \emph{type}
(divisorial, irrational or infinitely singular) is determined by its
\emph{rational rank} $\ratrk\nu$ and \emph{transcendence degree}
$\trdeg\nu$, see~\cite[Table~C.1]{valtree}. These invariants are
defined as follows.
The semigroup $\nu(R\setminus\{ 0\})\subset\R$
generates a group $\Gamma_\nu\subset\R$, the value group of $\nu$. 
One sets $\ratrk\nu=\dim_\Q(\Gamma_\nu\otimes_\Z\Q)$.
The residue field of the valuation is 
$k_\nu\=\{\phi\in R\ ;\ \nu(\phi)\ge0\}/\{\phi\in R\ ;\ \nu(\phi)>0\}$
and $\trdeg\nu$ is its transcendence degree over $\C$. 
%
%
\subsection{Skewness and infimum}\label{sec-skewness}
An important invariant of a valuation is its \emph{skewness} $\a$
defined by $\a(\nu)=\sup\{\nu(\phi)/m(\phi)\ ;\ \phi\in\fm\}$.
Skewness naturally \emph{parameterizes} the trees $\cVqm$ and $\cV$ in
the sense that $\a:\cVqm\to[1,\infty)$ is strictly increasing and
restricts to a bijection onto its image on any segment; indeed
$\a(\nu_{\phi,t})=t$ for any $\nu_{\phi,t}\in\cVqm$. 
Thus divisorial (irrational) valuations have rational
(irrational) skewness. Curve valuations have 
infinite skewness: see~\cite[Table~C.1]{valtree}.

The tree structure on $\cV$ implies that any collection 
$(\nu_i)_{i\in I}$ of valuations in $\cV$ admits an 
\emph{infimum} $\wedge_i\nu_i$. 
There are two useful formulas involving the infimum and skewness.
First, if $\nu\in\cV$ and $\phi\in\fm$ is irreducible, then 
$\nu(\phi)=m(\phi)\a(\nu\wedge\nu_\phi)$.
Second, if $C,D\in\cC$, then 
$C\cdot D=m(C)m(D)\a(\nu_C\wedge\nu_D)$.
%
%
\subsection{Tangent space and weak topology}\label{sec-tangent}
Let $\nu$ be a valuation in $\cV$. Declare
$\mu,\mu'\in\cV\setminus\{\nu\}$ to be equivalent if the segments
$]\nu,\mu]$ and $]\nu,\mu']$ intersect.  An equivalence class is
called a \emph{tangent vector} at $\nu$ and the set of tangent vectors
at $\nu$, the \emph{tangent space}, denoted by $T\nu$. If $\vv$ is a
tangent vector, we denote by $U(\vv)$ the set of points in $\cV$
defining the equivalence class $\vv$. The points in $U(\vv)$ are said
to \emph{represent} $\vv$.  Irrational valuations are \emph{regular
points} of $\cV$ in the sense that the tangent space has exactly two
elements, whereas divisorial points are \emph{branch points}:
see~\cite[Table~C.1]{valtree}.  We refer to Section~\ref{S21} for a
geometric interpretation of the (uncountable) tangent space at a
divisorial valuation.

We endow $\cV$ with the \emph{weak topology},
generated by the sets $U(\vv)$ over all tangent vectors $\vv$;
this turns $\cV$ into a compact (Hausdorff) space.
The weak topology on $\cV$ is
characterized by $\nu_k\to\nu$ iff $\nu_k(\phi)\to\nu(\phi)$
for all $\phi\in R$.
%
%
\subsection{Multiplicity and thinness}\label{sec-mult}
By setting $m(\nu)=\inf\{m(C)\ ;\ C\in\cC,\ \nu_C\ge\nu\}$ we
extend the notion of \emph{multiplicity} from $\cC$ to $\cV$.  
The infinitely
singular valuations are characterized as having infinite multiplicity.

As $m$ is increasing and integer valued, it is piecewise 
constant on any segment $[\nu_\fm,\nu_C]$, where $C\in\cC$.
This naturally defines $m(\vv)$ for any tangent
vector $\vv$. If $\nu$ is nondivisorial, then $m(\vv)=m(\nu)$ 
for any $\vv\in T\nu$. 
When $\nu$ is divisorial, the situation is more complicated.
Suffice it to say that there exists an integer $b(\nu)$,
divisible by $m(\nu)$, such that $m(\vv)=b(\nu)$ for all
but at most two tangent vectors $\vv$ at $\nu$. We call
$b(\nu)$ the \emph{generic multiplicity} of $\nu$.

A quantity that combines skewness and multiplicity
is \emph{thinness} $A$, defined by
$A(\nu)=2+\int_{\nu_\fm}^\nu m(\mu)\,d\alpha(\mu)$.
The integral is in fact a finite or countable sum as $\mu\mapsto
m(\mu)$ is increasing and integer-valued on $[\nu_\fm,\nu[$. Thinness
defines an increasing parameterization of $\cV$ with values in
$[2,+\infty]$.

If $\nu_C$ is a curve valuation, and $\nu_1<\nu_C$ is 
a quasimonomial valuation 
such that $m\=m(\nu_1)=m(C)$, then 
the multiplicity is constant on $[\nu_1,\nu_C]$, hence
$A(\nu)= m\a(\nu)+B$ for some rational number $B$.
%
%
\subsection{Tree potentials}\label{sec-treepot1}
The set of \emph{tree potentials} is by definition the smallest
set of real-valued functions on $\cVqm$ containing all functions of
the form $\nu\mapsto\nu(\phi)$ for $\phi\in R$, and closed under
positive linear combinations, minima, and pointwise limits. 
It is the analog for a real tree of the set of concave functions 
on the real line.

The tree potentials that we shall use are associated to
germs $\phi\in R$ 
or ideals $I\subset R$.
In the first case, write $\phi=\prod\phi_i^{m_i}$ with $\phi_i$
irreducible. Then $\nu(\phi)=\sum m_i\,\a(\nu\wedge\nu_{\phi_i})$. 
In particular, the tree potential $\nu\mapsto\nu(\phi)$
is locally constant outside the finite tree whose ends are
the curve valuations $\nu_{\phi_i}$;
and is increasing, concave and piecewise linear with integer slope 
as a function of skewness on any segment. 

These properties carry over to the tree potential $\nu\mapsto\nu(I)$
associated to an ideal $I$, as follows from 
$\nu(I)=\min_j\nu(\phi_j)$, where $(\phi_j)$ is any 
finite set of generators.
%
%
\subsection{Dual graphs}\label{sec-dual}
Every divisorial valuation $\nu$ arises as follows: there exists a
\emph{modification} (\ie a finite composition of point blowups)
$\pi:X\to(\C^2,0)$ and an \emph{exceptional component} $E$ (\ie an
irreducible component of $\pi^{-1}(0)$) such that
the center of $\nu$ on $X$ equals $E$.
More precisely, $\nu=\nu_E:=b_E^{-1}\pi_*\div_E$,  
where $b_E=\div_E(\pi^* \fm)\in\N^*$ is the generic multiplicity 
at $\nu$.
The ideal $\pi^*\fm$ is locally principal along $E$, 
of the form $(z^{b_E})$ if $E=\{z=0\}$. 
Further, the thinness of $\nu$ is $A(\nu)=a_E/b_E$, 
where $a_E=1+\div_E(\pi^*\omega)$
and $\omega$ is a nonvanishing holomorphic two-form, say
$\omega=dx\wedge dy$ in local coordinates at $(\C^2,0)$.

To a modification $\pi:X\to(\C^2,0)$ is associated a 
\emph{dual graph}: its set of vertices $\Gast_\pi$ is
the set of exceptional components $E\subset\pi^{-1}(0)$ and
$E,F\in\Gast_\pi$ are joined by an edge iff $E\cap F\ne\emptyset$.
Since the dual graph has no cycles, $\Gast_\pi$ admits a natural
poset structure, with unique minimal element given by the
blowup of the origin. 
Then $\Gast_\pi$ embeds in 
$\cV$ in the sense that $E\mapsto\nu_E$ is injective
and $E \in [F,F']$ iff $\nu_E\in [\nu_F,\nu_{F'}]$.

More generally, for any curve $C$ and any 
modification $\pi:X\to(\C^2,0)$
such that $\pi^{-1}(C)$ has \emph{simple normal crossings},
the set of irreducible components of $\pi^{-1}(C)$ has
a natural poset structure and embeds in $\cV$ as a set
of divisorial and curve valuations. 
See~\cite[Figure~9]{valtree}.
%
%
\subsection{Infinitely near points}\label{S21}
An \emph{infinitely near point} is a (closed) point
$p\in\pi^{-1}(0)$ where $\pi:X\to(\C^2,0)$ is
a modification. Such a point defines an open subset
$U(p)\subset\cV$, consisting of all valuations
$\nu\in\cV$ whose center on $X$ is $p$, \ie
$\nu=\pi_*\mu$, where $\mu$ is a valuation on 
the local ring $(\cO_p,\fm_p)$ of holomorphic germs 
at $p$ with $\mu\ge0$ on $\cO_p$ and $\mu>0$ on $\fm_p$.

If $p$ is a smooth point on $\pi^{-1}(0)$, belonging to
a unique exceptional component $E$, then
$U(p)=U(\vv)$, where $\vv$ is the tangent vector at $\nu_E$
represented by $\nu_F$ and $F$ is the exceptional divisor
obtained by blowing up $p$. 
If $p$ is a singular point on $\pi^{-1}(0)$, say 
$\{p\}=E\cap E'$, then $U(p)=U(\vv)\cap U(\vv')$, where 
$\vv$ is the tangent vector at $\nu_E$ represented by $\nu_{E'}$ and 
$\vv'$ is the tangent vector at $\nu_{E'}$ represented by $\nu_E$.
See~\cite[Figure~9]{valtree}.

It follows that if $p$ is any point on an exceptional component 
$E\subset\pi^{-1}(0)$, then all valuations in 
$U(p)$ represent a common tangent vector $\vv_p$ at $\nu_E$.
%
%
\subsection{Monomialization of valuations}\label{S22}
A valuation $\nu$ on $R$ is \emph{monomial} if it is 
monomial in some local coordinates $(x,y)$, \ie 
if there exist $s,t\ge0$ such that
$\nu(\sum_{i,j}a_{ij}x^iy^j)=\min\{si+tj\ ;\ a_{ij}\ne0\}$.
In particular $\nu(x)=s$, $\nu(y)=t$.

We shall repeatedly reduce questions regarding quasimonomial
valuations to the monomial case using modifications.

First, consider a modification $\pi:X\to(\C^2,0)$ and an
infinitely near point $p\in\pi^{-1}(0)$
at the intersection of two exceptional components $E=\{z=0\}$ and
$F=\{w=0\}$ in $X$.  Let $\nu_E,\nu_F\in\cV$ be the associated
divisorial valuations and $b_E,b_F$ their generic multiplicities.
Denote by $\mu_t$ the monomial valuation in $(z,w)$ with
$\mu_t(z)=(1-t)/b_E$ and $\mu_t(w)=t/b_F$.  The ideal $\pi^* \fm$ is
locally principal, generated by $z^{b_E}w^{b_F}$ at $p$
(see Section~\ref{sec-dual} above), so $\nu_t:=\pi_*\mu_t$ is normalized by
$\nu_t(\fm)=1$. Moreover $[0,1]\ni t\mapsto\nu_t$ gives a
parameterization of the segment $[\nu_E,\nu_F]$ in $\cV$.

Second, consider an irreducible curve $C$, a desingularization
$\pi:X\to(\C^2,0)$ of $C$, and pick local coordinates $(z,w)$ at the
intersection point $p$ between the strict transform $\tC=\{w=0\}$ and
the unique exceptional component $E=\{z=0\}$ intersecting $\tC$.
Define a monomial valuation $\mu_t$ in $(z,w)$ by $\mu_t(z)=1/b_E$,
$\mu_t(w)=t$.  Then $\nu_t:=\pi_*\mu_t$ is normalized and 
$[0,\infty]\ni t\mapsto\nu_t$ 
gives a parameterization of the segment $[\nu_E,\nu_C]$ in $\cV$.
%
%
\subsection{Monomialization of maps}\label{S230}
A holomorphic map $f:(\C^2,0)\to(\C^2,0)$ is 
\emph{monomial} if one can find coordinates $(z,w)$
at the source and $(z',w')$ at the target space such that 
$(z',w')=f(z,w)=(z^aw^b,z^cw^d)$ for some integers $a,b,c,d$. 
Then $f$ is dominant iff $ad\ne bc$.
\begin{Prop}\label{P-monomap}
  For any dominant holomorphic map $f:(\C^2,0) \to (\C^2,0)$ and for
  any curve $C$, there exist modifications $\pi:X\to(\C^2,0)$ and
  $\pi': X'\to(\C^2,0)$ such that the total transform of $C$ by $\pi$
  and of $f(C)$ by $\pi'$ both have simple normal crossings,
  and such that the lift $\hf:X\to X'$ of $f$ is
  holomorphic and monomial at any point $p\in\pi^{-1}(0)$.
\end{Prop}
A proof is given in~\cite[Section~3]{Cut}. 
While the result is not known in high dimensions, 
its proof is elementary in dimension two.
%
%
%
%
\section{Image of a single valuation}\label{sec-actsingle}
From now on and until Section~\ref{sec-affine}, 
$f:(\C^2,0)\to(\C^2,0)$ denotes
a holomorphic fixed point germ.
We will always assume $f$ \emph{dominant}, that is,
the Jacobian determinant $Jf$ is not
identically zero.
In this section, we describe how $f$ maps one valuation to another.
%
%
\subsection{Definitions}
An element $\nu\in\cV$ is by definition a 
valuation $\nu:R\to\Rbar$ centered at the origin 
and normalized by $\nu(\fm)=1$.
We let $f_*\nu(\phi) \= \nu ( f^* \phi)$.
This defines $f_*\nu$ as a centered valuation on $R$ 
with values in $\Rbar$. However, 
$f_*\nu$ need not be normalized, and can even be
identically $+\infty$ on $\fm$. 
The latter situation appears exactly when
$f^*\fm$ is included in the prime ideal 
$\fp_\nu=\{\nu=+\infty\}$. When nontrivial, $\fp_\nu$ 
is generated by an irreducible
element, $\fp_\nu=(\phi)$, and $\nu=\nu_\phi$ is a curve valuation.
Geometrically, the inclusion $f^*\fm\subset(\phi)$ 
means that the curve $C=\{\phi=0\}$ is contracted, $f(C)=0$.
We therefore call $\nu_\phi=\nu_C$ a 
\emph{contracted curve valuation} in this case.
Any contracted curve is critical, so 
there are most finitely many contracted curve valuations 
and they are all analytic.
In summary:
\begin{Lemma}\label{Lcontr}
  Fix $\nu\in \cV$. Then $f_*\nu$ defines a valuation on $R$ 
  with values in $\Rbar$, 
  which is identically $+\infty$ on $\fm$ iff $\nu\in\fC_f$,
  where $\fC_f$ denotes the 
  (finite) set of contracted curve valuations.
\end{Lemma}
\begin{Def}
  If $\nu\in\cV$ is not a contracted curve valuation,
  then $f_\bullet\nu\in\cV$ is the (unique)
  valuation equivalent to $f_*\nu$ and normalized by
  $(f_\bullet\nu)(\fm)=1$.
\end{Def}
\begin{Def}
  For $\nu\in\cV$ we define the 
  \emph{attraction rate of $f$ along $\nu$} to
  be the number $c(f,\nu)\=\nu(f^*\fm)\in[1, +\infty]$.
\end{Def}
The quantity $c(f)$ defined in the introduction
equals $c(f,\nu_\fm)$. Clearly
\begin{equation}\label{Edefcon}
  f_*\nu=c(f,\nu)\times f_\bullet\nu,
\end{equation}
for any $\nu\in\cV\setminus\fC_f$.

If $f$ and $g$ are dominant fixed point germs, then 
$(f\circ g)_*\nu=f_*g_*\nu$ for any $\nu\in\cV$, and
$(f\circ g)_\bullet\nu=f_\bullet g_\bullet\nu$ as long
as $\nu\not\in\fC_{f\circ g}$. 
This yields 
$c(f\circ g,\nu)=c(g,\nu)\times c(f,g_\bullet\nu)$. 
In particular, we get
\begin{equation}\label{Ecompconlocal}
  c(f^n)=\prod_{k=0}^{n-1}c(f,f_\bullet^k\nu_\fm).
\end{equation}
for all $n\ge 1$.
%
%
\subsection{Preservation of type}
Next we show that $f_\bullet$ preserves the type of
a valuation. 
\begin{Prop}\label{Pgen}
  The map $f_\bullet:\cV\setminus\fC_f\to\cV$ preserves 
  the set of divisorial (resp.\ irrational, infinitely singular, 
  or curve) valuations.
  Moreover, there exists an integer $N\ge1$ such that for any
  $\nu\in\cV$, the set $f_\bullet^{-1} \{ \nu\}$ has cardinality
  at most $N$.
\end{Prop}
We shall see in Section~\ref{Sintercurve} that $f_\bullet$ extends in
a natural way to $\fC_f$. When $\nu\in \fC_f$, however, 
$f_\bullet\nu$ is a divisorial valuation.
\begin{proof}
  Pick $\nu\in\cV\setminus\fC_f$. 
  Let us show that $\nu$ and $f_\bullet\nu$
  have the same type. We shall prove in 
  Section~\ref{Sintercurve} that if $\nu$
  is a curve valuation then so is $f_\bullet\nu$.
  Hence we may assume $\nu$ takes the value
  $\infty$ only on $0\in R$. 
  Since $f$ is dominant, $f_*\nu$ has the same property.
  As follows from the discussion in
  Section~\ref{sec-numerical}, we only have to 
  prove that $\nu$ and $f_*\nu$
  have the same rational rank and transcendence degree.

  Using Proposition~\ref{P-monomap}, we reduce to the case where
  $f$ is monomial, hence rational.
  Let $K$ be the field of rational functions on $\C^2$.  As $f$ is
  dominant, $f^*:K\to K$ is injective and $K$ is a \emph{finite}
  extension of $f^*K\= \{ \phi \circ f,\, \phi \in K \}$ of some
  degree $N\ge1$.

  Let us compare the two value groups $\Gamma_\nu$ and
  $\Gamma_{f_*\nu}$.  Clearly $\Gamma_{f_*\nu}\subset\Gamma_\nu$.
  Conversely consider $\gamma=\nu(\phi)\in\Gamma_\nu$.
  We may find $a_i\in f^*K$ such that $\sum_{i=0}^{N}a_i\phi^i=0$.  
  Thus there exist $i<j$ such that 
  $\nu(a_i\phi^i)=\nu(a_j\phi^j)$.
  In particular,
  $(j-i)\gamma=(j-i)\nu(\phi)=\nu(a_i)-\nu(a_j)\in\Gamma_{f_*\nu}$.
  We infer
  $\Gamma_\nu\otimes_\Z\Q\subset\Gamma_{f_*\nu}\otimes_\Z\Q$
  and hence $\Gamma_\nu\otimes_\Z\Q=\Gamma_{f_*\nu}\otimes_\Z\Q$.
  Thus the rational ranks of $\nu$ and $f_*\nu$ coincide. 
  The transcendence degrees must also coincide, as
  the extension of residue fields $k_{f_*\nu}\subset k_\nu$ 
  is finite.

  To prove the last assertion in the proposition, 
  we may again assume $f$ is rational.  
  Pick a finite set of distinct 
  valuations $\{\mu_i\}_{i\in I}\subset \cV$ for which $f_*\mu_i$ 
  are all proportional one to
  another. Suitable multiples of the $\mu_i$ thus define
  one and the same valuation $\mu$ on the
  field $f^* K$. But the field extension $K$ over $f^*K$ is finite of
  degree $N$, so $\mu$ admits at most $N$ extensions to $K$
  by~\cite[Theorem 19, p.55]{Z}.  Thus $\#I\le N$.  
\end{proof}
%
%
\subsection{Image of a divisorial valuation}\label{S-imagediv}
It is important for applications to understand geometrically what the
image of a valuation is. We first treat the case of a divisorial
valuation (cf.\ Section~\ref{sec-dual}).
\begin{Prop}\label{Pinterdiv} 
  Let $\nu$ be a divisorial valuation and set $\nu'=f_\bullet\nu$.
  Then there exist modifications 
  $\pi:X\to(\C^2,0)$, $\pi':X'\to (\C^2,0)$, and
  exceptional components $E\subset \pi^{-1}(0)$, 
  $E'\subset(\pi')^{-1}(0)$, 
  such that $\nu=\nu_E$, $\nu'=\nu_{E'}$ 
  and such that the map $f$ lifts to a holomorphic map 
  $\hf:X\to X'$ sending $E$ onto $E'$. 
  Moreover, 
  \begin{equation}\label{E103}
    c(f,\nu_E)=\frac{b_{E'}}{b_E}k.
  \end{equation}
  Here $b_E$ and $b_{E'}$ are the generic multiplicities of
  $\nu_E$ and $\nu_{E'}$, respectively and $k\ge1$ is the largest
  integer such that $\hf^*E'\ge kE$ as divisors.
\end{Prop}
\begin{proof}
  Proposition~\ref{Pgen} implies that $\nu'$ is divisorial.
  Choose modifications 
  $\pi:X\to(\C^2,0)$, $\pi':X'\to(\C^2,0)$, such that the center 
  $E$ of $\nu$ in $X$ and $E'$ of
  $\nu'$ in $X'$ are both one dimensional. By blowing up the
  source space further, 
  we may assume that the lift $\hf$ of $f$ is holomorphic. 
  By construction, a holomorphic function locally defined 
  at a point $p'\in E'$ vanishes along $E'$ iff its pull back by $\hf$
  vanishes along $E$. So $f$ maps $E$ onto $E'$. 
  Pick $k\ge1$ maximal with $\hf^*E'\ge kE$. Since
  $\pi^*f^*\fm$ and $\hf^*(\pi')^*\fm$ vanish along $E$ 
  to order $b_Ec(f,\nu_E)$ and $b_{E'}k$, respectively,
  we get~\eqref{E103}.
\end{proof}
%
%
\subsection{Image of a curve valuation}\label{Sintercurve}
Next we give an analogue of Proposition~\ref{Pinterdiv}
for analytic curve valuations. We distinguish between curves that
are contracted and those that are not.
\begin{Prop}\label{Pintercurve} 
  Let $C$ be an analytic irreducible curve such that 
  $f(C)\ne 0$ \ie $\nu_C\notin\fC_f$. 
  Then $C'\=f(C)$ is an analytic 
  irreducible curve and $f_\bullet\nu_C=\nu_{C'}$. Further,
  \begin{equation}\label{E104}
    c(f,\nu_C)=\frac{m(C')}{m(C)}e(f,C),
  \end{equation}
  where $e(f,C)\in \N^*$ denotes the topological degree of the
  restriction $f:C\to C'$.
\end{Prop}
As the proof shows, the proposition is also valid for
formal irreducible curves.
\begin{Prop}\label{Pintercurve2} 
  Suppose $C$ is an irreducible curve germ such that $f(C)=0$,
  \ie $\nu_C\in\fC_f$. Then $c(f,\nu_C)=\infty$. Further,
  the limit of $f_\bullet\nu$ as $\nu$ increases to $\nu_C$ 
  exists, and is a divisorial valuation that 
  we denote by $f_\bullet \nu_C$. It can
  be interpreted geometrically as follows. There exist 
  modifications $\pi:X\to(\C^2,0)$ and $\pi':X'\to(\C^2,0)$,
  such that $f$ lifts to a holomorphic map $\hf:X\to X'$ sending $C$
  to a curve germ included in an exceptional component $E'$ of
  $\pi'^{-1}(0)$, for which $f_\bullet\nu_C=\nu_{E'}$.
\end{Prop}
\begin{proof}[Proof of Proposition~\ref{Pintercurve}]
  Write $C=\phi^{-1}(0)$ and $C'=(\phi')^{-1}(0)$ for 
  $\phi,\phi'\in\fm$ irreducible.
  Then $(f_*\nu_C)(\phi')=\nu_C(\phi'\circ f)=\infty$ 
  as $\phi$ divides $\phi'\circ f$. 
  Hence $f_\bullet\nu_C=\nu_{C'}$.

  Pick parameterizations $h:\Delta\to C$ 
  and $h':\Delta\to C'$.
  The composition $(h')^{-1}\circ f\circ h$ has
  topological degree $e\=e(f,C)$ at the origin.
  We may assume $f\circ h(t)=h'(t^e)$.  
  Now pick 
  $x\in\fm$ of multiplicity $1$, and generic so that 
  $\nu_C(x)=\nu_{C'}(x)=1$. Then
  \begin{multline*}
    c(f,\nu_C)
    =(f_*\nu_C)(x)
    =\nu_C(x\circ f)
    =m(C)^{-1}\div_t(x\circ f\circ h(t))=\\
    =m(C)^{-1}\div_t(x\circ h'(t^e))
    =e\times m(C)^{-1}\div_t(x\circ h'(t))
    =e\times\frac{m(C')}{m(C)}.
  \end{multline*}
  This concludes the proof of the proposition.
\end{proof}
\begin{proof}[Proof of Proposition~\ref{Pintercurve2}]
  We know from Lemma~\ref{Lcontr} that $c(f,\nu_C)=\infty$.
  Let $\div_C$ be the \emph{divisorial} valuation on $R$
  associated to $C$. It is not centered at the origin, but 
  since $C$ is contracted, $f_*\div_C$ is equivalent to a 
  divisorial valuation $\nu'\in\cV$.

  Pick modifications $\pi:X\to(\C^2,0)$ and $\pi':X'\to(\C^2,0)$
  such that $\pi^{-1}(C)$ has simple normal crossings,
  the center of $\nu'$ on $X$ is an exceptional 
  component $E'\subset(\pi')^{-1}(0)$ and such that
  the lift $\hf:X \to X'$ of $f$ is holomorphic 
  and locally monomial at every point on $\pi^{-1}(0)$.

  Write $\tC$ for the strict transform of $C$: it intersects
  the exceptional divisor at a single point $p$ on 
  a unique exceptional component $E$. 
  Set $p'=\hf(p)$. 
  As in Proposition~\ref{Pinterdiv}, $\hf$ maps
  $\tC$ to the germ of $E'$ at $p'$. In suitable 
  coordinates $(z,w)$ at $p$ and $(z',w')$ at $p'$ 
  we have $\hf(z,w)=(z^a,z^cw^d)$ with $a,d\ge1$, $b\ge0$.
  Here $\tC=\{w=0\}$, $E=\{z=0\}$ and $E'=\{w'=0\}$.
  The ideals $\pi^*\fm$ and $(\pi')^*\fm$ are principal at 
  $p$ and $p'$, generated by $z^m$ 
  and $z'^kw'^l$, respectively, where $m=m(C)\ge1$,
  $k\ge0$, $l\ge1$.

  Let $\mu_{s,t}$ (resp. $\mu'_{s,t}$) be the monomial valuation with
  weights $(s,t)$ on the coordinates $(z,w)$ 
  (resp.\ $(z',w')$).
  The valuation $\nu_t\=\pi_*\mu_{1/m,t}$ is an element of $\cV$ 
  and increases to $\nu_C$ as $t\to\infty$.
  Clearly $f_*\nu_t=\pi'_*\mu'_{a/m,c/m+dt}$.
  As $t\to\infty$, $f_\bullet \nu_t$ converges to 
  $\pi'_*\mu'_{0,d/l}=\nu'$. The proof is complete.
\end{proof}
%
%
%
%
%
\section{Action on the valuative tree}\label{sec-acttree}
Having defined the map $f_\bullet$ pointwise on $\cV\setminus\fC_f$,
we now investigate how it interacts with the structure on $\cV$.
First we show that it preserves the tree structure in a strong sense
(Theorem~\ref{Ttreemap}).  We then show how $f_\bullet$ can be used to
analyze the lift of $f$ to a (meromorphic) map between surfaces
obtained by blowing up points over the origin
(Propositions~\ref{Pintergeo} and~\ref{Ptangentinter}).
\begin{Thm}\label{Ttreemap}
  The map $f_\bullet:\cV\setminus\fC_f\to\cV$ extends uniquely to a
  map $f_\bullet:\cV\to\cV$ which is regular in the following sense:
  for any $\nu\in\cV$, we can decompose the segment $[\nu_\fm,\nu]$
  into finitely many subsegments $I_j=[\nu_j,\nu_{j+1}]$, 
  $0\le j\le k$ with $\nu_0=\nu_\fm$, $\nu_{k+1}=\nu$ and $\nu_i$
  divisorial, such that $f_\bullet$ is a monotone homeomorphism of
  $I_j$ onto its image.
  
  For $t\in [1,\a(\nu)]$ denote by $\nu_t$ the unique
  valuation dominated by $\nu$ and of skewness~$t$.
  Then on any segment
  $[\a(\nu_i), \a(\nu_{i+1})]$ we have
  $\a(f_\bullet\nu_t)=\frac{a+bt}{c+dt}$ for nonnegative integers
  $a,b,c,d$ with $ad\ne bc$.
\end{Thm}
\begin{Prop}\label{Pintergeo}
  Let $\pi:X\to(\C^2,0)$, $\pi':X'\to(\C^2,0)$ be 
  modifications, and let $\hf:X\dashrightarrow X'$ 
  be the lift of $f$.
  For an infinitely near point $p\in\pi^{-1}(0)\subset X$, 
  let $U(p)\subset\cV$ be the open set of valuations whose center 
  on $X$ is $p$ (see Section~\ref{S21}).
  Then $\hf$ is holomorphic at $p$ iff $f_\bullet U(p)$ 
  does not contain any divisorial valuation associated
  to an exceptional component of $\pi'$. When $\hf$ is holomorphic at
  $p$, the point $p'=\hf(p)\in X'$ 
  is characterized by $f_\bullet U(p)\subset U(p')$.
\end{Prop}
It follows from Theorem~\ref{Ttreemap} that $f$ induces a
\emph{tangent map} $Df_\bullet$ between tree tangent spaces.
Indeed, if $\nu\in\cV$ is any valuation and $\vv$ is a tangent
vector at $\nu$, then $\vv$ is represented by a valuation 
$\mu$ such that $f_\bullet$ is a homeomorphism of $[\nu,\mu]$
onto $[f_\bullet\nu,f_\bullet\mu]$. We let $Df_\bullet\vv$ be the
tangent vector at $f_\bullet\nu$ represented by $f_\bullet\mu$;
it clearly does not depend on the choice of $\mu$.

At a divisorial valuation, the tangent map
has the following interpretation.
\begin{Prop}\label{Ptangentinter}
  Consider a divisorial valuation $\nu$ and set $\nu'=f_\bullet\nu$. 
  Let $\pi:X\to(\C^2,0)$, $\pi':X'\to (\C^2,0)$ be 
  modifications such that 
  $\nu$ and $\nu'$ are associated to exceptional components 
  $E\subset\pi^{-1}(0)$,
  and $E'\subset(\pi')^{-1}(0)$, respectively. 
  The lift $\hf:X\dashrightarrow X'$ sends $E$ surjectively 
  onto $E'$. Let $f_E$ be its restriction to $E$. 

  For any point $p\in E$
  let $\vv_p$ be the tangent vector at $\nu$ defined by
  $p$ as in Section~\ref{S21}.
  Use the same notation for points on $E'$ 
  (tangent vectors are now at $\nu'$).
  Then $Df_\bullet(\vv_p)=\vv_{p'}$ where $p'=f_E(p)$.
\end{Prop}
When $f_\bullet\nu=\nu$ and $\pi=\pi'$, 
$f_E$ is a rational selfmap of $E$, hence admits a 
noncritical fixed point. 
This fact will be used in the proof of Theorem~C.
%
%
\subsection{The critical tree}
As a step towards Theorem~\ref{Ttreemap} 
we first investigate how $f_\bullet$ interacts with the partial
ordering on $\cV$.
Note that $f_*$ is order preserving:
if $\nu\ge\nu'$ in the sense that 
$\nu(\phi)\ge\nu'(\phi)$ for all $\phi$, 
then $f_*\nu\ge f_*\nu'$. 
The situation concerning $f_\bullet$ is more
complicated as the function $\nu\mapsto c(f,\nu)$ is also
increasing; thus $f_\bullet \nu $ is the quotient of
two increasing functions, and in general not increasing.
\begin{Prop}\label{Pcrittree}
  The subset $\cT_f$ of $\cV$ where $c(f,\cdot)$ is not 
  locally constant is a finite closed subtree of $\cV$.  
  Its ends (\ie its maximal elements) are exactly the maximal 
  elements in the finite set $\cE_f$ consisting of 
  divisorial valuations $\nu$ with $f_\bullet\nu=\nu_\fm$ 
  and of contracted curve valuations.
\end{Prop}
In particular, all contracted curve valuations are
ends in $\cT_f$.
We call $\cT_f$ the \emph{critical tree} of $f$.
Note that $f_\bullet$ is order preserving on $\cV\setminus\cT_f$.
\begin{Example}
  Let $f(x,y)=(x^2y,xy^3)$. The critical tree consists of all
  monomial valuations, and its ends are the two contracted curve
  valuations $\nu_x$ and $\nu_y$.  The monomial valuation $\nu_{x,2}$
  with weight $2$ on $x$ and $1$ on $y$ is the only preimage of
  $\nu_\fm$, so $\cE_f=\{\nu_x,\nu_y,\nu_{x,2}\}$.  Note that
  $f_\bullet$ is not surjective onto $\cV$, 
  as $\nu_x\notin f_\bullet(\cV)$.
\end{Example}
\begin{proof}[Proof of Proposition~\ref{Pcrittree}]
  The function $c(f,\cdot)$ 
 is a tree potential
  (see Section~\ref{sec-treepot1}). Indeed, 
  if we write $f=(g,h)$ in coordinates, then 
  $c(f,\nu)=\min\{\nu(g),\nu(h)\}$.
  It is then a general fact that the locus where $c(f,\cdot)$ is
  not locally constant is a subtree of $\cV$. In our case, this
  locus $\cT_f$ is a finite subtree as it is contained in the
  subtree $\cS$ consisting of valuations dominated by curve 
  valuations $\nu_\phi$ such that $\phi$ divides $gh$.
  
  Let $\nu$ be an end of $\cT_f$. We shall prove that it is a
  maximal element in $\cE_f$.
  Note that $\nu$ is quasimonomial or a curve valuation, as
  $\cS$ does not contain infinitely singular valuations.  First
  suppose that $\nu=\nu_\phi$ is a curve valuation.  The tree
  potential $c(f,\mu)=\mu(f^*\fm)$ is piecewise affine on any segment
  with \emph{integer} slope in terms of skewness.
  Thus $c(f,\mu)\to\infty$ as $\mu$
  increases to $\nu_\phi$, so $\phi$ divides any
  element in $f^*\fm$. This means $\nu_\phi$ is a
  contracted curve valuation, hence a maximal element of $\cE_f$.
  
  Now assume the end $\nu$ of $\cT_f$ is quasimonomial. 
  Let $\cL$ denote the set of functions 
  $L_{a,b}=ax+by$, $(a,b)\ne(0,0)$. 
  Suppose there exists $L\in\cL$ such that $\nu(f^*L)>c(f,\nu)$.
  Then the same inequality holds in a
  neighborhood $U$ of $\nu$.  
  Note that $c(f,\mu)=\min\{\mu(f^*L),\mu(f^*L')\}$ for any
  fixed $L'\neq L$ and any $\mu\in\cV$.  
  Thus $c(f,\mu)=\mu(f^*L')$ for any $\mu\in U$.  
  Now a tree potential
  of the type $\mu\to\mu(\phi)$, $\phi\in R$ cannot attain a local
  maximum at a point in $\cVqm$ except if it is locally constant
  there.  Hence there exists $\mu>\nu$ for which $c(f,\mu)>c(f,\nu)$, 
  a contradiction. 
  Thus $f_*\nu(L)=\nu(f^*L)=c(f,\nu)$ is a constant 
  independent of $L\in\cL$. This implies that 
  $f_*\nu$ is proportional to $\nu_\fm$, so $\nu\in\cE_f$.
  Now the set $\{\mu>\nu\}$ does not meet $\cE_f$ as
  $c(f,\cdot)$ is constant and $f_\bullet$ order-preserving on it,
  so $\nu$ is a maximal element of $\cE_f$.

  Conversely, let $\nu$ be a maximal element of $\cE_f$.
  If $\nu$ is a contracted curve valuation, then 
  $c(f,\mu)\to\infty$ as $\mu$ increases to $\nu$, so
  $\nu$ is an end in $\cT_f$.
  If $\nu$ is divisorial and $f_\bullet\nu=\nu_\fm$,
  then $f_\bullet$ cannot be order preserving near $\nu$, so
  $\nu\in\cT_f$. 
  Since $\nu$ is maximal in $\cE_f$, as are all ends of $\cT_f$,
  $\nu$ must be an end of $\cT_f$.
\end{proof}
%
%
\subsection{Regularity}
We are now ready to prove Theorem~\ref{Ttreemap}

Note first that the image of a monomial valuation $ \nu_{s,t}$ with
$\nu_{s,t}(z) =s$, $\nu_{s,t}(w)=t$ 
by a monomial map $f(z,w)=(z^aw^b,z^cw^d)$ 
equals $\nu_{as+bt,cs+dt}$. In particular, when $f$
is dominant the induced map $f_*$ on the set of monomial valuations is
injective. More precisely, $f_* \nu_{s,t}$ and $f_*\nu_{s',t'}$ are
proportional iff the two vectors $(s,t)$, $(s',t')$ are.

Now consider a general fixed point germ $f$. We first treat the case 
of an analytic curve valuation $\nu = \nu_C$ 
and refer to the end of proof for the other cases. 
By Proposition~\ref{P-monomap}, one can find modifications
$\pi:X\to(\C^2,0)$, $\pi': X' \to (\C^2,0)$ such that
$\pi^{-1}(C)$ has simple normal crossings and
the lift $\hf:X \to X'$ of $f$ is holomorphic and locally monomial at
any point on $\pi^{-1}(0)$.  

By Section~\ref{sec-dual}, the dual graph $\Gamma$ of the
(reducible) curve $\pi^{-1}(C)$ embeds in $\cV$. 
It contains
$\nu_\fm$ and $\nu_C$, hence the full segment 
$I=[\nu_\fm,\nu_C]$. 
Let $E_0,\dots, E_n$ be the vertices of $\Gamma$
(\ie irreducible components of $\pi^{-1}(C)$) whose associated
valuations belong to $I$. We order them by 
$\nu_\fm=\nu_{E_0}<\nu_{E_1}\dots<\nu_{E_n}=\nu_C$.  
The segment $[\nu_{E_i},\nu_{E_{i+1}}]$ consists of valuations in $\cV$ that are
monomial in suitable coordinates at 
$p_i\=E_i\cap E_{i+1}$: see Section~\ref{S22}.
As $\hat{f}$ is monomial in these coordinates, 
the restriction of $f_\bullet$ to
$[\nu_{E_i},\nu_{E_{i+1}}]$ is injective. 
This proves the first part of the theorem.

To prove the second part, we need to control skewness.
Pick $\nu_0,\nu_1\in\cV$ such that $\nu_0<\nu_1$, 
$f_\bullet$ is monotone on $[\nu_0,\nu_1]$, and
$\nu_1$ is not infinitely singular.
Then $f_\bullet[\nu_0,\nu_1]$ is a segment containing no
infinitely singular valuation so we can pick a curve valuation 
$\nu_\psi\in\cV$ such that $\nu_\psi\ge f_\bullet\mu$
for all $\mu\in[\nu_0,\nu_1]$.
Then 
\begin{equation}\label{E201}
  \a(f_\bullet\mu)
  =\frac{f_\bullet\mu(\psi)}{m(\psi)}
  =\frac{f_*\mu(\psi)}{c(f,\mu)\cdot m(\psi)}
  =\frac{\mu(f^*\psi)}{c(f,\mu)\cdot m(\psi)}.
\end{equation}
Now $\mu\mapsto\mu(f^*\psi)=\sum_i\a(\mu\wedge\nu_{\psi_i})$ 
is piecewise linear in $\a(\mu)$ 
with nonnegative integer coefficients. 
(Here $f^*\psi=\prod\psi_i$ with $\psi_i$ irreducible.)
The same is true for $c(f,\mu)=\min\{\mu(f^*x),\mu(f^*y)\}$. 
Thus $\a(f_\bullet\mu)$ is piecewise a M{\"o}bius function 
of $\a(\mu)$ with nonnegative integer coefficients, as desired.

\smallskip
Finally we consider the statement in Theorem~\ref{Ttreemap} when $\nu$
is not an analytic curve valuation. When $\nu$ is quasimonomial it is
dominated by an analytic curve valuation, so the previous argument
applies. Otherwise it is either a formal curve or infinitely singular
valuation.  But then $\nu$ does not belong to the critical tree $\cT_f$,
so we may find a divisorial valuation $\nu_1\le \nu$ such that
$c(f,\cdot)$ is constant on $[\nu_1,\nu]$, equal to $c\=c(f,\nu)$. In
particular, $f_\bullet$ is increasing on $[\nu_1,\nu]$, and gives a
bijection onto $[f_\bullet\nu_1,f_\bullet\nu]$.  We claim that
$\a(f_\bullet \mu)$ is in fact an \emph{affine} function of $\a(\mu)$
on $[\nu_1,\nu]$, at least if $\nu_1$ is chosen large enough.
To prove this, pick $\nu_2 \in [\nu_1,\nu]$. Fix $\psi\in\fm$
irreducible, such that $\nu\wedge\nu_\psi>\nu_2$. As $c(f,\cdot)$ is
constant,~\eqref{E201} shows that $\a(f_\bullet\mu)$ is piecewise
linear with rational coefficients on each segment $[\nu_1,\nu_2]$.
But the number of points in $[\nu_1,\nu_2]$ at which
$\a(f_\bullet\mu)$ is not smooth is bounded from above by a constant
depending only on $f$ (and not on $\psi$) by Proposition~\ref{Pgen}.  
We can hence decompose $[\nu_1,\nu]$ into finitely many segments, 
on which $\a(f_\bullet\mu)$ is piecewise linear in $\a(\mu)$.

This completes the proof of Theorem~\ref{Ttreemap}.
%
%
\subsection{Geometric interpretation of indeterminacy points}
We now prove Proposition~\ref{Pintergeo}.
Suppose first that $p$ is not an indeterminacy point of $\hf$
and write $p'=\hf(p)$. Pick a valuation $\nu\in U(p)$. 
This means $\nu=\pi_*\mu$,
for some centered valuation $\mu$ on $\cO_p$.
As $\hf$ is holomorphic at $p$ we have 
$\hf^*\cO_{p'}\subset\cO_p$ and $\hf^*\fm_{p'}\subset\fm_p$.
By duality, $\mu':=\hf_*\mu$ is a centered valuation on $\cO_{p'}$.
Now $\pi'_*\mu'$ is proportional to $f_\bullet\nu$, so
$f_\bullet\nu\in U(p')$.
Thus $f_\bullet U(p)\subset U(p')$.
On the other hand, the center of a divisorial valuation $\nu'$
associated to an exceptional component $E'\subset(\pi')^{-1}(0)$
is equal to $E'\ne p'$, so $\nu'\not\in U(p')$.

Conversely, let $p$ be an indeterminacy point of $\hf$.
There exists a modification $\varpi:Y\to X$
such that $\hf$ lifts to a holomorphic map, $g:Y\to X'$
and the image of $\varpi^{-1}(p)$ under $g$ 
contains a curve $E'\subset X'$.  
We may thus find an exceptional component 
$E\subset\varpi^{-1}(p)$ that is
sent by $g$ surjectively onto $E'$. 
Hence the center of $\nu_E$ on $X$ is $p$, so
$\nu_E\in U(p)$. By construction, $f_\bullet\nu_E=\nu_{E'}$. 
The proof is complete.
%
%
\subsection{The tangent map}
Finally we prove Proposition~\ref{Ptangentinter}. 
The lift $\hf:X\dashrightarrow X'$ need not be holomorphic, 
but we may find a modification
$\varpi:Y\to X$ such that $\hf$ lifts to a holomorphic
map $g:Y\to X'$.
Pick a point $q\in\varpi^{-1}(p)$ on the
strict transform of $E$.

If a valuation in $\cV$ has center $q$ on $Y$, then it
has center $p$ on $X$. Hence $U(q)\subset U(p)$. Moreover,
$\nu$ lies in the closure of $U(q)$
and all valuations in $U(p)$ represent the same tangent 
vector $\vv_p$ at $\nu$. 
Similarly, $\nu'\in\overline{U(p')}$
and all valuations in $U(p')$ represent the same
tangent vector $\vv_{p'}$ at $\nu'$.

Proposition~\ref{Pintergeo} shows that $f_\bullet U(q)\subset U(p')$.
By the definition of the tangent map $Df_\bullet$, we conclude that
$Df_\bullet\vv_p=\vv_{p'}$.
This concludes the proof.
%
%
%
%
\section{Eigenvaluations}\label{sec-fixed}
The goal of this section is to show that the induced map
$f_\bullet:\cV\to\cV$ admits a fixed 
point and use the latter to analyze the attraction rate 
of $f$.
\begin{Def}
  The \emph{attraction rate} of $f$ is
  $c(f):=c(f,\nu_\fm)=\nu(f^*\fm)$.
  The \emph{asymptotic attraction rate} is
  the limit $c_\infty\=\lim_{n\to \infty} c(f^n)^{1/n}$.
\end{Def}
The limit defining $c_\infty$ exists as the sequence
$c(f^n)$ is supermultiplicative. 
The main result of this section is
\begin{Thm}\label{Teigen}
  Any dominant fixed point germ $f:(\C^2,0)\to(\C^2,0)$ 
  with $c_\infty(f)>1$
  admits a valuation $\nu_\star\in\cV$ with 
  $f_\bullet\nu_\star=\nu_\star$ and $c(\nu_\star)=c_\infty$.
  The valuation $\nu_\star$ may be chosen quasimonomial,
  infinitely singular, or as a 
  noncontracted, analytic curve valuation.
  When $\nu_\star$ is not quasimonomial, 
  there exists $\nu_0<\nu_\star$ such
  that $c(f,\nu_0)=c_\infty$, $f_\bullet$ is order-preserving on
  $\{\nu\ge\nu_0\}$ and $f_\bullet\nu>\nu$ for
  $\nu\in[\nu_0,\nu_\star[\,$.
  
  Finally, we may find $0<\delta\le 1$ such that $\delta\cdot
  c_\infty^n\le c(f^n)\le c_\infty^n$ for all $n\ge 1$.
\end{Thm}
\begin{Def}\label{Deigen}
  A valuation $\nu_\star$ as above
  is called an \emph{eigenvaluation} for $f$.
  Note that $\nu_\star$ may be viewed as an eigenvector for
  $f_*$ with eigenvalue $c_\infty$.
\end{Def}

We do not claim that the eigenvaluation is unique.
The proof of Theorem~\ref{Teigen} will be given 
in Section~\ref{S24} after we have discussed the main ingredient:
a purely tree-theoretic fixed point theorem.
%
%
\subsection{Fixed point theorem on trees}\label{S23}
Let us recall some general definitions on trees taken
from~\cite[Ch.3]{valtree}.  A (rooted, nonmetric) \emph{complete tree}
is a poset $(\cT,\le)$ with a unique element minimal 
$\tau_0$ (the root) such
that for any $\tau,\tau'\in \cT$ with $\tau<\tau'$ there exists an
increasing bijection from $[0,1]$ (endowed with the standard order
relation) onto the \emph{segment} $[\tau,\tau']:=\{\sigma\ ;\ 
\tau\le\sigma\le\tau'\}$ and such that any totally ordered subset of
$\cT$ has a majorant in $\cT$.  To any segment
$I=[\tau,\tau']\subset\cT$ with $\tau \le \tau'$ is associated a
\emph{retraction} $\pi_I:\cT\to I$: $\pi_I(\sigma)$ is the maximal
element in $I\cap[\tau_0,\sigma]$ if this set is nonempty, and
$\pi_I(\sigma)=\tau$ otherwise.  The map $\pi_I$ is continuous and
surjective onto $I$, and $\pi_I|_I=\id$. We endow $\cT$ with the
weakest topology in which $\pi_I$ is continuous for every segment $I$.
It is generated by sets $U(\vv)$ over tangent vectors $\vv$ as in
Section~\ref{sec-tangent}.
\begin{Def}
  A map $\cF:\cT\to\cT'$ between complete trees
  $\cT$, $\cT'$ is a \emph{tree map} if it is
  continuous on branches, \ie if the restriction
  $\cF|_I:I\to\cT'$ is continuous 
  for every segment $I$ in $\cT$.

  A tree map $\cF$ is \emph{regular} if any segment $I$ in
  $\cT$ can be decomposed into finitely many segments 
  on each of which $\cF$ is a homeomorphism onto its image.
\end{Def}
There is a natural way in which selfmaps of trees induce selfmaps of
segments. 
Let $\cF:\cT\to\cT$ be a tree map and let $I$ 
be a segment in $\cT$.
Set $\cF_I:=\pi_I\circ\cF|_I$, where $\pi_I:\cT\to I$ is the
retraction defined above. Then $\cF_I$
is a continuous selfmap of $I$. 
Notice that $\cF\ne\cF_I$ on $I$ in general:
the set $\{\tau\in I\ ;\ \cF(\tau)\ne\cF_I(\tau)\}$ 
is an open subset of $I$ on which $\cF_I$ is locally constant.

We define an \emph{end} $\tau\in\cT$, \ie a maximal element 
under $\le$, to be \emph{weakly attracting} if
there exists a segment $I=[\tau',\tau]$ such that 
$\tau$ is a globally attracting fixed point 
for the induced map $\cF_I$
in the sense that $\cF_I(\sigma)>\sigma$ for every
$\sigma\in[\tau',\tau[$. 
The end is \emph{strongly attracting} when in addition
the segment $I$ can be chosen $\cF$-invariant.

We are now ready to state the following fixed point 
property for tree maps.
\begin{Thm}\label{MT1}
  Let $\cT$ be a complete tree and $\cF:\cT\to\cT$ a tree map 
  (resp.\ a regular tree map). 
  Then one of the following two statements hold:
  \begin{itemize}
  \item $\cF$ admits a fixed point $\tau$ which is not an end;
  \item $\cF$ admits a weakly (resp.\ strongly) attracting end.
\end{itemize}
\end{Thm}
\begin{Remark}
  Theorem~\ref{MT1} was proved independently by
  Rivera~\cite[Section~8.3]{Rivera2} in his work on $p$-adic dynamics.
  It may be viewed as a ``non-invertible'' analog of the description
  of isometries on metric $\R$-trees (see~\cite{MorganShalen,Paulin}).
\end{Remark}
\begin{proof}[Proof of Theorem~\ref{MT1}]
  Suppose first $\cF$ is a tree map. 
  We may assume that $\cF\tau_0\ne\tau_0$, or else there is nothing
  to prove. Pick any end $\tau'_0$ such that
  $\tau_0<\cF\tau_0\le\tau'_0$ and let
  $I_0=[\tau_0,\tau'_0]$. We consider the induced selfmap
  $\cF_0=\cF_{I_0}$ of $I_0$ defined above.

  Since $I_0$ is homeomorphic to a compact interval, $\cF_0$ admits a
  fixed point $\tau_1\in I_0$. Further, $\tau_1>\tau_0$ as
  $\cF_0\tau_0>\tau_0$.  Notice that we may assume that $\tau_1$ is
  not an end unless $\cF_0\tau>\tau$ for every
  $\tau\in[\tau_0,\tau_0'[$, in which case $\tau_1=\tau_0'$ is a
  weakly attracting fixed point for $\cF$.
  
  If $\tau_1\in]\tau_0,\tau'_0[$, then $\tau_1$ is either a
  fixed point of $\cF$, in which case we are done, or satisfies
  $\cF\tau_1>\tau_1$. In the latter case pick an end $\tau'_1$
  with $\tau_1<\cF\tau_1\le\tau'_1$ and let
  $I_1=[\tau_1,\tau'_1]$.  We write $\cF_1:=\cF_{I_1}$ for
  the induced selfmap of $I_1$.
  The map $\cF_1$, being a continuous selfmap of an compact
  interval, has a
  fixed point $\tau_2\in I_1$.  Since $\cF\tau_1>\tau_1$ we
  have $\cF_1\tau_1>\tau_1$, which implies $\tau_2>\tau_1$.  
  We may assume $\tau_1<\tau_2<\tau'_1$ unless $\tau'_1$ is a weakly
  attracting fixed point for $\cF$.
  If $\tau_2\ne\tau'_1$ is an end, then we assume
  $\cF\tau_2>\tau_2$, or else we are done. Then the inductive
  procedure continues.
  
  If this procedure stops, we obtain a fixed point for $\cF$: 
  a non-end or a weakly attracting end. 
  Suppose the procedure continues indefinitely.  
  We get a sequence of points
  $\tau_0<\tau_1<\tau_2<\dots$ such that $\cF\tau_n>\tau_n$
  for every $n$. Let $\tau_\infty=\sup\tau_n$. We claim that
  $\cF\tau_\infty=\tau_\infty$.
  To see this, consider $\tau'_n$, the inductively constructed ends.
  Since $\cF_n\tau_n=\tau_n$, the segments $[\tau_n,\tau_\infty]$
  and $[\tau_n,\tau'_n]$ intersect only at $\tau_n$.  In
  particular, $\cF\tau_n\not\ge\tau_\infty$. By continuity this
  implies that $\cF\tau_\infty\not>\tau_\infty$. On the other hand,
  $\cF\tau_n>\tau_n$, so continuity also gives
  $\cF\tau_\infty\ge\tau_\infty$.  Thus $\tau_\infty$ is a fixed
  point for $\cF$.

  If $\tau_\infty$ is an end, then set $I=[\tau_0,\tau_\infty]$
  and $I^*=[\tau_0,\tau_\infty[$.
  We claim that either $\cF$ has a fixed point $\tau\in I^*$, or
  $\cF_I\tau>\tau$ for all $\tau\in I^*$.  In either
  case the proof is complete.
  
  Therefore assume that $\cF$ has no fixed point in $I^*$ and that
  $\cF_I\tau\le\tau$ for some $\tau\in I^*$ Since
  $\cF_I\tau_n>\tau_n$ and $\tau_n\to\infty$ there then
  exists a fixed point $\tau'\in I$ for $\cF_I$.  By
  assumption, $\cF\tau'\ne\tau'$, so we must have $\cF\tau'>\tau'$
  and $\cF\tau'\notin I$.  Thus $\cF_I$ is locally constant at
  $\tau'$.  Pick $\tau''\in I$ maximal such that
  $\cF_I\tau<\tau$ on $]\tau',\tau''[$.  We have
  $\tau''<\tau_\infty$ since $\cF_I\tau_n>\tau_n$ for every
  $n$. Clearly $\cF_I\tau''=\tau''$ and we claim that in fact
  $\cF\tau''=\tau''$. Indeed, if not, then $\cF_I$ is locally
  constant at $\tau''$, which implies that $\cF_I\tau>\tau$
  for $\tau\approx\tau''$, $\tau<\tau''$.  This is a
  contradiction, completing the proof when $\cF$ is a tree map.

  Finally suppose  $\cF$ is regular, and $\tau$ is a weakly attracting
  end.  Consider a segment $I=[\tau_0,\tau[$, such that
  $\cF_I\tau'>\tau'$ for all $\tau'\in I$. Since $\cF$ is
  regular, we may in fact assume that $\cF$ is a homeomorphism on
  $I$.  This implies that $\cF_I=\cF$ on $I$.
\end{proof}
%
%
\subsection{Proof of Theorem~\ref{Teigen}}\label{S24}
We know from~\cite[Theorem~3.14]{valtree} that $(\cV,\le)$ is a
complete tree.  Theorem~\ref{Ttreemap} implies that the induced map
$f_\bullet:\cV\to\cV$ is a regular tree map.
Thus Theorem~\ref{MT1} applies, and provides a fixed point
$\nu_\star\in\cV$.

If $\nu_\star$ is not an end, then it is quasimonomial, with finite
skewness $\a_\star:=\a(\nu_\star)$.  This implies
$\nu_\fm(\phi)\le\nu_\star(\phi)\le\a_\star\nu_\fm(\phi)$ for any
$\phi\in R$. But $c(f^n,\nu)=\nu(f^{n*}\fm)$, so we get
$c(f^n,\nu_\fm)\le c(f^n,\nu_\star)\le\a_\star c(f^n,\nu_\fm)$.  As
$f_\bullet\nu_\star=\nu_\star$, we have
$c(f^n,\nu_\star)=c(f,\nu_\star)^n$.  This proves that $c(f^n)\le
c(f,\nu_\star)^n\le \a_\star c(f^n)$.  In particular,
$c_\infty=c(f,\nu_\star)$ and $c(f^n)/c_\infty^n\ge \delta$ for all
$n$ with $\delta=\a_\star^{-1}$.

Now suppose $\nu_\star$ is a strongly attracting end.  By definition
this implies that $f_\bullet\nu>\nu$ on $[\nu_0,\nu_\star[$ for some
$\nu_0<\nu_\star$.  We claim that $\nu_\star$ cannot belong to
$\cT_f$, the critical tree of $f$. Indeed, by
Proposition~\ref{Pcrittree} the ends of $\cT_f$ are either
quasimonomial or contracted curve valuations. But the latter are
mapped by $f_\bullet$ to divisorial valuations by
Proposition~\ref{Pintercurve}, hence cannot be fixed points. Thus
$\nu_\star\notin\cT_f$. Since $\cT_f$ is closed, we may assume
$\nu_0\not\in\cT_f$ so that $c(f,\cdot)$ is constant and $f_\bullet$
order-preserving on $\{\nu\ge\nu_0\}$.  This implies
$c(f^n,\nu_0)=\prod_0^{n-1}c(f,f_\bullet^k\nu)=c_\star^n$, where
$c_\star=c(f,\nu_\star)$.  The skewness $\a_0$ of $\nu_0$ is
finite and $c(f^n)\le c(f^n,\nu_0)\le\a_0c(f^n)$.  As before
we conclude that $c_\star=c_\infty$ and that $c(f^n)/c_\infty^n\ge\a_0^{-1}$.

If $\nu_\star$ is a strongly attracting end, then it
is infinitely singular or a curve valuation. 
Suppose it is a curve valuation. 
We saw above that it cannot be contracted.
To prove that it is analytic, 
we rely on the following useful formula.
\begin{Lemma}\label{LJaclocal}
  For any $\nu\in\cV\setminus\fC_f$ we have 
  \begin{equation}\label{eJaclocal}
    c(f,\nu)A(f_\bullet\nu)=A(\nu)+\nu(Jf).
  \end{equation}
\end{Lemma}
Here $A$ denotes thinness and $Jf$ is the Jacobian determinant of
$f$. To prove~\eqref{eJaclocal}, we may by continuity 
assume $\nu$ is divisorial and use the geometric interpretation
in Section~\ref{sec-dual}. With the help of Proposition~\ref{Pinterdiv}, 
the proof reduces to the change of variables formula.
The details are left to the reader.

Continuing the proof of the theorem, if $\nu_\star=\nu_C$ is a curve
valuation, pick $\nu_0$ as above such that $m(\nu_0)=m(C)$.
Apply~\eqref{eJaclocal} with $\nu=\nu_{C,t}$ for large $t$. Then
$A(\nu)=mt+B$ for some constant $B$ (see Section~\ref{sec-mult}),
$c(f,\nu)=c_\infty>1$ and $A(f_\bullet\nu)>A(\nu)$, so we get
$\nu(Jf)>(c_\infty-1)(mt+B)\to\infty$ as $t\to\infty$. Thus $C$
belongs to the critical set of $f$ and in particular is an analytic
curve.

This completes the proof of Theorem~\ref{Teigen}.
%
%
%
%
\section{Local normal forms}\label{sec-normal}
In this section, we prove a quite precise version of
the rigidification statement in Theorem~C.
As an easy consequence we will also obtain Theorem~A.
The key step in achieving rigidification involves 
finding suitable basins of attraction of the eigenvaluation: 
see Proposition~\ref{Plocalcon}.

Recall that a fixed point germ is called \emph{rigid} 
when its critical set is contained in a totally invariant set
with normal crossings.
\begin{Thm}\label{Tnormallocal}
  Let $f$ be a dominant fixed point germ with $c_\infty(f)>1$. 
  Then one can find a modification $\pi:X\to(\C^2,0)$
  and an infinitely near point $p\in\pi^{-1}(0)$ such that
  the lift $\hf$ of $f$ is holomorphic at $p$, 
  $\hf(p)=p$ and $\hf:(X,p)\to(X,p)$ is rigid.
  Moreover, there exist local coordinates $(z,w)$ at $p$ in which
  $\hf$ takes one of the following forms:
  \begin{itemize}
  \item[(i)] $\hf(z,w)=(z^aw^b(1+\phi),\lambda w(1+\psi))$ with
    $a\ge2$, $b\ge1$, $\lambda\in \C^*$ and $\phi(0)=\psi(0)=0$;
    then $c_\infty(f)=a$;
  \item[(ii)]
    $\hf(z,w)=(z^a w^b, z^c w^d)$ with $a,b,c,d\in\N$, $ad\ne bc$;
    then $c_\infty(f)$ is the spectral radius of the $2\times2$
    matrix with entries $a$, $b$, $c$ and $d$;
  \item[(iii)]
    $\hf(z,w)=(z^a,\lambda z^c w^d)$ with $d\ge a\ge 2$, $c\ge1$,
    $\lambda \in \C^*$;
    then $c_\infty(f)=a$;
  \item[(iv)]
    $\hf(z,w) = (z^a, \lambda z^cw + P(z))$ with $a \ge 2$, $c\ge 1$,
    $\lambda \in \C^*$ and $P\not\equiv 0$ is a polynomial;
    then $c_\infty(f)=a$.
  \end{itemize}
  The exceptional divisor is included in $\{ zw= 0\}$.  
  In the first case, it equals $\{ zw=0 \}$. 
  In the other three cases,
  it coincides with the curves that are contracted to $0$ by $\hf$.
\end{Thm}
The four cases arise according to the type of the 
eigenvaluation $\nu_\star$. 
In order, they appear when
$\nu_\star$ is divisorial\footnote{When $|\lambda|<1$
  one can make $\phi\equiv 0$.} 
(the exceptional component is $\{z=0\}$),
irrational, a curve valuation (the curve equals $\pi\{w=0\}$)
or an infinitely singular valuation.

As we do not claim that the eigenvaluation is unique, 
the four cases are a priori not mutually exclusive.
\begin{proof}[Proof of Theorem~A]
  We know from Theorem~\ref{Teigen} that
  the quotient $c_\infty^n/c(f^n)$ 
  is bounded from above and below by positive constants.
  By inspecting the local forms in Theorem~\ref{Tnormallocal} 
  it is clear that $c_\infty$ is a quadratic integer.
\end{proof}
%
%
\subsection{Basins of attraction}
\begin{figure}[ht]
  \begin{center}
    \includegraphics[width=\textwidth]{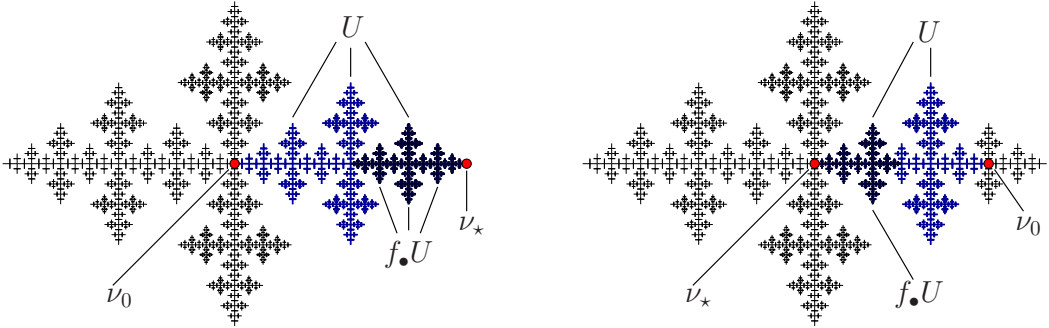}
  \end{center}
  \caption{Basins of attraction of an infinitely
    singular or curve eigenvaluation (left)
    and a divisorial eigenvaluation (right).
    See cases~(i) and~(ii) in Proposition~\ref{Plocalcon}.}\label{F3}
\end{figure}

\begin{figure}[ht]
  \begin{center}
    \includegraphics[width=\textwidth]{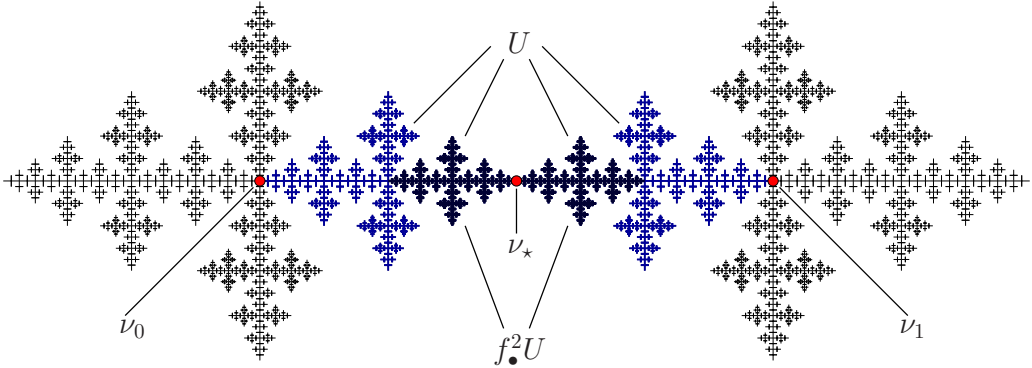}
  \end{center}
  \caption{A basin of attraction of an irrational
    eigenvaluation. See case~(iii) 
    in Proposition~\ref{Plocalcon}.}\label{F4}
\end{figure}

The first step in the proof of Theorem~\ref{Tnormallocal} is
to find small invariant regions in $\cV$ that serve as
basins of attraction for the eigenvaluation.
We need some flexibility in the choice of these basins; 
this explains the somewhat technical appearance of the result below.
\begin{Prop}\label{Plocalcon}
  Let $f$ is a dominant fixed point germ with $c_\infty(f)>1$. 
  Let $\nu_\star$ be an eigenvaluation for $f$ 
  whose existence is guaranteed by Theorem~\ref{Teigen}.
  \begin{itemize}
  \item[(i)]
    If $\nu_\star$ is infinitely singular 
    or an analytic curve valuation, then for any
    $\nu_0\in\cV$ with $\nu_0<\nu_\star$ and $\nu_0$
    sufficiently close to $\nu_\star$, $f_\bullet$ maps 
    the segment $I=[\nu_0,\nu_\star]$ strictly into itself
    and is order-preserving there.
    Moreover, set $U=U(\vv)$ 
    where $\vv$ is the tangent vector at $\nu_0$
    represented by $\nu_\star$.
    Then $f_\bullet$ also maps the open set $U$
    strictly into itself and $f^n_\bullet\to\nu_\star$
    as $n\to\infty$ on $U$.
  \item[(ii)]
    If $\nu_\star$ is divisorial, then there exists 
    a tangent vector $\ww$ at $\nu_\star$ such that for
    any $\nu_0\in\cV$ representing $\ww$ and sufficiently close
    to $\nu_\star$, $f_\bullet$ maps 
    the segment $I=[\nu_\star,\nu_0]$ strictly into itself
    and is order-preserving there.
    Moreover, set $U=U(\vv)\cap U(\ww)$,
    where $\vv$ is the tangent vector at $\nu_0$ 
    represented by $\nu_\star$.
    Then $f_\bullet(U)\subset U$.
    Further  $f^n_\bullet\to\nu_\star$ as $n\to\infty$ on $U$.
  \item[(iii)]
    If $\nu_\star$ is irrational, then 
    there exist $\nu_0,\nu_1\in\cV$, arbitrarily
    close to $\nu_\star$,  with
    $\nu_0<\nu_\star<\nu_1$ such that $f_\bullet$ maps 
    the segment $I=[\nu_0,\nu_1]$ into itself.
    Let $\vv_1$ ($\vv_2$) 
    be the tangent vector at $\nu_0$ ($\nu_1$) 
    represented by $\nu_\star$ and set 
    $U=U(\vv_1)\cap U(\vv_2)$. 
    Then $f_\bullet U\subset U$.
    Further, either $f^2_\bullet=\id$ on $I$ or
    $f^n_\bullet\to\nu_\star$ as $n\to\infty$ on $U$.
  \end{itemize}
\end{Prop}
The different basins $U$ are illustrated in 
Figures~\ref{F3} and~\ref{F4}.
\begin{Remark}
  In the case~(iii) of an irrational eigenvaluation, $f_\bullet$ may be
  order-reversing near $\nu_\star$ and not all segments containing
  $\nu_\star$ are $f_\bullet$-invariant: an $f_\bullet$-invariant
  segment has to be roughly symmetric around $\nu_\star$. 
  However, $f_\bullet^2$ is always 
  order-preserving and it is in fact much
  easier to prove Theorem~\ref{Tnormallocal} for $f^2$ rather 
  than $f$ (Lemma~\ref{L-quadra} and Lemma~\ref{L-genmul} below 
  are not needed.)
\end{Remark}

\begin{Example}
  Let $f(x,y)=(x^2y^3+x^7,x^7)$. 
  The segment $I=[\nu_\fm,\nu_{y,3/2}]$ consists of 
  monomial valuations of the form $\nu_{y,t}$, $1\le t\le3/2$,
  with $\nu_{y,t}(x)=1$ and $\nu_{y,t}(y)=t$. 
  We have $f_\bullet I\subset I$ and the restriction of 
  $f_\bullet$ to $I$ is given by $t\mapsto 7/(2+3t)$.
  One can check that the open subset 
  $U\subset\cV$ consisting of all valuations $\nu\in\cV$
  with $\nu(x)=1$ and $1<\nu(y)<3/2$ is invariant 
  under $f_\bullet$.

  There exists an eigenvaluation $\nu_\star$ belonging to $I$, which
  corresponds to $t=t_\star=\frac13 (\sqrt{22}-1)$.  We have
  $c_\infty(f)=c(f,\nu_\star)=3t_\star+2=1+\sqrt{22}$.  The basin of
  attraction $U$ is of the form $U=U(p)$, where $p$ is an infinitely
  near point on a model obtained by three blowups. In suitable
  coordinates at $p$, the lift of $f$ takes the form
  $\hf(z,w)=(z^2w^3,z^7)$.
\end{Example}
\begin{proof}[Proof of Proposition~\ref{Plocalcon}]
  When $\nu_\star$ is not quasimonomial,
  Theorem~\ref{Teigen} immediately implies all
  the relevant assertions in~(i) above.
  
  For the case of a quasimonomial eigenvaluation, 
  we rely on the fact that $f_\bullet$ is piecewise
  M{\"o}bius in the parameterization of $\cV$ by skewness.
  In doing so, we shall use the following elementary
  result, whose proof is left to the reader.
  \begin{Lemma}\label{Lmc2}
    Set $M(t)=\frac{at+b}{ct+d}$ where $a,b,c,d\in\N$, 
    $ad\ne bc$. Suppose $M(t_\star) = t_\star$ with 
    $t_\star>0$. 
    Then either $M\circ M\equiv\mathrm{id}$, or $t_\star$ 
    is locally attracting \ie $|M'(t_\star)|<1$.
  \end{Lemma}
  Now suppose $\nu_\star$ is irrational. 
  Pick valuations $\mu_0$, $\mu_1$ close to $\nu_\star$ 
  with $\mu_0<\nu_\star<\mu_1$.
  By Theorem~\ref{Ttreemap}, 
  we may assume that $f_\bullet$ is monotone on
  the segment $J:=[\mu_0,\mu_1]$, 
  and $\a(f_\bullet\nu)$ is a M{\"o}bius function of
  $\a(\nu)$ with nonnegative integer coefficients. 
  The segment $f_\bullet J$ contains $\nu_\star$.  
  The intersection $J\cap f_\bullet J$ 
  is thus a nontrivial segment $I'$
  and $f_\bullet I'\subset J$. 
  Apply Lemma~\ref{Lmc2} on $I'$. 
  
  Suppose first $f^2_\bullet=\id$ on $I'$. If $f_\bullet = \id$, we
  set $\nu_0 \= \mu_0$, and $\nu_1 \= \mu_1$. The segment 
  $[\nu_0,\nu_1]$ is clearly invariant. Otherwise $f^2_\bullet = \id$ but
  $f_\bullet \neq \id$, so that $f_\bullet$ is order-reversing. 
  We choose $\nu_0 < \nu_\star$ in the segment $I'$, and we set 
  $\nu_1=f_\bullet\nu_0$.
  The segment $I \= [ \nu_0, \nu_1]$ then contains $\nu_\star$
  and is $f_\bullet$-invariant (and $f_\bullet^2=\id$ on $I$).  
  When instead $f^2_\bullet \neq \id$, $\nu_\star$ is attracting and 
  for any two
  valuations $\nu_0< \nu_\star< \nu_1$ at the same distance to
  $\nu_\star$, the segment $I \= [\nu_0, \nu_1]$ is invariant. 
  Moreover, $f^n_\bullet\to\nu_\star$ on $I$ in this case.  
  Note also that in all cases, 
  we may choose $\nu_0$ and $\nu_1$ arbitrarily
  close to $\nu_\star$.

  Now consider the region $U$ associated to $I$. 
  Pick a curve valuation $\nu_\psi>\nu_1$. 
  Then $U=\{\nu\ ; \ \nu_0<\nu\wedge\nu_\psi<\nu_1\}$.
  Write $f^{*}\psi=\prod\psi_k$ with $\psi_k\in\fm$ irreducible.
  The curve valuations $\nu_{\psi_k}$ are the preimages of
  $\nu_\psi$. After shrinking $I$ we may assume that none of them
  lies in $U$. 
  We may also assume that the critical tree $\cT_{f}$
  of $f$ intersects $U$ only along $I$ (or not at all).
  Thus $\mu\mapsto\mu(\psi_k)$ and $\mu\mapsto c(f,\mu)$
  are constant on any segment $[\tilde{\nu},\nu]$,
  where $\nu\in U$ and $\tnu\=\nu\wedge\nu_\psi\in I$. 
  We infer
  \begin{equation*}
    (f_\bullet \nu)(\psi) 
    =\frac{\nu(f^{*}\psi) }{c(f,\nu)} 
    =\frac{\sum \nu (\psi_k)}{c(f,\nu)} 
    =\frac{\sum \tilde{\nu}(\psi_k)}{c(f,\tilde{\nu})}
    =(f_\bullet \tilde{\nu})(\psi).
  \end{equation*}  
  Thus 
  $f_\bullet\nu\wedge\nu_\psi=f_\bullet\tnu\wedge\nu_\psi
  \in f_\bullet I\subset I$.
  This implies $f_\bullet U\subset U$. 
  When $f_\bullet^2\neq\id$ on $I$, 
  $f^n_\bullet\to\nu_\star$ on $I$, hence also on $U$.

  \smallskip
  Suppose finally that $\nu_\star$ is divisorial.  By the remark
  following Proposition~\ref{Ptangentinter} there exists a tree
  tangent vector $\ww$ at $\nu_\star$ invariant and non-critical under
  the tangent map. Now pick $\mu_0\in\cV$ representing this tangent
  vector and close enough to $\nu_\star$ so that $f_\bullet$ is
  monotone on the segment $J:=[\nu_\star,\mu_0]$.  This is always
  possible by Theorem~\ref{Ttreemap}.  
  Then $I':=f_\bullet J\cap J\ne\emptyset$.
  If $f_\bullet\not\equiv\id$ on $I'$, then we 
  may argue as in the irrational case to construct 
  $I$ and $U$ with the desired properties.

  To conclude, let us therefore assume that
  $f_\bullet=\id$ on $I'$ and derive a contradiction,
  previewing arguments to appear in Section~\ref{Srigidif}.
  The function $c(f,\cdot)$ must be constant on $I'$,
  say $c(f,\cdot)\equiv k$ for some integer $k>1$.
  (See Proposition~\ref{Pinterdiv}.) 
  We may find a valuation $\nu_0\in I'$ near $\nu_\star$
  and a modification
  $\pi:X\to(\C^2,0)$ such that $\nu_\star$ and $\nu_0$ are
  associated to exceptional components $E$, $E_0$ intersecting
  in some point $p$ corresponding to $\ww$. 
  Then the lift 
  $\hf:X\dashrightarrow X$ of $f$ is holomorphic at $p$ and
  $\hf(p)=p$. Moreover, $\hf^*E=kE$ and $\hf^*E_0=kE_0$, forcing
  $\hf|_E$ to be critical at $p$, a contradiction.
\end{proof}
%
%
\subsection{Rigidification}\label{Srigidif}
We now prove Theorem~\ref{Tnormallocal}.  Let $\nu_\star$ be an
eigenvaluation for $f$ as in Theorem~\ref{Teigen}.  
The idea is to show that the basins of attractions found in
Proposition~\ref{Plocalcon} can be chosen as $U=U(p)$ for an
infinitely near point $p$.  This will imply that a suitable lift $\hf$
is holomorphic at $p$. Choosing $U$ small enough will allow us to
control the critical set and make $\hf$ rigid. The normal forms follow
from the analysis in~\cite{F-rigid} and
the formulas for $c_\infty$ from direct computations.
We split into cases depending on the type of $\nu_\star$.

\smallskip
First suppose $\nu_\star=\nu_C$ is a noncontracted analytic 
curve valuation. 
Pick $\nu_0$ as in Proposition~\ref{Plocalcon}.  
By increasing $\nu_0$, we may assume it is divisorial.  Pick a
modification $\pi:X\to(\C^2,0)$ such that $\nu_0=\nu_{E_0}$ for some
exceptional component $E_0\subset\pi^{-1}(0)$ and
such that $\pi^{-1}(C)$ has simple normal crossings.  There is a
unique exceptional component $E\subset\pi^{-1}(0)$ that intersects
the strict transform $\tC$ of $C$, say at the point $p\in E$.  
It follows from~\cite[Proposition~6.32]{valtree} that the divisorial
valuation $\nu_E$ associated to $E$ satisfies $\nu_E<\nu_C$ and
$\nu_F\notin\,]\nu_E,\nu_C[$ for any $F\subset\pi^{-1}(0)$.  In
particular $\nu_0\le\nu_E<\nu_C$.  The region $U=U(p)\subset\cV$
consisting of valuations whose center on $X$ is $p$ coincides with the
set of valuations $\nu\in\cV$ representing the same tangent vector as
$\nu_\star$ at $\nu_E$.  By Proposition~\ref{Plocalcon} we have
$f_\bullet U\Subset U$.

It now follows from Proposition~\ref{Pintergeo} that the lift 
$\hf\=\pi^{-1}\circ f\circ\pi$ is holomorphic at
$p$ and $\hf(p)=p$.
Notice that $\nu_0$ could be chosen arbitrarily close to $\nu_C$. 
Hence we may assume that $U$ does not contain any 
curve valuations associated to an irreducible component 
of the critical set, nor any preimages of 
$\nu_C$, except $\nu_C$ itself.
The critical set of $\hf$ is therefore 
equal to $E\cup\tC$, which is totally invariant and
has simple normal crossings. Thus $\hf$ is rigid at $p$.

To get the normal form of $\hf$ at $p$, we note that $\hf$ is
superattracting at $p$. Indeed, $E$ is contracted to $p$
since $f_\bullet\nu_E>\nu_E$, and the restriction of
$\hf$ to $\tC$ has topological degree 
$c(f,\nu_C)=c_\infty>1$ by Proposition~\ref{Pintercurve}.
We may thus apply the results from~\cite{F-rigid}. 
The critical set of $\hf$ has two branches:
$E$, which is contracted to $p$, and $\tC$, which is fixed. 
Thus $\hf$ is of Type~6, see~\cite[Table I, p.478]{F-rigid},
so $\hf(z,w)=(z^a,z^cw^d)$ with $a\ge2$, $c\ge 1$
and $d\ge2$. In fact, since the curve valuation $\nu_C$
is attracting, a simple calculation gives $d\ge a$.
The argument also shows that $c_\infty=1$.

\smallskip
Now suppose $\nu_\star$ is infinitely singular. Again pick
$\nu_0<\nu_\star$ divisorial as in Proposition~\ref{Plocalcon}.  The
multiplicity function $m=m(\nu)$ is nondecreasing, integer-valued and
unbounded on the segment $[\nu_\fm,\nu_\star[$, 
see Section~\ref{sec-mult}.
We may therefore find a divisorial valuation
$\nu\in\,]\nu_0,\nu_\star[$ where the multiplicity jumps.
By~\cite[Proposition~6.40]{valtree} there exists a modification
$\pi:X\to(\C^2,0)$ such that $\nu=\nu_E$ is associated to an
exceptional component $E\subset\pi^{-1}(0)$ and such that the center
of $\nu_\star$ on $X$ lies on $E$ and is a \emph{smooth} point
$p\in\pi^{-1}(0)$ (in~\cite{valtree}, such a point is called a free
point).  As above we may set $U=U(p)$ and conclude that 
$f_\bullet U\Subset U$.

Again the lift 
$\hf\=\pi^{-1}\circ f\circ\pi$ is holomorphic at
$p$ and $\hf(p)=p$. 
We may ensure that $U$ does not contain any 
curve valuations associated to irreducible components 
of the critical set. This implies that $\hf$ is rigid
at $p$, and that the critical set of $\hf$ is
equal to $E$, which is contracted to $p$.
Comparing with~\cite{F-rigid} 
we see that $\hf$ is of Type~4, \ie
$\hf(z,w)=(z^a,\lambda z^cw+P(z))$ for some polynomial $P$
in suitable coordinates.
Since $f^n_\bullet\to\nu_\star$ on $U$, $\hf$ does not
leave any curve invariant, hence $P\not\equiv0$.

A direct computation gives $c(f,\nu_E)=a$ and in fact also
$c(f^n,\nu_E)=a^n$, for any $n\ge 1$.
But $c(f^n,\nu_E)/ c(f^n)$ is bounded from above and below 
by positive constants, so $c_\infty=\lim_{n\to\infty}c(f^n,\nu_E)^{1/n}=a$.

\smallskip
Next suppose $\nu_\star$ is divisorial.  Pick $\nu_0$ as in
Proposition~\ref{Plocalcon}; we may assume it is divisorial.  Then
pick a modification $\pi:X\to(\C^2,0)$ such that both
$\nu_\star=\nu_E$ and $\nu_0=\nu_{E_0}$ are associated to exceptional
components $E,E_0\subset\pi^{-1}(0)$.  There is then a unique
exceptional component $F\subset\pi^{-1}(0)$ intersecting $E$
and such that $\nu_F\in\,]\nu_E,\nu_0]$
(see~\cite[Corollary~6.32]{valtree}).  Replace $\nu_0$ by $\nu_F$.
The set $U$ in Proposition~\ref{Plocalcon} is of the form
$U=U(p)$, where $E\cap F=\{p\}$.

As before, $f_\bullet U\subset U$, which implies that
the lift $\hf$ of $f$ is holomorphic at $p$ and $\hf(p)=p$. 
By moving the original divisorial valuation $\nu_0$ closer
to $\nu_\star$ we may assume that $U$ does not
contain any curve valuations associated with 
irreducible components of the critical set of $f$.
Then the critical set of $\hf$ is locally 
included in the exceptional divisor $E\cup F$, which is
totally invariant since $f_\bullet U\subset U$.
Hence $\hf$ is rigid at $p$.

As for the normal form at $p$, we may assume
that $p$ is a noncritical fixed point for the restriction of 
$\hf$ on $E$: see the remark after Proposition~\ref{Ptangentinter}.
Choose local coordinates $(z,w)$ at $p$ 
such that $E=\{z=0\}$, and $F=\{w=0\}$. 
Using the fact that $\hf^{-1}E\subset E\cup F$, 
$\hf E=E$, and $\hf F=p$, we easily get that $\hf$ can be 
written in the form in~(i).
Proposition~\ref{Pinterdiv} with $E'=E$ then
implies $c_\infty=c(f,\nu_E)=a$. 

\smallskip
Finally suppose $\nu_\star$ is irrational.  We rely on the
following result, which strengthens Proposition~\ref{Plocalcon} above.
\begin{Lemma}\label{L-key-quadra}
  There exists a modification $\pi:X\to(\C^2,0)$ and exceptional
  components $E_0,E_1\subset\pi^{-1}(0)$, intersecting 
  in some point $p\in\pi^{-1}(0)$, such that:
  \begin{itemize}
  \item
    the divisorial valuations $\nu_0$ and $\nu_1$, associated to
    $E_0$ and $E_1$, respectively, satisfy 
    $\nu_0<\nu_\star<\nu_1$;
  \item
    the open set $U\=\{\nu\ ;\ \nu_0<\nu\wedge\nu_1<\nu_1\}$ 
    is $f_\bullet$-invariant.
  \end{itemize}
  Moreover, $\nu_0$ and $\nu_1$ can be chosen to be arbitrarily 
  close to $\nu_\star$.
\end{Lemma}
We give a proof below.
The set $U$ of the lemma is of the form $U=U(p)$. 
As $U$ is $f_\bullet$-invariant, 
we conclude as before that the lift $\hf$ of
$f$ is holomorphic and rigid at $p$.

To get the normal form of $\hf$ at $p$, we proceed as follows.  
Pick local coordinates $(z,w)$ such that $E_0=\{z=0\}$ and $E_1=\{w=0\}$. 
Since $E_0\cup E_1$ is locally totally invariant under $\hf$
we have $\hf(z,w)=(z^aw^b\phi,z^cw^d\psi)$, where $a,b,c,d\in\N$ and
$\phi,\psi\in\cO_p\setminus\fm_p$.  For $s,t\ge0$ let $\mu_{s,t}$ be
the monomial valuation on $\cO_p$ in coordinates $(z,w)$ with
$\mu_{s,t}(z)=s$ and $\mu_{s,t}(w)=t$.  Then the segment
$I=[\nu_0,\nu_1]$ is parameterized by $\pi_*\mu_{s,t}$ as
$b_{E_0}s+b_{E_1}t=1$.  
Now $\pi_*\mu_{1,0}=b_{E_0}\nu_0$ and
$\pi_*\mu_{0,1}=b_{E_1}\nu_1$: see Section~\ref{S22}.  
Also note that
$\hf_*\mu_{1,0}=\mu_{a,b}$ and $\hf_*\mu_{0,1}=\mu_{c,d}$.  Since
$f_\bullet I\subset I$, we get $a,b,c,d\ge0$.  But $f_\bullet$ is
injective on $[\nu_0,\nu_1]$, so $\mu_{a,c}$ and $\mu_{b,d}$ are
not proportional, \ie $ad\ne bc$.  We have proved that $\hf$ is of
Type~6 in~\cite{F-rigid}, hence can be conjugated
to a monomial map.

Denote by $\rho$ the spectral radius 
of the $2\times 2$ matrix $M$ with entries $a,b,c,d$.  
If $a_n,b_n,c_n,d_n$ are the entries of $M^n$, then
for any $s,t>0$, the sequence $\{(sa_n+tc_n)/\rho^n\}_1^\infty$
is bounded from above and below by positive constants.
We have 
\begin{multline*}
  c(f^n,\nu_0)
  =\nu_0(f^{n*}\fm)
  =b_{E_0}^{-1}\div_z(\pi^* f^{n*}\fm)
  =b_{E_0}^{-1}\div_z(\hf^{n*}\pi^*\fm)=\\
  =b_{E_0}^{-1}\div_z(\hf^{n*}(z^{b_{E_0}}w^{b_{E_1}}))
  =b_{E_0}^{-1}(a_n b_{E_0}+c_n b_{E_1}).
\end{multline*}
Thus $c_\infty=\lim c(f^n,\nu_0)^{1/n}\to\rho$.

\smallskip
To conclude the proof of Theorem~\ref{Tnormallocal},
we now prove Lemma~\ref{L-key-quadra}.
\begin{proof}[Proof of Lemma~\ref{L-key-quadra}]
  We first apply Proposition~\ref{Plocalcon} and pick valuations
  $\mu_0<\nu_\star<\mu_1$ (arbitrarily) close to $\nu_\star$, such
  that $f_\bullet I\subset I$ with $I\=[\mu_0,\mu_1]$.  On $I$ the
  multiplicity is constant, equal to $m\ge1$.  
  We rely on the following two lemmas.
  \begin{Lemma}\label{L-quadra}
    Let $\abar,\bbar,\cbar,\dbar\ge0$ be integers.  Suppose
    $M(t)\=(\abar t+\bbar)/(\cbar t+\dbar)$ has an irrational fixed
    point $t_+>0$. Then there exist arbitrarily large
    integers $p,q,p',q'$ such that $p/q<t_+<p'/q'$, $p'q-pq'=1$ and
    $M$ maps the interval $(p/q, p'/q')$ into itself.
  \end{Lemma}
  \begin{Lemma}\label{L-genmul}
    Let $\nu$ be a divisorial valuation of multiplicity $m\ge1$, of
    skewness $p/m^2q$ with $p,q \in \N^*$, and generic multiplicity
    $b$. Then $b\le mq$.
  \end{Lemma}
  Lemma~\ref{L-genmul} is a direct consequence 
  of~\cite[Lemmas~3.53, 3.54]{valtree}; the details 
  are left to the reader. 
  A proof of Lemma~\ref{L-quadra} is given below.

  We continue the proof of Lemma~\ref{L-key-quadra}.
  Denote by $\mu_t$ the unique valuation in $I$ of skewness $t$. 
  By Theorem~\ref{Ttreemap}, we have 
  $\a(f_\bullet\mu_t)=(at+b)/(ct+d)$ for some integers $a,b,c,d\ge0$. 
  We apply Lemma~\ref{L-quadra} to the $4$-tuple 
  $(\abar,\bbar,\cbar,\dbar)=(am^2,bm^4,c,dm^2)$. 
  Note that $t_+\=m^2\a(\nu_\star)>1$ is irrational
  and a fixed point of $M(t)\=(\abar t + \bbar)/(c\bar t+\dbar)$. 
  This gives us integers $p,q,p',q'$.
  Define $\nu_0\=\mu_{p/m^2q}$ and $\nu_1\=\mu_{p'/m^2q'}$. 
  The condition $p/q<t_+<p'/q'$ implies 
  $p/m^2q<\a(\nu_\star)<p'/m^2 q'$, hence $\nu_0<\nu_\star<\nu_1$.  
  The fact that $M$ maps $(p/q,p'/q')$ into itself shows that 
  $f_\bullet]\nu_0,\nu_1[\,\subset\,]\nu_0,\nu_1[\,$.
  When $\nu_0$ and $\nu_1$ are close enough to $\nu_\star$,
  this yields $f_\bullet U\subset U$.

  On $I$, the multiplicity
  equals $m$ and
  $\a(\nu_1)-\a(\nu_0)=(p'q-pq')/(m^2qq')= 1/m^2qq'$.  
  Lemma~\ref{L-genmul} 
  gives $b(\nu_0)\le mq$ and 
  $b(\nu_1)\le mq'$, so that $\a(\nu_1)-\a(\nu_0)\le 1/b(\nu_0)b(\nu_1)$.
  We now apply~\cite[Proposition~6.38]{valtree}: 
  there exists a modification $\pi: X \to (\C^2,0)$ 
  and exceptional components $E_0,E_1\subset\pi^{-1}(0)$, 
  such that $\nu_0$ (resp.\ $\nu_1$) is associated to
  $E_0$ (resp.\ $E_1$), and $E_0$, $E_1$ intersect in a point. 
\end{proof}

\begin{proof}[Proof of Lemma~\ref{L-quadra}]
  The proof was suggested to us by R.~A.~Mollin. 

  For simplicity, write $a,b,c,d$ instead of
  $\abar,\bbar,\cbar,\dbar$. 
  After conjugating by $t\mapsto t^{-1}$ if necessary, we
  may assume $t_+>1$.  
  Write $\delta \= a-d$. Then $t_+$ is a
  root of the quadratic polynomial $ct^2-\delta t-b=0$.
  The other root $t_-$ is negative, since $t_+t_-=-b/c<0$.  
  Define the integer $n\ge0$ by $t_-+n\in(-1,0)$. 
  Then $t_+ +n$ has a periodic continued fraction
  expansion, see~\cite[p.241]{mollin}: $t_++n = [a_0, a_1,\dots]$ with
  $a_{j+l}=a_j$ for all $j$, where $l\ge1$.  Write
  $P_j/Q_j=[a_0,\dots,a_j]$.  Then $P_j/Q_j\to t_++n$ as $j\to\infty$.

  Take $k$ to be an \emph{even} multiple of $l$ in the sequel.  As
  in~\cite[p.240]{mollin} we have $t_++n=[a_0,\dots,a_{k-1};t_++n]
  =\frac{P_{k-1}(t_++n)+P_{k-2}}{Q_{k-1}(t_++n)+Q_{k-2}}$, thus there
  exists $\lambda_k>0$ rational such that 
  $\lambda_kb=P_{k-2}+n(P_{k-1}-Q_{k-2}-n Q_{k-1})$; 
  $\lambda_k\delta = P_{k-1}-Q_{k-2}-2nQ_{k-1}$; 
  and $ \lambda_kc =Q_{k-1}$.
  From~\cite[p.225]{mollin}, we infer 
  $P_{k-1}Q_{k-2}-P_{k-2}Q_{k-1}=(-1)^{k}=1$.  Set
  \begin{equation*}
    \begin{pmatrix}
      p' \\ q'
    \end{pmatrix}
    =   
    \begin{pmatrix}
      P_{k-1} - n Q_{k-1} \\
      Q_{k-1} 
    \end{pmatrix}
    \qand
    \begin{pmatrix}
      p \\ q
    \end{pmatrix}
    = 
    \begin{pmatrix}
      P_{k-2} + n P_{k-1} - n Q_{k-2} -n^2 Q_{k-1} \\
      Q_{k-2} + n Q_{k-1}
    \end{pmatrix}
  \end{equation*}
  A direct computation yields 
  $q'=\lambda_k c$, $p'-q=\lambda_k\delta$, $p'q-pq'=1$ 
  and $p=\lambda_k b$, so
  \begin{equation}\label{e-quadra}
    q^2+q\lambda_k\delta-\lambda_k^2 bc=1.  
  \end{equation}
  Clearly $q', p$ and $q$ are positive, and $p'q = 1 + pq'$
  implies $p'>0$.
  Notice that $p'/q'\to t_+$ as $k\to\infty$, hence
  $p/q=p'/q'-1/qq'\to t_+$ as well.
  From the fact $P_{k-1}/Q_{k-1} > t_+ + n $, we
  infer $p'/q' > t_+$, thus 
  $p/q=\lambda_k b/(p'-\lambda_k\delta)
  <\lambda_kb/(\lambda_kct_+-\lambda_k\delta)=t_+$. 
  To conclude, it suffices to prove that $M$ maps
  $(p/q, p'/q')$ into itself, or more precisely
  $p'/q'> M(p'/q') > p/q$, and $p'/q'> M(p/q) > p/q$. 
  Now 
  \begin{equation*}
    M(p'/q') - p/q =
    \frac{aq^2 + a \delta \lambda_k q - \lambda_k^2 bc (\delta +d)}{q(
      cp' + dq')}\ \mathop{=}\limits^{\text{by~\eqref{e-quadra}}}\
    \frac{a}{q( cp' + dq')}>0~,
  \end{equation*}
  \begin{equation*}
    p'/q'-M(p'/q')
    = \frac{bq'^2 - cp'^2  - \delta  p'q'}{q'( cp' + dq')} 
    =\frac{q' \left(- b + c(p'/q')^2  - \delta  (p'/q')\right)}{cp' + dq'}.
  \end{equation*}
  But $p'/q'>t_+$ and $t\mapsto ct^2-\delta t-b$ is
  increasing near $t_+$, thus $p'/q'-M(p'/q')>0$. 
  Hence $p'/q'>M(p'/q')> p/q$.  
  Similarly, $p'/q'-M(p/q)=\frac{d}{q'(cp+dq)}>0$
  and $M(p/q)-p/q>0$ (for $k$ large).
  Thus $M(p/q,p'/q')\subset(p/q,p'/q')$.
\end{proof}
%
%
%
%
\section{Proof of Theorem~B}\label{sec-green}
For background on plurisubharmonic (psh) functions, 
see \eg~\cite{dem0}.
In~\cite{pshsing} we showed that any
(quasimonomial) valuation $\nu$ 
may be naturally evaluated on any psh function $u$
and that $\nu\mapsto\nu(u)$ defines a tree 
potential---the tree transform of $u$.
Here we shall only need the case when 
$u=\log\max_i|\phi_i|^{t_i}$,
for holomorphic germs $\phi_i$, and $t_i>0$, $i=1,2$,
in which case $\nu(u)=\min_it_i\nu(\phi_i)$.

We prove Theorem~B in several steps, the proofs of which are given at
the end of the section. The first lemma below is an improvement
of~\cite[Proposition~4.2]{pshsing}.
\begin{Lemma}\label{Lpullpot}
  If $g:\cVqm\to[0,\infty[$ is a tree potential, then
  so is the function $f^*g:\cVqm\to[0,\infty[\,$, defined by
  $(f^*g)(\nu):=c(f,\nu)g(f_\bullet\nu)$.
\end{Lemma}
\begin{Lemma}\label{Lfindgood}
  There exists a quasimonomial valuation $\nu_0\in\cV$ 
  such that $f_\bullet\nu_0\ge\nu_0$ and $c(f,\nu_0)=c_\infty$.
  In particular, $f_*\nu_0\ge c_\infty\nu_0$.
\end{Lemma}
Define tree potentials $g_0$ and $g_1$ by
$g_0(\nu)=\a(\nu\wedge\nu_0)$ and $g_1:=c_\infty^{-1}f^*g_0$.
\begin{Lemma}\label{Lestg0g1}
  We have $g_1\ge g_0$.
\end{Lemma}
Since $\nu_0$ is quasimonomial, we can write $\nu_0=\nu_{\phi,t}$ 
with $t>0$, $\phi\in\fm$ irreducible and $m(\phi)=m(\nu)=:m$.
Pick a coordinate $x$ transverse to $\phi$, 
\ie $x\in\fm$ is irreducible, $m(x)=1$ and $\nu_x\wedge\nu_\phi=\nu_\fm$.
Define the psh function 
$u_0\=\log\max\{|\phi|^{1/m},|x|^t\}$
and set $u_1=c_\infty^{-1}u_0\circ f$.
\begin{Lemma}\label{Ltransfg0g1}
  The tree transforms of $u_0$ and $u_1$ are
  $g_0$ and $g_1$, respectively. 
  Hence $\nu(u_0)\le\nu(u_1)$ for every $\nu\in\cVqm$.
\end{Lemma}
\begin{Lemma}\label{Lgoodpsh}
  We have $u_1\le u_0+O(1)$ and
  $\delta\log\|\cdot\|\le u_0\le\log\|\cdot\|+O(1)$
  for some $\delta>0$.
\end{Lemma}
By subtracting a constant from $u_0$ (\ie multiplying 
$x$ and $\phi$ by constants), we may assume $u_1\le u_0$.
It is then clear that $u_n\=c_\infty^{-n}u_0\circ f^n$ defines
a decreasing sequence of psh functions, hence 
converges to a function $u_\infty$. 
\begin{Lemma}\label{Lnontriv}
  The function $u_\infty$ is not identically $-\infty$. 
  In particular it defines a psh function satisfying the 
  invariance relation $u_\infty\circ f=c_\infty u_\infty$.
\end{Lemma}
This lemma completes the proof of Theorem~B.
\begin{proof}[Proof of Lemma~\ref{Lpullpot}]
  First, if $g(\nu)=\nu(\psi)$ for some $\psi\in R$, then 
  $f^*g$ is a tree potential as 
  $(f^*g)(\nu)=c(f,\nu)(f_\bullet\nu)(\psi)=\nu(f^*\psi)$.
  Second, if $g=g_I$ with $g_I(\nu)=\nu(I)$ for some ideal 
  $I\subset R$, then $f^*g$ is the minimum of tree potentials,
  hence a tree potential. 
  Now, a general tree potential $g$ is
  the (pointwise) limit of tree potentials of the form 
  $c\,g_I$, where $c>0$ and $I$ is an ideal.
  It follows that $f^*g$ is a
  limit of tree potentials, hence a tree potential.
\end{proof}
\begin{proof}[Proof of Lemma~\ref{Lfindgood}]
  Let $\nu_\star$ be an eigenvaluation as in Theorem~\ref{Teigen}.
  When $\nu_\star$ is quasimonomial, set $\nu_0:=\nu_\star$. 
  Otherwise, $\nu_\star$ is an end and we may pick
  $\nu_0$ as in the statement of Theorem~\ref{Teigen}.
\end{proof}
\begin{proof}[Proof of Lemma~\ref{Lestg0g1}]
  Clearly,
  $g_1(\nu_0)=c_\infty^{-1}c_\infty\a(f_\bullet\nu_0\wedge\nu_0)
  =\a(\nu_0)=g_0(\nu_0)$. 
  As $g_1$ is a tree potential, 
  $g_1(\nu_0)\le\a(\nu_0)g_1(\nu_\fm)$, so
  $g_1(\nu_\fm)\ge1=g_0(\nu_\fm)$.
  Since $g_1$ is concave and $g_0$ affine on the segment
  $[\nu_\fm,\nu_0]$, we get $g_1\ge g_0$ there.
  Outside this segment, $g_0$ is locally constant
  and $g_1$ is nondecreasing, so $g_1\ge g_0$ everywhere.
\end{proof}
\begin{proof}[Proof of Lemma~\ref{Ltransfg0g1}]
  The tree transform of $u_0$ is $g_0$. Indeed, as 
  $\nu_0=\nu_{\phi,t}$:
  \begin{equation*}
    \nu(u_0)
    =\min\{m^{-1}\nu(\phi),t\nu(x)\}
    =\min\{\a(\nu\wedge\nu_\phi),t\a(\nu\wedge\nu_x)\}
    =\a(\nu\wedge\nu_0)
    =g_0(\nu)
  \end{equation*}
  for any $\nu\in\cVqm$.
  The special form of $u_0$ directly implies that
  the tree transform of $u_0\circ f$ is $f^*g_0$
  (see also~\cite[Proposition~4.2]{pshsing}),
  which completes the proof.
\end{proof}
\begin{proof}[Proof of Lemma~\ref{Lgoodpsh}]
  Notice that $u_1=\log\max_i|\phi_i|^{t_i}$ for holomorphic
  germs $\phi_1,\phi_2$ and $t_1,t_2>0$.
  Pick a modification $\pi:X\to(\C^2,0)$ 
  such that the total transform of the curve 
  $\{\phi\phi_1\phi_2=0\}$ has simple normal crossings.
  For any point $p\in\pi^{-1}(0)$, we can find
  local coordinates $(z,w)$ at $p$ such that 
  $\phi(z,w)=z^kw^{l}$ for some $k,l\in\N$.
  This implies
  \begin{equation*}
    u_0\circ\pi(z,w)
    =c_0\log|z|+d_0\log|w|+\log\max\{|z|^{r_0},|w|^{s_0}\}+O(1)
    \end{equation*} 
  for some constants $c_0,d_0,r_0,s_0\ge 0$. 
  These constants can be computed in terms of the tree transform
  of $u_0$. For instance, if $p$ is the intersection point between
  two exceptional components $E=\{z=0\}$ and $F=\{w=0\}$ with
  associated divisorial valuation $\nu_E$, $\nu_F$ and 
  generic multiplicities $b_E$, $b_F$,
  then $c_0=b_E\nu_E(u_0)$, $d_0=b_F\nu_F(u_0)$.
  Further, $r_0$ and $s_0$ may be computed in terms of 
  the restriction of the tree transform to $[\nu_E,\nu_F]$,
  using the monomialization procedure in Section~\ref{S22}.

  Similarly, we have
  \begin{equation*}
    u_1\circ\pi(z,w)
    =c_1\log|z|+d_1\log|w|+\log\max\{|z|^{r_1},|w|^{s_1}\}+O(1)
    \end{equation*} 
  for some constants $c_1,d_1,r_1,s_1\ge 0$. 
  The fact that $\nu(u_0)\le\nu(u_1)$ for all quasimonomial $\nu$
  implies that $c_0\le c_1$, $d_0\le d_1$, $r_0\le r_1$ and $s_0\le s_1$.
  We conclude that $u_1\circ\pi\le u_0\circ\pi+O(1)$ 
  in a neighborhood of $p$. As this is true for any $p$
  on the exceptional divisor of $\pi$, this inequality can be pushed
  down to a neighborhood of the origin in $\C^2$, in which we get
  $u_1\le u_0+O(1)$.

  The same proof applied to 
  $u_0$ and $u=\log\|\cdot\|$ 
  gives $\a(\nu_0)^{-1} u+O(1)\le u_0\le u+O(1)$.
  Indeed, $1=\nu(u)\le\nu(u_0)\le\a(\nu_0)$ for 
  all $\nu\in\cVqm$.
\end{proof}

\begin{proof}[Proof of Lemma~\ref{Lnontriv}]
  This is the key point in the proof of Theorem~B. 
  Let $\nu_\star$ be an eigenvaluation as in Theorem~\ref{Teigen}.
  We use the normal forms in Theorem~\ref{Tnormallocal}.

  When $\nu_\star=\nu_C$ is an analytic curve valuation, the normal
  form is $\hf(z,w)=(z^a,\lambda z^cw^d)$.  Here $E:=\{z=0\}$ is an
  exceptional component, whereas $\{w=0\}$ is the strict transform of
  $C$.  Write $\hf^n(z,w)=(z_n,w_n)$.  We have $c_\infty=a$
  so it is clear that
  $c_\infty^{-n}\log|z_n|=\log |z|$ for any point in a neighborhood
  $\Omega$ of $p$.  Moreover, $\log\|\pi(z,w)\|=b_E\log|z|+O(1)$ and
  $u_0\ge\delta\log\Vert\cdot\Vert+O(1)$ so $u_\infty=\lim
  c_\infty^{-n}u_0\circ f^n$ is finite in $\pi(\Omega\setminus E)$.

  The case when $\nu_\star$ is infinitely singular 
  is completely analogous
  and left to the reader. When $\nu_\star$ is irrational, the normal
  form is $\hf(z,w)=(z^aw^b,z^cw^d)$, and one checks that
  $\|\hf^n(z,w)\|\ge c(z,w)\, \rho^n$ when $zw\neq0$, 
  for some $c(z,w)>0$, where $\rho = c_\infty(f)$ is the 
  spectral radius of the $2\times2$ matrix with entries $a,b,c,d$.

  \smallskip
  The final---and most delicate---case is when $\nu_\star$ is
  divisorial. The problem is that we cannot apply
  Theorem~\ref{Tnormallocal} directly: for the normal form~(i) there need
  not exist any invariant open region attracted to the origin.
  
  By Proposition~\ref{Ptangentinter}, 
  the tangent map $Df_\bullet$ at $\nu_\star$ may be viewed as
  a surjective rational selfmap of $\P^1$.
  We consider three different cases.

  In the first case, the tangent map has topological degree at least
  two. It then has infinitely many \emph{repelling} periodic
  orbits~\cite{CG}.  Pick one of them, such that for any tangent
  vector $\vv$ in this orbit, the open set $U(\vv)$ does not intersect
  the critical tree of $f$.  Replacing $f$ by a suitable iterate we
  may assume that $\vv$ is a fixed point for the tangent map.  Then
  $f_\bullet$ maps $U(\vv)$ into itself.  We may assume that $\vv$ has
  multiplicity $m(\vv)=b(\nu_\star)$ (this is true for all but at most
  two $\vv$).  By~\cite[Proposition~6.40]{valtree} there exists a
  modification such that $\nu_\star$ is associated with an exceptional
  component $E$ of $\pi$, $\vv$ is associated with a point $p\in E$
  that is a smooth point on $\pi^{-1}(0)$. Then $U(\vv)=U(p)$.
  Hence $\hf$ is holomorphic at $p$. 
  Now $\hf(E)=E$, and $p$ is a repelling fixed point for the
  restriction $\hf|_E$. On the other hand, $\hf^*E=c_\infty E$ and
  $c_\infty>1$, so $E$ is superattracting. Thus there exists a stable
  manifold through $p$, \ie a holomorphic curve $V$ which is smooth
  and transverse to $E$ and invariant under $\hf$, 
  see~\cite[Ex~III.2, p.68]{Shub}.  
  Notice that $V$ is \emph{not} contracted to $p$.  It
  is now easy to put $\hf$ in the form
  $\hf(z,w)=(z^a(1+\phi),\lambda w(1+\psi))$, where
  $\phi(0)=\psi(0)=0$, $a=c_\infty$, $|\lambda|>1$, $E=\{z=0\}$ and
  $V=\{w=0\}$.  Arguing as before we see that $u_\infty\ne-\infty$ on
  $\pi(V\setminus E)$.

  In the second case, the tangent map is a M{\"o}bius map with
  two distinct fixed points.
  At least one of these must be 
  non-repelling. Use this tangent vector in the proof of 
  Theorem~\ref{Tnormallocal} to get the normal form
  \begin{equation*}
    \hf(z,w)=\left(z^aw^b (1+\phi(z,w)),\lambda
    w(1+z\psi(z,w))\right),
  \end{equation*}
  with $a=c_\infty\ge2$, $b\ge1$, $0<|\lambda|\le 1$ 
  and $\phi(0)=0$. Here we have used the fact that
  a M{\"o}bius map of $\P^1$ with two fixed points
  is \emph{linear} in some affine coordinate.
  Fix $\e>0$ small and pick $\e_0>0$ such that
  $|\phi(z,w)|\le\e$, $|\psi(z,w)|\le \e/4$ when $|z|,|w|\le\e_0$.
  A simple induction shows that if $0<|z_0|,|w_0|\le\e$, then
  $|z_n|\le 2^{-n}\e$ and $|2\lambda|^{-n}|w_0|\le|w_n|\le 2|w_0|$.
  It then follows that 
  $a^{-(n+1)}\log|z_{n+1}|=a^{-n}\log|z_n|+O(na^{-n})$.
  Thus $c_\infty^{-n}\log|z_n|\notto-\infty$ and
  $c_\infty^{-n}\log|w_n|\to0$.
  This implies $u_\infty(\pi(z_0,w_0))\ne-\infty$.

  In the third and final case, the tangent map is a M{\"o}bius
  map with a single (parabolic) fixed point.
  We can then obtain the normal form
  \begin{equation*}
    \hf(z,w)=\left(z^aw^b(1+\phi(z,w)),\frac{w}{1+w}(1+z\psi(z,w))\right),
  \end{equation*}
  with $a=c_\infty>1$, $b\ge1$ and $\phi(0)=0$.
  A direct computation shows that for small $\e$, the
  region defined by $0\le|\mathrm{Im}\,w|<|\mathrm{Re}\,w|<\e$
  and $0<|z|<|w|^3$ is invariant and that 
  $|w_n|\ge2^{-n}$ in this region.
  As in the previous case, this gives
  $a^{-(n+1)}\log|z_{n+1}|=a^{-n}\log|z_n|+O(na^{-n})$,
  which allows us to conclude.
\end{proof}
%
%
%
%
\section{The affine case}\label{sec-affine}
Now we turn to polynomial maps and prove Theorems~A' and~B'. 
The approach mimics the one in the local case
and is based on the action on the set $\cV_0$ of 
normalized valuations centered at infinity.
This set has a tree structure similar to that of the
valuative tree $\cV$, but the induced map on
$\cV_0$ is not everywhere defined. 
We therefore also work with an invariant
subtree $\cV_1$ of $\cV_0$:
this is the main difference to the local case.
The tree structure of the spaces $\cV_0$ and $\cV_1$ is
detailed in the appendix.
%
%
\subsection{Induced action: domain of definition}\label{affine1}
For the rest of the paper, $R$ denotes the ring of polynomials
in two complex variables. 
The role of the maximal ideal $\fm$ in the local case is played 
by the vector space $\cL$ of affine functions on $\C^2$.
We refer to $X,Y\in\cL$ as \emph{coordinates} if $R=\C[X,Y]$.

The valuation space $\cV_0$ is by definition the set of valuations 
$\nu:R\to(-\infty,\infty]$ normalized by 
$\nu(\cL):=\min\{\nu(L)\ ;\ L\in\cL\}=-1$. 
We refer to Section~\ref{sec-full} for
a description of its tree structure. Suffice it to say that
$\cV_0$ comes equipped with a partial ordering $\le$, 
a multiplicity function $m$, 
a decreasing parameterization 
$\a:\cV_0\to[1,-\infty]$ by skewness and
an increasing parameterization 
$A:\cV_0\to[-2,\infty]$ by thinness.
The valuation $\affroot$ is the root of $\cV_0$.

If $F$ is a dominant polynomial map of $\C^2$ 
and $\nu$ is a non-constant valuation on $R$, 
then $(F_*\nu)(P)\=\nu(F^*P)$, $P\in R$ defines
a nonconstant valuation $F_*\nu$.
The proof of Proposition~\ref{Pgen} applies also 
in this context, so that $F_*$ sends divisorial (resp.\ irrational,
curve or infinitely singular) valuations to valuations of the same
type.

For $\nu\in\cV_0$ set
\begin{equation}\label{e-defd}
  d(F,\nu)
  \=-\nu(F^*\cL)\=-\min\{\nu(F^*L)\ ;\ L\in\cL\}.
\end{equation}
Thus $d(F,\nu)\ge0$ and $F_*\nu$ 
is centered at infinity exactly when $d(F,\nu)>0$.
We then write $F_*\nu=d(F,\nu)F_\bullet\nu$ 
with $F_\bullet\nu\in\cV_0$. 
Clearly $d(F,\affroot)=\deg(F)$.
\begin{Def}
  We denote by $\cD_F^o$ the set of valuations $\nu\in\cV_0$ 
  such that $F_*\nu$ is centered at infinity, 
  and let $\cD_F$ be its closure in $\cV_0$.
\end{Def}
Thus $\cD_F^o=\{d(F,\cdot)>0\}$.
We summarize in the following proposition the structure 
of the set $\cD_F$ and the behavior of the function
$d(F,\cdot)$ on $\cV_0$. Compare with Proposition~\ref{Pcrittree}
in the local case.
\begin{Prop}\label{prop-fntd}
  The function $d(F,\cdot)$ is continuous and decreasing on $\cV_0$.
  The locus $\cT_F$ where $d(F,\cdot)$ is not locally constant is a
  finite subtree of $\cV_0$.  The ends of $\cT_F$ are the maximal
  elements of the finite set $\cE_F$
  consisting of divisorial valuations $\nu\in\cV_0$
  satisfying one of the following properties:
  \begin{itemize}
  \item[(i)]
    $d(F,\nu)>0$ and $F_\bullet\nu=\affroot$;
  \item[(ii)]
    $d(F,\nu)= 0$, $F_*\nu=c\div_C$, where $c>0$ and
    $C\subset\C^2$ is a rational 
    curve with one place at infinity.
   \end{itemize} 
   The valuations in~(ii) are always maximal in $\cE_F$. 
   Moreover, $\cD_F^o$ and $\cD_F$ are subtrees of $\cV_0$,
   and $\cD_F\setminus\cD_F^o$ is a finite set
   consisting of all divisorial valuations satisfying~(ii).
\end{Prop}

Recall that a curve $C\subset \C^2$ has 
\emph{one place at infinity} when
its closure in $\P^2$ intersects the line at infinity at a
single point $p$, and admits a single branch there.

\begin{proof}[Proof of Proposition~\ref{prop-fntd}]
  We have 
  $d(F,\nu)=-\min\{\nu(F^*X),\nu(F^*Y),0\}$
  in coordinates $(X,Y)$. 
  Thus the function
  $d(F,\cdot)$ is clearly decreasing and continuous.
  It is moreover a tree potential in the sense 
  of Section~\ref{sec-pot-affine}; this implies that $\cT_F$
  is a tree.

  The elements in $\cE_F$ are divisorial by Proposition~\ref{Pgen};
  let us prove that they are finite in number.
  Pick a modification $\pi:S\to\P^2$ such
  that $F$ lifts to a holomorphic map $\hF:S\to\P^2$. 
  The image of a divisorial valuation $\nu$ centered at 
  a (closed) point $p\in S$ is a
  divisorial valuation centered at the point $\hF(p)$.
  In particular, $F_*\nu$ is proportional to neither
  $\affroot$ nor $\div_C$, where $C\subset\C^2$.
  So any valuation in $\cE_F$
  is associated to an irreducible component of 
  $\pi^{-1} (L_\infty)$, where $L_\infty$ 
  denotes the line at infinity. 
  Thus $\cE_F$ is finite.

  We now show that the ends of $\cT_F$
  are exactly the maximal elements of $\cE_F$. 
  Pick an end $\nu$ of $\cT_f$. 
  If $d(F,\nu)>0$, then $F_*\nu$ is centered at
  infinity and we may argue as in the proof of
  Proposition~\ref{Pcrittree} to 
  show that $F_\bullet\nu=\affroot$
  and that $\nu$ is a maximal element of $\cE_F$.
  Hence assume $d(F,\nu)=0$.
  Then $F_*\nu$ is not centered at infinity, so its center 
  is either a point $p\in\C^2$  or an irreducible 
  curve $C\subset\C^2$. 

  In the first case, pick affine coordinates 
  $(X,Y)$ such that $p$ is the origin. 
  Then $\nu(F^*X),\nu(F^*Y)>0$. The same
  inequality must hold in a neighborhood of $\nu$. Thus
  $d(F,\cdot)\equiv0$ on this neighborhood, contradicting
  $\nu\in\cT_F$.
  
  If instead $F_*\nu$ is centered at an irreducible curve
  $C\subset\C^2$, then $F_*\nu$ is proportional to $\div_C$. 
  It remains to be shown that $C$ is a rational
  curve with one place at infinity.
  Pick a modification $\pi:S\to\P^2$ such that $\nu$ 
  is associated to an exceptional component $E$ of $\pi$ and
  $F$ lifts to a holomorphic map $\hF:S\to\P^2$. 
  Then $\hF(E)$ equals $\bar{C}$, the closure of $C$ in $\P^2$.
  As $E$ is a rational curve, so is $C$.
  If $\bar{C}$ has several branches at infinity, one can
  find two distinct points $p,q\in E$ such that 
  $\hF(p),\hF(q)\in L_\infty$.
  Denote by $U(p)$ (resp.\ $U(q)$) the set of valuations 
  in $\cV_0$ whose center is $p$ (resp.\ $q$). 
  These are disjoint open subsets of $\cV_0$ containing $\nu$ in
  their closures,
  so at least one of them contains a valuation $\mu>\nu$.
  Then $F_*\mu$ is still centered at infinity
  but $d(F,\mu)\le d(F,\nu)=0$, so $d(F,\mu)=0$, a
  contradiction.
  We conclude that $C$ has one place at infinity.
  In particular, $\nu\in\cE_F$. 

  Next we show that all valuations satisfying~(ii) 
  are maximal elements of $\cE_F$ as
  well as ends in $\cT_f$.
  This will in particular complete the proof that all 
  ends in $\cT_F$ are maximal elements in $\cE_F$.
  Hence pick $\nu$ satisfying~(ii).
  Clearly $d(F,\nu)=0$. It suffices to show that 
  $d(F,\cdot)\not\equiv0$ at $\nu$.
  In the notation above, let $p'$ 
  be the unique intersection point of $\bar{C}=\hF(E)$ with the
  line at infinity and pick $p\in\hF^{-1}(p')\cap E$.
  We may approximate $\nu$ by valuations $\mu$ whose center on 
  $S$ is the point $p$. Since $\hF$ is holomorphic,
  $F_*\mu$ is centered at $p'\in L_\infty$, so $d(F,\mu)>0$. 
  Hence $d(F,\cdot)\not\equiv0$ at $\nu$.

  Next consider a divisorial valuation $\nu$ satisfying~(i) 
  and maximal in $\cE_F$.  
  As $F_\bullet\nu=\affroot$, $F_\bullet$ cannot be order
  preserving near $\nu$, thus $\nu\in\cT_F$. 
  It is an end of $\cT_F$, since we already proved that
  all ends of $\cT_F$ belong to $\cE_F$.

  That $\cD_F^o$, and hence $\cD_F$, is a subtree of $\cV_0$
  follows from the fact that $d(F,\cdot)$ is
  decreasing. If $\nu\in\cD_F\setminus\cD_F^o$, then $d(F,\nu)=0$,
  but $d(F,\cdot)$ is not locally constant at $\nu$. This implies
  that $\nu$ is an end of $\cT_F$.
\end{proof}
%
%
\subsection{Induced action: regularity}\label{affine2}
Having described its domain of definition, we now show that
the induced map on the valuation space exhibits many of the
features present in the local case.
The following two results together generalize 
Theorem~\ref{Ttreemap}.
\begin{Thm}\label{thm-deffpt}
  The map $F_\bullet:\cD_F^o\to\cV_0$ extends uniquely
  to a surjective, regular tree map $F_\bullet:\cD_F\to\cV_0$.
  The valuations in $\cD_F\setminus\cD_F^o$ are mapped to
  curve valuations associated with rational curves having one place
  at infinity.
\end{Thm}
The set $\cV_1$ is 
the closure of the set $\cV_1^o$
of valuations $\nu\in \cV_0$ such that 
$A(\nu)<0$ and $\nu(P)<0$ for all nonconstant $P$.
See Section~\ref{sec-v1} for details.
\begin{Thm}\label{thm-v1inv}
  The subtree $\cV_1\subset\cV_0$ 
  is included in $\cD_F$, and $F_\bullet\cV_1\subset \cV_1$.
  As a consequence, 
  $F_\bullet:\cV_1\to\cV_1$ is a regular tree map.

  Moreover, for any $\nu\in\cV_1$, we can decompose the segment
  $[\affroot,\nu]$ into finitely many subsegments
  $I_j=[\nu_j,\nu_{j+1}]$, $0\le j\le k$, with 
  $\nu_0=\affroot$, $\nu_{k+1}=\nu$
  and $\nu_i$  divisorial, such that $F_\bullet$ is a 
  monotone homeomorphism of
  $I_j$ onto its image, and $\a(F_\bullet\cdot)$ is a piecewise
  M{\"o}bius function of $\a(\cdot)$ with nonnegative 
  integer coefficients.
\end{Thm}
\begin{Remark}
  As follows from the proof of Theorem~\ref{thm-deffpt},
  the map $F_\bullet:\cD_F\to\cV_0$ is also piecewise M{\"o}bius
  along segments, but the coefficients may be negative in general.
  Indeed, let $F(X,Y)=(X,Y^2(XY-1))$ and 
  denote by $\nu_t$ the minimal valuation with $\nu_t(X)=-1$,
  $\nu_t(Y)=1$ and $\nu_t(XY-1)=t>0$. It has skewness $-t-1$.
  Now $F_\bullet\nu_t$ is monomial, with values $-1$ and $2+t$ on 
  $X$ and $Y$ respectively. Its skewness is $-t-2$.
  Hence $\a(F_\bullet\nu_t)=\a(\nu_t)-1$.
\end{Remark}
\begin{proof}[Proof of Theorem~\ref{thm-deffpt}]
  The proof is similar to the local case, so we shall
  skip some of the details and focus on the new features in
  the affine situation.

  Let $\nu$ be an end in $\cD_F$. 
  We need to study the behavior of $F_\bullet$ 
  on the closed segment $[\affroot,\nu]$ when $\nu\in\cD_F^o$,
  and on the half-open segment $[\affroot,\nu[$ 
  when $\nu\in\cD_F\setminus\cD_F^o$. 
  In both cases we need to decompose the segment into
  finitely many subsegments on each of which $F_\bullet$
  is a monotone homeomorphism.

  First suppose $\nu\in\cD_F\setminus\cD_F^o$. 
  We know from Proposition~\ref{prop-fntd} that $\nu$ is divisorial
  and that $F_*\nu$ is proportional to $\div_C$, where
  $C\subset\C^2$ is a rational curve with one place at infinity.
  Let $\bar{C}$ denote the closure of $C$ in $\P^2$.
  We may pick modifications $\pi:S\to\P^2$ and 
  $\pi':S'\to\P^2$ such that $\pi^{-1}(C\cup L_\infty)$
  has simple normal crossings and such that
  the lift $\hF:S\to S'$ of $F$ is holomorphic and 
  locally monomial at any point on $\pi^{-1}(L_\infty)$.
  We may also assume that 
  $\nu$ is associated to an irreducible component $E$
  of $\pi^{-1}(L_\infty)$ that is mapped by $\hF$ onto $\tC$, 
  the strict transform of $\bar{C}$ by $\pi'$.
  
  Let $\affroot<\nu_1<\dots<\nu_n=\nu$ be the
  divisorial valuations in $[\affroot,\nu]$
  whose centers on $S$ are one-dimensional. 
  As in the proof of Theorem~\ref{Ttreemap},
  $F_\bullet$ is a monotone homeomorphism of 
  $[\nu_i,\nu_{i+1}[$ onto its image for $0\le i<n$.
  Moreover, $F_\bullet\mu\to\nu_C$ as $\mu$ increases 
  to $\nu$.

  Next consider the case when $\nu$ is an analytic 
  curve valuation in $\cD_F$.
  The proof just given still applies
  in exactly the same way, the only difference being
  that $\nu$ is mapped to a curve valuation of the same type.

  If instead $\nu$ is a formal curve valuation or an infinitely singular
  valuation, then $\nu$ does not belong to the critical tree $\cT_F$
  of $F$. Thus there exists a divisorial valuation $\nu_0<\nu$ such that
  $F_\bullet$ is increasing and continuous, hence a homeomorphism, on
  $[\nu_0,\nu]$.  This completes the proof that $F_\bullet$ is a
  regular tree map, as $\nu_0$ is dominated by an analytic curve
  valuation in $\cD_F$.

  \smallskip
  To see that $F_\bullet$ is surjective, it suffices to show
  that any divisorial valuation $\nu'\in\cV_0$ has a preimage. 
  Pick modifications $\pi:S\to\P^2$, $\pi':S'\to \P^2$ such that 
  $\nu'$ is associated to some exceptional component 
  $E'\subset S'$, and $F$
  lifts to a holomorphic map $\hF:S\to S'$. 
  As $\P^2$ is compact, $\hF$ is surjective, 
  so one can find an irreducible curve $E\subset S$ 
  such that $\hF(E)=E'$. 
  Since $F$ is a polynomial map of $\C^2$, the curve $E$ 
  is mapped by $\pi$ into the line at infinity in $\P^2$,
  hence the divisorial valuation $\nu$ 
  associated to $E$ is centered at infinity. 
  We can normalize it by 
  $\nu(F^*\cL)=-1$, so that $\nu\in\cV_0$. 
  Then $F_\bullet\nu=\nu'$.
  This finishes the proof of the theorem.
\end{proof}
\begin{proof}[Proof of Theorem~\ref{thm-v1inv}]
  If $\nu\in\cV_1^o$, then $\nu(P)<0$ for any 
  nonconstant $P$. This immediately implies 
  $d(F,\nu)>0$ so $\nu\in D_F^o$. 
  Taking closures yields $\cV_1\subset\cD_F$.

  The fact that $F_\bullet$ leaves $\cV_1$ invariant relies
  on the following affine version of the Jacobian
  formula~\eqref{eJaclocal}.
  As in the local case, it ultimately reduces to the
  change of variables formula.

  \begin{Lemma}\label{LJacaffine}
    For any valuation $\nu\in\cD_F$, we have
    \begin{equation}\label{eJacaffine}
      \nu(JF)+A(\nu)=d(F,\nu)\,A(F_\bullet\nu),
    \end{equation}
    where $JF$ denotes the Jacobian determinant of $F$ in $\C^2$.
  \end{Lemma}
  If $\nu\in\cV_1^0$, then $\nu(JF)<0$, $A(\nu)<0$ and $d(F,\nu)>0$,
  hence $A(F_\bullet\nu)<0$. On the other hand, it is clear
  that $(F_\bullet\nu)(P)=\nu(F^*P)/d(F,\nu)<0$ for
  all nonconstant polynomials $P$.
  Thus $F_\bullet\nu\in\cV_1^o$, so $F_\bullet\cV_1^o\subset\cV_1^o$.
  By continuity, $F_\bullet\cV_1\subset\cV_1$.

  To control skewness, we proceed as in the local case. First assume
  that $\nu\in\cV_1$ is not infinitely singular.  
  By Theorem~\ref{thm-deffpt},
  $[\affroot,\nu]$ can be decomposed into finitely many segments 
  $[\nu_i,\nu_{i+1}]$ on which $F_\bullet$ is monotone. 
  All the $\nu_i$ can be chosen divisorial. 
  The key remark is that since $F_\bullet[\nu_i,\nu_{i+1}]$ is a
  segment containing no infinitely singular valuation, 
  we can pick an affine curve $C=\{Q=0\}$ with one place at infinity 
  such that the
  associated curve valuation $\nu_C$ dominates $F_\bullet\mu$ for all
  $\mu\in[\nu_i,\nu_{i+1}]$ (see Theorem~\ref{thm-struct-affine}). 
  Then 
  \begin{equation}
    \a(F_\bullet\mu)
    =\frac{-F_\bullet\mu(Q)}{\deg(Q)}
    =\frac{-F_*\mu(Q)}{d(F,\mu)\cdot\deg(Q)}
    =\frac{-\mu(F^*Q)}{d(F,\mu)\cdot\deg(Q)}.
  \end{equation}
  The function
  $\mu\mapsto-\mu(F^*Q)$ is piecewise linear in $\a(\mu)$ 
  with nonnegative integer coefficients, see Lemma~\ref{L-act-pol}. 
  The same is true for the function 
  $d(F,\mu)=-\min\{\mu(F^*X),\mu(F^*Y),0\}$. 
  Thus $\a(F_\bullet\mu)$ is a M{\"o}bius function 
  of $\a(\mu)$ with nonnegative integer coefficients, as claimed.

  When $\nu$ is infinitely singular, $F_\bullet$ is monotone on
  $[\nu_0, \nu]$ for $\nu_0$ divisorial close enough to $\nu$. 
  For $\mu \in [\nu_0, \nu]$ divisorial, we have seen that the segment
  $[\affroot, \mu]$ can be decomposed into finitely many segments on
  which the function $\a(F_\bullet \cdot)$ is M{\"o}bius with
  nonnegative integer coefficients. The number of such segments is
  bounded by the number of irreducible components of $F^*C$ (with $C$
  as above), which via Bezout's Theorem can be bounded 
  by $\deg(F)^2$. 
  As this bound is uniform in $\mu$, the segment
  $[\affroot, \nu]$
  can also be decomposed into finitely many segments on which
  $\a(F_\bullet \cdot)$ is M{\"o}bius. This completes the proof.
\end{proof}
%
%
\subsection{Normal forms}\label{affine3}
We now work towards the proof of Theorems~A' and~B'.
Let $F:\C^2\to\C^2$ be a dominant polynomial map.
The analogue of~\eqref{Ecompconlocal} is 
\begin{equation}\label{Ecompconaffine}
  \deg(F^n)=\prod_{k=0}^{n-1}d(F,F^k_\bullet(\affroot)).
\end{equation}
Since $d(F,\cdot)$ is decreasing, 
the sequence $\deg(F^n)$ is submultiplicative
(as can also be verified directly), so
the limit $d_\infty:=\lim\deg(F^n)^{1/n}$ exists.

The quantity $d_\infty$ is invariant under conjugacy:
if $F_2=G\circ F_1\circ G^{-1}$ for a polynomial automorphism $G$,
then $D^{-1}\le\deg(F_1^n)/\deg(F_2^n)\le D$ for all $n$, where 
$D=\deg(G)\deg(G^{-1})$, hence $d_\infty(F_1)=d_\infty(F_2)$.

In order to control the degrees $\deg(F^n)$ we show that 
$F$ can be made rigid.
\begin{Thm}\label{Tnormalaffine}
  Let $F:\C^2\to\C^2$ be a dominant polynomial map
  with $d_\infty(F)>1$. 
  Assume $F$ is not conjugate to a skew product.
  Then one can find a modification $\pi:S\to\P^2$, 
  and $p\in\pi^{-1}(L_\infty)$ such that
  the lift $\hF$ of $F$ is holomorphic at $p$, 
  $\hF(p)=p$ and $\hF:(S,p)\to(S,p)$ is rigid.
  Further, there exist local coordinates $(z,w)$ at $p$ in which
  $\hF$ takes one of the following forms:
  \begin{itemize}
  \item[(i)]
    $\hF(z,w)=(z^aw^b(1+\phi),\lambda w^d(1+\psi))$ with
    $b\ge1$, $a>d\ge1$, $\lambda\in \C^*$ and $\phi(0)=\psi(0)=0$; 
    then $d_\infty=a$;
  \item[(ii)]
    $\hF(z,w)=(z^aw^b,z^cw^d)$ with $a,b,c,d\in\N$, 
    $ad\ne bc$;
    then $d_\infty$ is the spectral radius of the 
    $2\times2$ matrix with entries $a,b,c,d$;
  \item[(iii)]
    $\hF(z,w)=(z^a,\lambda z^cw+P(z))$ with $a\ge 2$, $c\ge 1$,
    $\lambda \in \C^*$ and $P\not\equiv 0$ is a polynomial;
    then $d_\infty=a$.
  \end{itemize}
\end{Thm}
\begin{Remark}\label{Rehoui}
  As in the local case, the proof is based on the construction of 
  an eigenvaluation $\nu_\star\in\cV_1$ (not necessarily unique)
  and the normal form above depends on the type of $\nu_\star$,
  namely divisorial, irrational, or infinitely singular.
  In the normal form~(i) we can assume $d=1$ and $|\lambda|\ge1$
  except if $\nu_\star$ is divisorial with $\a(\nu_\star)>0$
  and $A(\nu_\star)=0$ 
  (a case which we suspect may in fact never appear).
\end{Remark}
\begin{proof}
  By Theorem~\ref{thm-v1inv}, $F_\bullet:\cD_F\to\cV_0$ 
  restricts to a regular tree map $F_\bullet:\cV_1\to\cV_1$.
  Theorem~\ref{MT1} implies the existence of a
  fixed point (eigenvaluation) $\nu_\star$, which is 
  either an interior point of $\cV_1$ or 
  a strongly attracting end on $\cV_1$ 
  in the sense of Section~\ref{S23}.
  We now proceed according to the nature of $\nu_\star$:
  it can be divisorial, irrational, or infinitely singular, but
  not a curve valuation.
  We further subdivide the divisorial case into the 
  subcases $\a(\nu_\star)>0$ and $\a(\nu_\star)=0$. 
  The latter will not be discussed here, as we will
  show later that $F$ is then conjugate to a skew product.

  \smallskip
  First assume that $\nu_\star\in\cV_1$ is divisorial, with skewness
  $\a_\star:=\a(\nu_\star)>0$. Then 
  \begin{equation*}
    \deg(F^n)
    \le d(F^n,\nu_\star)
    =d(F,\nu_\star)^n
    =-\nu_\star(F^{n*}\cL)
    \le\a_\star^{-1}\deg(F^n)
  \end{equation*}
  by Theorem~\ref{thm-struct-affine}.
  Thus $d(F,\nu_\star)=d_\infty$ and 
  $\deg(F^n)/d_\infty^n\in[1,\a_\star^{-1}]$ for all $n$.

  We may now prove that a suitable lift $\hF$ of $F$
  may be put in the form~(i) by arguing as in the
  local case. 
  First suppose $\nu_\star$ is not an end of $\cV_1$. 
  Then $F_\bullet$ induces an action on the tree 
  tangent space at $\nu_\star$, and 
  this action admits a noncritical fixed point:
  see the remark after Proposition~\ref{Ptangentinter}.
  By the M{\"o}bius property in Theorem~\ref{thm-v1inv},
  we may construct basins of attraction $U$ of $\nu_\star$
  as in Proposition~\ref{Plocalcon}. Notice that these basins
  are subsets of $\cV_0$, and not $\cV_1$.
  They may be chosen as $U=U(p)$ for some infinitely near point $p$, 
  and by choosing them small enough, the lift $\hF$ of $F$ will
  be holomorphic and rigid at $p$. The normal form~(i) with $d=1$
  is an easy consequence. 

  If instead $\nu_\star$ is an end of $\cV_1$, then the argument 
  above breaks down as the tangent vector constructed may not
  be represented by any valuation in $\cV_1$, so the M\"obius
  property may fail. However, the
  tangent vector represented by the valuation
  $\affroot$ is invariant since $\cV_1$ is invariant. 
  Proceeding as above using this
  tangent vector yields the normal form in~(i).
  We have $a>d$ since $\nu_\star$ is an attracting end in $\cV_1$.

  That $a=d_\infty$ (in both cases) is a simple computation. 
  Indeed, $\hF^*\pi^*\cL$ and $\pi^*F^*\cL$ have poles of order
  $ab_\star$ and $d(F,\nu_\star)b_\star$, respectively,
  where $b_\star$ is the generic multiplicity of $\nu_\star$.
  Hence $d_\infty=d(F,\nu_\star)=a$.

  \smallskip
  Next suppose $\nu_\star$ is infinitely singular. 
  It is by construction a strongly attracting end, 
  so there exists
  $\nu_0<\nu_\star$ such that $d(F,\cdot)$ is constant on 
  $[\nu_0,\nu_\star]$. Arguing as in the local case we see that
  $d_\infty=d(F,\nu_\star)$ and
  $\deg(F^n)/d_\infty^n\in[1,\a(\nu_0)^{-1}]$ for all $n$.
  The construction of the basin $U=U(p)$ and the local
  normal form~(ii) are also proved as in the local case.
  We get $d_\infty=a$ for the same reason.

  \smallskip
  Finally assume $\nu_\star$ is irrational. Then
  $\a_\star:=\a(\nu_\star)>0$ and the argument in the divisorial case
  shows that $d_\infty=d(F,\nu_\star)$ and
  $\deg(F^n)/d_\infty^n\in[1,\a_\star^{-1}]$ for all $n$.  
  Rigidification can be proved as in the local case:
  the most delicate situation is when $F_\bullet$ is 
  order-reversing at $\nu_\star$.
  The proof relies on the analogue of Lemma~\ref{L-key-quadra}
  with details being left to the reader.
  The monomial form of $\hF$ follows from~\cite{F-rigid} and the
  expression for $d(F,\nu_\star)=d_\infty$ is proved as in the local
  case, by pulling back $\cL$. 
\end{proof}
%
%
\subsection{Proof of Theorem~A'}\label{affine4}
We have already done most of the work by finding the normal forms in
Theorem~\ref{Tnormalaffine}.  Indeed, assuming that the eigenvaluation
$\nu_\star$ was not divisorial with $\a(\nu_\star)=0$ we saw that
$\deg(F^n)/d_\infty^n\in[1,D]$ for all $n$ for some $D\ge1$.
Moreover, $d_\infty$ is an integer except when $\nu_\star$ is
irrational. In the latter case, it is the spectral radius of a 
$2\times2$ matrix with integer coefficients, hence a quadratic integer.

Let us therefore assume $\nu_\star$ is a divisorial
valuation with $\a(\nu_\star)=0$. It is then an end of $\cV_1$, 
and, by Theorem~\ref{thm-struct-affine},
associated to a rational pencil of curves with one place at infinity.
The Line Embedding Theorem~\cite{abh-moh, Suzuki}
shows that there exists a change of coordinates by a polynomial 
automorphism such that the pencil becomes $\{X=\mathrm{const}\}$.  
Writing $F=(P,Q)$, we infer
\begin{equation*}
  0=d(F,\nu_\star)\nu_\star(X)
  =\nu_\star(F^*X) 
  =\nu_\star(P) 
  =\min_{c\in \C}-\deg_Y(P(c,Y)).
\end{equation*}
Hence $P$ depends only on $X$, and 
$F(X,Y)=(P(X), Q(X,Y))$ is a skew product. 
As $\nu_\star$ was chosen to be a strongly
attracting end, we get $d_\infty(F) = \deg(P)\ge \deg_Y(Q)$.
It is easy to see that $\deg (F^n)/ d_\infty^n$ is unbounded iff
$\deg_Y(Q) =\deg(P)> \deg_Y(Q(X_0,Y))$ for some $X_0\in\C$.
This completes the proof of Theorem~A'.
%
%
%
%
\subsection{Proof of Theorem~B'}\label{affine5}
The proof follows the same steps as in Theorem~B.

By adapting the proof of Lemma~\ref{Lfindgood}
and using the assumption that $F$ is not conjugate to a skew
product we find a quasimonomial valuation
$\nu_1\in\cV_1$ with $F_\bullet\nu_1\ge\nu_1$,
$\a(\nu_1)>0$ and $d(F,\nu_1)=d_\infty$.
Indeed, either the eigenvaluation $\nu_\star$ has positive skewness,
in which case we take $\nu_1=\nu_\star$, or $\nu_\star$ is an
infinitely singular valuation $\nu_\star$, in which case we take
$\nu_1<\nu_\star$ sufficiently close to~$\nu_\star$.

The function $g_0:\cV_0\to\R$ given by $g_0(\nu)=-\a(\nu\wedge\nu_1)$
is a tree potential. We then prove that 
$g_1(\nu):=d_\infty^{-1}d(F,\nu)g_0(F_\bullet\nu)$
defines a tree potential and that $g_1\ge g_0$.

We know that $\nu_1$ is dominated by a pencil valuation $\nu_{|P|}$
for some polynomial $P$ of degree $m$. Set
$U_0=\log^+\max\{|P|^{1/m},|X|^t\}$, where $t=\a(\nu_1)$ and $X$ is a
generic affine function.  By Corollary~\ref{C-hopelast}, $U_0$ is a
psh function whose tree transform is $g_0$.  Now
$U_1:=d_\infty^{-1}F^*U_0$ is also psh, with tree transform $g_1$.
Since $g_1\ge g_0$ we get $U_1\le U_0$, after
subtracting a constant from $U_0$, if necessary.  We also obtain
$\delta\log^+\|\cdot\|\le U_0\le\log^+\|\cdot\|+O(1)$ for some
$\delta>0$.  Thus $U_n:=d_\infty^{-n} U_0 \circ F^n$ defines a
decreasing sequence of psh functions which are uniformly bounded from
below.  Hence $U_n$ converges to a nonnegative psh function
$U_\infty$.  The proof is completed by finding one point 
$q\in\C^2$ such that $U_\infty(q)>0$, or, equivalently,
$\log\|F^n(q)\|\ge\e d_\infty^n$ for some $\e>0$, and this
is done exactly as in the proof of Theorem~B. 
%
%
%
%
\appendix
\section{Valuations centered at infinity}\label{sec-app}
In this appendix we define and analyze valuation spaces 
adapted to the analysis of polynomial maps of $\C^2$ at infinity.
%
%
\subsection{The full valuation space}\label{sec-full}
Let $R$ be the polynomial ring in two complex variables.
Denote by $\cL$ the vector space of affine functions on $\C^2$.
We refer to $X,Y\in\cL$ as \emph{coordinates} if $R=\C[X,Y]$.
By a \emph{valuation} on $R$ we mean a function 
$\nu:R\to(-\infty,\infty]$ satisfying
$\nu(PQ)=\nu(P)+\nu(Q)$ and 
$\nu(P+Q)\ge\min\{\nu(P),\nu(Q)\}$ for all $P,Q\in R$;
and such that $\nu(0)=\infty$, $\nu|_{\C^*}=0$.

Any valuation extends to the fraction field of $R$, and then
restricts to a valuation on the local ring at any point in 
the compactification $\P^2$ of $\C^2$. Its \emph{center}
is the unique (not necessarily closed) point such that 
the valuation is nonnegative on the corresponding local 
ring and (strictly) positive on its maximal ideal.

In particular, a valuation is \emph{centered at infinity} if there 
exists a polynomial $P\in R$ with $\nu(P)<0$. 
Then $\nu(\cL):=\min\{\nu(L)\ ;\ L\in\cL\}<0$ and
$\nu(\cL)=\min\{\nu(X),\nu(Y)\}$ for any coordinates $(X,Y)$.
\begin{Def}
  We let $\cV_0$ be the set of valuations 
  $\nu:R\to(-\infty,\infty]$ centered at infinity, 
  normalized by $\nu(\cL)=-1$.
\end{Def}
Our first goal is to show that $\cV_0$ is a complete 
tree which can be endowed with a multiplicity function as well
as two natural parameterizations (called skewness and 
thinness as in the local case), see Theorem~\ref{T101} below.

We equip $\cV_0$ with the topology of pointwise convergence 
(in which it is compact), and with the partial ordering defined by
$\mu\le\nu$ iff $\mu(P)\le\nu(P)$ for all polynomials $P$.  The unique
minimal element is the valuation $\affroot$ whose
value on $P$ is $-\deg(P)$.
It is divisorial by definition and 
its center on $\P^2$ is $L_\infty$, the line at infinity.

Any other valuation is centered at a (closed) point $p$ on 
$L_\infty$ and is classified as a divisorial, quasimonomial, irrational,
curve, or infinitely singular valuation according to the analysis in
Section~\ref{sec-background}. 
Pick local coordinates $(z,w)$ at
$p$ such that $(1/z,w/z)$ are affine coordinates on $\C^2$. 
Note that $L_\infty=\{z=0\}$. 
The set of valuations $\nu\in\cV_0$ centered at $p$ 
can be identified with
the set $\tcV_p$ of centered valuation on $\cO_p$,
normalized by $\nu(z)=1$. 
The union of $\tcV_p$ and $\div_z=\affroot$ is referred to as 
\emph{the relative valuation space} in~\cite[Section~3.9]{valtree}.  
It has a natural tree structure induced by the natural partial
ordering, with root $\div_z$. Since this partial ordering is
compatible to the partial ordering on $\cV_0$ we see that 
$\cV_0$ is naturally a complete tree, rooted in $\affroot$.

The relative valuation space $\tcV_p$ comes equipped with 
a multiplicity function and two natural
parameterizations that we now recall.
First, the \emph{relative multiplicity} $\tm(C)$ 
of a (formal) curve $C$ is the intersection multiplicity 
at $0$ of $C$ with $\{z=0\}$. For $P\in\hat\cO_p$ 
set $\tm(P)=\tm(P^{-1}(0))$; this is the 
order of vanishing of $P(0,w)$ at $0$. 
The relative multiplicity of a (quasimonomial) valuation
$\nu$ is the infimum of $\tm(C)$ over all irreducible formal
curves $C$ such that the curve valuation $\nu_C$ 
associated to $C$ dominates $\nu$.  
The \emph{relative skewness} is defined by
$\tilde{\a}(\nu)=\sup\nu(P)/\tm(P)$, and the
\emph{relative thinness} by 
$\tilde{A}(\nu)=1+\int_{\div_z}^\nu\tm(\mu)\,d\tilde{\a}(\mu)$.

Thus we have defined $\tm(\nu)$, $\tilde{\a}(\nu)$ and $\tA(\nu)$
for any $\nu\in\cV_0$, $\nu\ne\affroot$. It turns out to be
more natural to modify these quantities slightly.
\begin{Def}\label{D-glob}
  We define the \emph{multiplicity} $m(\nu)$,
  \emph{skewness} $\a(\nu)$ and \emph{thinness} $A(\nu)$
  of a valuation $\nu\in\cV_0$ as follows. 
  If $\nu\ne\affroot$, then
  \begin{equation*}
    m(\nu):=\tm(\nu),
    \quad
    \a(\nu)\=1-\tilde{\a}(\nu)
    \quad\text{and}\quad    
    A(\nu)\=\tA(\nu)-3,
  \end{equation*}
 and for $\nu=\affroot$ we set $m=1$, $\a=1$, $A=-2$.
\end{Def}
The definition of these invariants may seem arbitrary, 
but are justified by Theorem~\ref{thm-struct-affine}; see
also formula~\eqref{e-726}.
\begin{Thm}\label{T101}
  The valuation space $(\cV_0,\le)$ is a complete nonmetric tree
  rooted in $\affroot$. Moreover, 
  $\a:\cV_0 \to [-\infty, +1]$ (resp.\ $ A:\cV_0\to[-2,+ \infty]$) 
  gives a decreasing (resp.\ increasing)
  parameterization of $\cV_0$; and $m:\cV_0\to\overline{\N}$ 
  is nondecreasing. 
  Further, for any $\nu\in\cV_0$ we have
  $A(\nu)=-2-\int_{\affroot}^\nu m(\mu)\,d\a(\mu)$.
\end{Thm}
\begin{proof}
  These statements follow from the properties of the 
  relative valuation space: 
  see~\cite[Proposition~3.61, Corollary~3.66]{valtree}.
  The integral formula reduces to the formula for $\tA$ 
  in terms of $\tilde{\a}$ and $\tm$.
\end{proof}
We conclude this section by giving a geometric interpretation of 
thinness in the case of a divisorial valuation.  
Any divisorial valuation $\nu\in\cV_0$ is obtained as
follows. First, $\affroot(P)$ is the order of vanishing of the (rational)
function $P$ along the line $L_\infty$ at infinity. 
Second, if $\nu\ne\affroot$, there
exists a modification $\pi:S\to\P^2$ 
and an exceptional component $E\subset\pi^{-1}(L_\infty)$ such that
$\nu(P)=b^{-1}\div_E(\pi^*P)$. Here $\div_E$ denotes order of
vanishing along $E$ and
$b=b(\nu)=-\div_E(\pi^*\cL)$. 
Set $a(\nu):=1+\div_E(\pi^*\Omega)$, where $\Omega=dX\wedge dY$.  
Then
\begin{equation}\label{e-726}
  A(\nu) =a(\nu)/b(\nu).
\end{equation}
Compare this with the corresponding local result in 
Section~\ref{sec-dual}. Formula~\eqref{e-726} can easily be
reduced to the local case, using~\cite[Theorem~6.50]{valtree}.
%
%
\subsection{Pencils}\label{sec-pencil}
An affine  curve $C$ is said to have \emph{one place at infinity} 
if its closure $\P^2$ intersects the line at infinity in a single
point and is analytically irreducible there.  
Such curves have received considerable attention: 
see \eg~\cite{abh-moh,Suzuki,CPR}, and play a crucial role in our study.
The branch of $C$ at infinity defines a curve valuation $\nu_C\in\cV_0$.
If $C=P^{-1}(0)$ for $P\in R$, then we also say that $P$ has one
place at infinity and write $\nu_C=\nu_P$.  
Notice that $\nu_P(P)=\infty$.

By a theorem of Moh~\cite{Moh} (see also~\cite[p.565]{CPR}), the curve
$P^{-1}(\lambda)$ has one place at infinity for any $\lambda\in\C$. 
Denote by $|C|$ or $|P|$ the corresponding pencil of curves in
$\C^2$ and define the function $\nu_{|C|}=\nu_{|P|}$ by
\begin{equation}\label{e104}
  \nu_{|C|}(Q)=\min_{\lambda\in\C}\nu_{P+\lambda}(Q).
\end{equation}
\begin{Prop}\label{L102}
  The function $\nu=\nu_{|C|}$ defines a divisorial valuation centered
  at infinity.  We have $\nu_{|P|}(Q)\le 0$ for all $Q\in R$, and
  $\nu_{|P|}(Q)=0$ iff $Q=c(P+\lambda)$ for some 
  $c\in\C^*$, $\lambda\in\C$.
  Moreover, 
  $A(\nu)=(2g-1)/b$, where $b=b(\nu)$ is the generic
  multiplicity of $\nu$ and $g$ is the genus of a 
  generic element in the pencil $|C|$.
  In particular, $A(\nu)\le 0$ iff $|C|$ is a pencil of 
  rational curves.
\end{Prop}
We shall call any such  valuation a \emph{pencil valuation} and, 
when the generic element of the associated pencil is rational, 
a \emph{rational pencil valuation}.

In the proof we shall use the following lemma, which can
be proved using Bezout's theorem and the interpretation 
of a curve valuation as a local intersection number.
\begin{Lemma}
  If $C$ is a curve with one place at infinity and $\nu_C\in\cV_0$ 
  is the associated curve valuation, then
  \begin{equation}\label{e106}
    \nu_C(Q)=-\frac{(C\cdot D)_{\C^2}}{\deg(C)},
  \end{equation}
  where $D=Q^{-1}(0)$ and $(C\cdot D)_{\C^2}$ denotes the
  total number of intersection points in $\C^2$, counting
  multiplicities.
\end{Lemma}
\begin{proof}[Proof of Proposition~\ref{L102}]
  Pick a modification $\pi:S\to\P^2$ such that the pencil has no base
  point on $S$.
  The strict transform $\tC_\lambda$ of a generic 
  member $C_\lambda$ of
  the pencil is then smooth, and intersects the exceptional divisor
  transversely.  But $C_\lambda$ has one place at infinity, so
  $\tC_\lambda$ intersects the exceptional divisor at a \emph{unique}
  smooth point, say on an exceptional component $E$, independent of
  the choice of (generic) element of the pencil. From this geometric
  representation, it is easy to conclude that $\nu = \nu_{|C|}$ is the
  divisorial valuation associated to $E$. The fact that $\nu$ takes
  only non-positive values follows from~\eqref{e106}. To compute the
  thinness of $\nu$, write $a=1+\div_E(\pi^*\Omega)$, 
  and $b=\div_E (\pi^*\cL)$ as in the previous section.  By the genus
  formula~\cite[p.471]{GH}, the genus $g$ of $\tC_\lambda$ satisfies
  $2g-2=K_S\cdot\tC_\lambda+\tC_\lambda\cdot\tC_\lambda$, where $K_S$
  is a canonical divisor of $S$.  Now $\tC_\lambda\cdot\tC_\lambda=0$
  as $\tC_\lambda$ is an element of a pencil on $S$ without base
  points. Moreover, $K_S$ can be chosen as the divisor of the 2-form
  $\pi^*\Omega$ on $S$. This form has a zero of order
  $a-1$ along $E$. Thus $K_S\cdot\tC_\lambda=a-1$,
  so~\eqref{e-726} gives $2g-1=a=A(\nu)b$.
\end{proof}      
%
%
\subsection{A smaller valuation space}\label{sec-v1}
We shall now define a subtree $\cV_1$ of $\cV_0$ which
is crucial to our study of (possibly nonproper) 
polynomial maps of $\C^2$.
\begin{Def}
  Let $\cV_1^o$ be the set of quasimonomial valuations
  $\nu\in\cV_0$ with negative thinness $A(\nu)<0$ and 
  such that $\nu(P)<0$ for all nonconstant polynomials $P$.
  Let $\cV_1$ be the closure of $\cV_1^o$ in
  $\cV_0$.
\end{Def}
\begin{Thm}\label{thm-struct-affine}
  The space $\cV_1$ is a complete subtree of $\cV_0$. Moreover:
  \begin{itemize}
  \item[(i)]
    every quasimonomial valuation $\nu\in\cV_1$ is dominated by 
    a pencil valuation; 
    the multiplicity $m(\nu)$ is the smallest degree
    of any such pencil;
  \item[(ii)]
    skewness $\a:\cV_1\to[0,1]$ defines a decreasing
    parameterization of $\cV_1$;
    further, for any $\nu\in\cV_1$, we have 
    $\a(\nu)=-\sup\frac{\nu(P)}{\deg(P)}$,
    taken over nonconstant polynomials $P$;
    when $\nu$ is quasimonomial, the supremum is attained by any
    polynomial with one place at infinity 
    whose pencil valuation dominates $\nu$;
  \item[(iii)]
    thinness $A:\cV_1\to[-2,0]$ defines an increasing 
    parameterization of $\cV_1$; 
    further, for any $\nu\in\cV_1$, we have
    $A(\nu)=-2-\int_{\affroot}^\nu m(\mu)\,d\a(\mu)$;
  \item[(iv)]
    the ends of $\cV_1$, none of which is
    contained in $\cV_1^o$, are of the form:
    \begin{itemize}
    \item[(a)]
      rational pencil valuations;
      these have $\a=0$ and $A<0$;
    \item[(b)]
      divisorial valuations with $\a>0$ and $A=0$;
    \item[(c)]
      infinitely singular valuations  
      with $\a\ge0$ and $A\le0$.
    \end{itemize}
  \end{itemize}
\end{Thm} 
We postpone the proof of Theorem~\ref{thm-struct-affine} 
to Section~\ref{sec-skp}, but point out that the 
main difficulty is to prove~(i). Specifically we need to
construct affine curves with one place at infinity.

Notice that the description above of the tree structure
(partial ordering, multiplicity, skewness, thinness)
on $\cV_1$ is purely in terms of polynomials on $\C^2$.
%
%
\subsection{Potential theory}\label{sec-pot-affine}
We can develop a potential theory on $\cV_0$ similar to the one on the
valuative tree $\cV$ as described in Section~\ref{sec-treepot1}. 
In fact, doing so is a purely tree-theoretic endeavor, 
and we refer to~\cite[Chapter~7]{valtree} 
for the precise construction of a
potential theory on an arbitrary tree.  The main feature
distinguishing $\cV_0$ from $\cV$ is that the parameterization
$\a:\cV_0\to[-\infty,1]$ by skewness is \emph{decreasing} in the
partial ordering, and that a general 
(irreducible) polynomial does not determine a
unique valuation in $\cV_0$. This accounts for the
different properties of tree potentials in these two contexts.

The set of tree potentials is the smallest set of functions on
$\cV_0$ that is closed under sums, multiplication by positive real
constants, under minima  and pointwise limits, and contains all
functions of the form $\nu\mapsto -\a(\mu\wedge\nu)$ for any valuation
$\mu\in \cV_0$. 
With this definition any tree potential is
nonpositive on $\cV_1$, but may be positive elsewhere on $\cV_0$.
Alternatively, tree potentials can be
characterized as increasing functions on $\cV_0$ with certain (strong)
concavity properties. In particular, a tree potential restricts to 
an increasing, concave function on any segment $[\affroot,\nu]$ 
parameterized by skewness.

An important property that we shall use is that, 
for any polynomial $P\in R$, the function $\nu\mapsto \nu(P)$ is a
tree potential. More precisely, we have
\begin{Lemma}\label{L-act-pol}
  Let $P\in R$ be a polynomial and let $\nu_i\in\cV_0$, $i=1,\dots,d$
  be the curve valuations associated the to irreducible branches of
  $P^{-1}(0)$ at infinity. 
  Then there exist integers $m_i\ge1$ 
  such that $\sum_1^d m_i=\deg(P)$ and 
  \begin{equation}\label{e109}
    \nu(P)=-\sum_1^d m_i\,\a(\nu\wedge\nu_i)
    \quad\text{for all $\nu\in \cV_0$}.
  \end{equation}
\end{Lemma}
The proof is given below. 
In particular, on any segment, 
the function $\nu\mapsto\nu(P)$ 
is piecewise linear with integer coefficients as a
function of skewness. Moreover, it does not attain its local
maximum at a quasimonomial valuation unless it is locally
constant there.

For the proof of Theorem~B', we need
\begin{Cor}\label{C-hopelast}
  Suppose a polynomial $P\in R$ has one place at infinity.
  For $t\in[0,1]$, let $\mu_t$ be the unique valuation in
  $[\affroot,\nu_P]$ with skewness $t$.
  Then
  \begin{equation}\label{e107}
    \min\{-\nu(P)/\deg (P),t\}=-\a(\nu\wedge\mu_t)
    \quad\text{for any $\nu\in\cV_0$}.
  \end{equation}
\end{Cor}
\begin{proof}
  Lemma~\ref{L-act-pol} shows that
  $\nu(P)=-\deg(P)\a(\nu\wedge\nu_P)$ for all $\nu\in\cV_0$.
  The proof reduces to showing that
  $\a(\nu\wedge\mu_t)=\min\{\a(\nu\wedge\nu_P),t\}$.
  Both sides of this equation are 
  locally constant off the segment $[\affroot,\nu_P]$,
  and in fact also off the subsegment $I=[\affroot,\mu_t]$.
  On $I$ they both coincide with $\a(\nu)$, which completes
  the proof.
\end{proof}
\begin{proof}[Proof of Lemma~\ref{L-act-pol}]
  By linearity, we may assume $P$ to be irreducible and reduced.  
  Denote by $C_1, \dots, C_d$ the (analytic) branches of the curve 
  $C\=P^{-1}(0)\subset\P^2$ at infinity, 
  set $m_i\=C_i\cdot L_\infty$, 
  and let $\nu_i\in \cV_0$ be the curve valuation associated to $C_i$. 

  Then $\sum_1^dm_i=\deg(P)$ by Bezout's theorem. 
  In particular~\eqref{e109} holds for $\nu=\affroot$.
  Now consider $\nu\ne\affroot$ and let 
  $p\in L_\infty$ be the center of $\nu$ on $\P^2$,
  so that $\nu\in\tcV_p$.
  As above, choose coordinates $(z,w)$ at $p$
  such that $(1/z, w/z)$ are affine coordinates in $\C^2$. 
  Thus $L_\infty=\{z=0\}$. 
  Set $\tP(z,w) = z^{\deg(P)} P(1/z, w/z)$.
  Then $\nu(P)=-\deg(P)+\nu(\tP)$. 
  Suppose $C_1,\dots,C_k$ are exactly those 
  branches of $C$ that contain $p$. 
  Pick a defining equation (possibly transcendental)
  $\phi_i\in\hat\cO_p$ for each branch $C_i$. 
  As $P$ is reduced, we have, locally at $p$, 
  $\tP=\psi\prod_1^k\phi_i$ where $\psi(0)\ne 0$. 
  Thus 
  $\nu(\tP)=\sum_1^k\nu(\phi_i)
  =\sum m_i\,\tilde{\a}(\nu\wedge\nu_{\phi_i})$ 
  by~\cite[Proposition~3.25, Proposition~3.65]{valtree}. 
  On the other hand, for any $i>k$, 
  we have $\nu\wedge\nu_i=\affroot$ so
  $\a(\nu\wedge \nu_i)= 1$. This yields
  \begin{multline*}
    \nu(P)
    =-\deg(P)+\nu(\tP) 
    =-\sum_1^dm_i+\sum_1^km_i\,\tilde{\a}(\nu\wedge\nu_{\phi_i})
    =\\
    =-\sum_{k+1}^dm_i\,\a(\nu\wedge \nu_i)
    +\sum_1^km_i\,(\tilde{\a}(\nu\wedge\nu_{\phi_i})-1) 
    =-\sum_1^dm_i\,\a(\nu\wedge \nu_i),
  \end{multline*}
  and completes the proof.
\end{proof}
%
%
\subsection{Key polynomials}\label{sec-skp}
The rest of the appendix is devoted to the (rather technical)
proof of Theorem~\ref{thm-struct-affine}
on the structure of the valuation space $\cV_1$.
If $\nu\in\cV_0$ is any valuation centered at some point $p$
at infinity, then $\nu$ can be completely described in terms 
of its values on special polynomials:  the \emph{key polynomials}
of $\nu$, in local coordinates at $p$.
When $\nu$ belongs to $\cV_1$, the key polynomials define
curves with one place at infinity, see Lemma~\ref{L104}(ii). 
This fact holds the key to Theorem~\ref{thm-struct-affine}.
\subsubsection{Definitions}
Let us first review some material taken from~\cite{valtree}.  

Roughly speaking, an SKP is a sequence of polynomials 
$(\tU_j)_0^k$, $1\le k\le\infty$ (in some coordinates $(z,w)$)
together with a sequence $(\tbeta_j)_0^k$ of positive real numbers,
satisfying certain combinatorial relations. 
To any SKP is attached a valuation $\nu$ with $\nu(\tU_j)=\tbeta_j$ 
for all $j$. 
Conversely, to any valuation $\nu$ is attached a unique SKP.

More precisely, an SKP is a sequence 
$[(\tU_j)_0^k;(\tbeta_j)_0^k]$ such that:
\begin{itemize}
\item[(i)]
  $ \tU_0=z,\tU_1=w,\tbeta_0=1,\tbeta_1>0$;
\item[(ii)]
  if $k>1$, then $\tbeta_j\in\Q$ for $1\le j<k$,
  and $\tbeta_{j+1}>n_j\tbeta_j=\sum_{l=0}^{j-1}m_{j,l}\tbeta_l$,
  where $n_j=\min\{l\ ; \ l\tbeta_j\in\sum_0^{j-1}\Z\tbeta_i\}$ 
  and $0\le m_{j,l}<n_l$ for $1\le l<j$;
\item[(iii)]
  for $1\le j<k$:
  $\tU_{j+1}=
  \tU_j^{n_j}-\theta_j\cdot\tU_0^{m_{j,0}}\dots\,\tU_{j-1}^{m_{j,j-1}}$,
  where $\theta_j\in\C^*$.
\end{itemize}
\begin{Lemma}\label{L-keypol}
  The polynomial $\tU_j$ is irreducible
  both in $\C[z,w]$ and in $\C\{z,w\}$. In particular, the curve
  $\tU_j^{-1}(0)$ is globally irreducible in $\P^2$, and analytically
  irreducible at $(z,w)=(0,0)$.
  Moreover, if $j\ge1$, $\tU_j$ monic in $w$, 
  and its intersection multiplicity
  with $\{ z=0 \}$ is given by $d_j\=n_1\dots n_{j-1}$, with the
  convention $d_1=1$.
\end{Lemma}
We refer to~\cite[Lemma~2.4]{valtree} for a proof.

Let us indicate how an SKP determines a valuation.  
We construct by induction on $k$ a valuation $\nu_k$ associated to the
data $[(\tU_j)_0^k;(\tbeta_j)_0^k]$. 
First, $\nu_1$ is the
monomial valuation with $\nu_1(z)=\tbeta_0=1,\nu_1(w)=\tbeta_1$.
Now assume $j>1$, that the SKP $[(\tU_j)_0^k;(\tbeta_j)_0^k]$ is given
and that $\nu_1,\dots,\nu_{k-1}$ have been defined. 
Consider a polynomial $\phi\in \C [z,w]$. As $\tU_k$ is monic
in $w$, we can divide $\phi$ by $U_k$ in $\C(z)[w]$:
$\phi=\phi_0+\tU_k\psi$ with $\deg_z(\phi_0)<d_k=\deg_w(\tU_k)$ and
$\psi\in \C[z,w]$. Iterating the procedure we get a unique
decomposition $\phi=\sum_j\phi_j\tU_k^j$
with $\phi_j\in \C[z,w]$ and $\deg_z(\phi_j)<d_k$. 
Define
\begin{equation}\label{min}
  \nu_k(\phi)
  \=\min_j\nu_k(\phi_j\tU_k^j)
  \=\min_j\{\nu_{k-1}(\phi_j)+j\tbeta_k\}.
\end{equation}
It is clear that $\nu_k(\tU_k) = \tbeta_k$, so
$\nu_k(\tU_j)=\tbeta_j$ for all $j\le k$.  
One can prove that $\nu_k$ 
is indeed a valuation, see~\cite[Theorem~2.8]{valtree}.

For the converse construction, suppose we are given a valuation 
$\nu$ with $\nu(z)=1$, $\nu(w)>0$.
Set $\tbeta_0=1$, $\tbeta_1 \= \nu(w)$, and suppose
$\tU_j$, $\tbeta_j$ have been defined for $j=1, \dots, k$.  
If $\nu$ coincides with the valuation $\nu_k$ attached to 
$[(\tU_j)_0^k;(\tbeta_j)_0^k]$, then we are done. 
Otherwise, we still have $\nu(\tU_j)=\nu_k(\tU_j)=\tbeta_j$, 
and it follows from the
definition of $\nu_k$ that $\nu \ge \nu_k$. In particular the set
$\mathcal{D} \= \{ \phi\ ; \ \nu (\phi) > \nu_k (\phi) \}$ is an
ideal of the graded ring of $\nu_k$ (over $\C(z)[w]$). This graded
ring can be shown to be Euclidean, and its irreducible
elements are all of the form 
$\tU_k^{n_k}-\theta_k\cdot\tU_0^{m_{k,0}}\dots\,\tU_{k-1}^{m_{k,k-1}}$, 
where $\theta_k\in\C^*$, 
$n_k=\min\{l\ ;\ l\tbeta_k\in\sum_0^{k-1}\Z\tbeta_j\}$, 
and $n_k\tbeta_k=\sum m_{k,j}\tbeta_j$,
see~\cite[Theorem~2.29]{valtree}. 
We choose $\tU_{k+1}$ to be the generator of $\mathcal{D}$, 
and set $\tbeta_{k+1}\=\nu(\tU_{k+1})$.

\subsubsection{Global properties}
Let $[(\tU_j)_0^k ; (\tbeta_j)_0^k]$ be an SKP in coordinates
$(z,w)$.  
Write $D_j$ for the degree of $\tU_j$ as a polynomial in $(z,w)$.  
By Lemma~\ref{L-keypol} and Bezout's theorem we have 
$d_j\le D_j$ for $1\le j<k+1$.
Moreover, if we identify $\{z=0\}$ with the line at infinity in $\P^2$,
then the curve $\tU_j^{-1}(0)$ has
one place at infinity iff $D_j=d_j$.  The following technical result
exhibits one important case when this occurs.
\begin{Lemma}\label{L103}
  Consider a finite SKP $[(\tU_j)_0^{k+1};(\tbeta_j)_0^{k+1}]$. 
  Assume $\tbeta_j\le d_j$ for $1\le j\le k$ and
  $\tbeta_1+\sum_1^{k-1}(\tbeta_{j+1}-n_j\tbeta_j)<2$.
  Then $D_j=d_j$ for $1\le j\le k+1$.
\end{Lemma}
It is important to notice that the assumptions in the lemma 
involve only the indices $j\le k$ but the conclusions 
apply also for $j=k+1$.
\begin{proof}
  Recall that $\tU_1=w$, so that $d_1=D_1=1$. 
  By induction it suffices to prove $d_{k+1}=D_{k+1}$
  under the assumption $d_j=D_j$ for $1\le j<k+1$.
  Now
  \begin{equation*}
    \tU_{k+1}
    =\tU_k^{n_k}-\theta_k\cdot\tU_0^{m_{k,0}}\dots\tU_{k-1}^{m_{k,k-1}}.
  \end{equation*}

  The second term in the right hand side has degree
  $D':=\sum_0^{k-1}m_{k,j}D_j$. The first term has degree
  $n_kD_k=n_kd_k=d_{k+1}=\prod_1^kn_j$, and $\tU_k(0,w)^{n_k}$ has
  degree $d_{k+1}$ by Lemma~\ref{L-keypol}. 
  It hence suffices to prove the inequalities
  $D'\le d_{k+1}$, and $m_{k,0}\ge1$. 
  For the first, we use $m_{k,j}<n_j$ for
  $0<j<k$ together with the inductive assumptions:
  \begin{multline*}
    D'
    =m_{k,0}+\sum_1^{k-1}m_{k,j}d_j
    =\sum_0^{k-1}m_{k,j}\tbeta_j
    +\sum_1^{k-1}m_{k,j}(d_j-\tbeta_j)\le\\
    \le n_k\tbeta_k
    +\sum_1^{k-1}(n_j-1)(d_j-\tbeta_j)
    =(n_k-1)\tbeta_k+d_k-1
    +\tbeta_1+\sum_1^{k-1}(\tbeta_{j+1}-n_j\tbeta_j)<\\
    <(n_k-1)\tbeta_k+d_k+1
    \le (n_k-1)d_k+d_k+1
    =d_{k+1}+1.
  \end{multline*}
  Hence $D'\le d_{k+1}$. The fact that $m_{k,0}\ge1$ 
  is proved in a similar way. Namely:
  \begin{equation*}
    \sum_1^{k-1}m_{k,j}\tbeta_j
    \le 
    \sum_1^{k-1}(n_j-1) \tbeta_j
    \le n_{k-1}\tbeta_{k-1}-\tbeta_1<\tbeta_k,
  \end{equation*}
  so we get $m_{k,0} =  n_k \tbeta_k - \sum_1^{k-1}m_{k,j}\tbeta_j>0$.
\end{proof}

\subsubsection{The SKP in affine coordinates}
Consider a valuation $\nu\in\cV_0$ centered at $p\in L_\infty$
so that $\nu\in\tcV_p$. Pick local coordinates $(z,w)$ at $p$
such that $(X,Y)=(1/z,w/z)$ are affine coordinates on $\C^2$.
Suppose $\nu$ is given by the SKP $[(\tU_j)_0^k;(\tbeta_j)_0^k]$,
$1\le k\le\infty$.  
Define $U_j(X,Y)=X^{D_j}\tU_j(1/X,Y/X)$ for $1\le j<k+1$, 
with $D_j=\deg(\tU_j)$.
\begin{Lemma}\label{L104}
  If $A(\nu)\le0$ and $\nu(P)\le0$ for all polynomials $P$, 
  then:
  \begin{itemize}
  \item[(i)]
    $\tbeta_j\le d_j$ for $1\le j<k+1$;
    and $\a(\nu)=1-\talpha(\nu)=1-\sup_j\frac{\tbeta_j}{d_j}\ge0$;
  \item[(ii)]
    $d_j=D_j$, \ie $U_j^{-1}(0)$ has one place 
    at infinity for $1\le j<k+1$;
  \end{itemize}
  If in addition $\nu$ is quasimonomial, then:
  \begin{itemize}
  \item[(iii)]
    the SKP of $\nu$ is finite, \ie $k<\infty$, 
    the curve valuation $\nu_{U_k}\in\cV_0$ dominates $\nu$,
    and $m(\nu)=D_k=\deg(U_k)$.
\end{itemize}
\end{Lemma}
\begin{proof}
  Since $\nu(Y)\le 0$ we have $\tbeta_1=\nu(w)\le1=d_1=D_1$.
  Suppose the inequality in~(i) does not always 
  hold and pick $j$ minimal,
  $1\le j<k$ such that $\tbeta_{j+1}>d_{j+1}$.
  Let $\mu\in\tcV_p$ be the valuation defined by the truncated
  SKP $[(\tU_l)_0^j;(\tbeta_l)_0^j]$. 
  Then $\mu<\nu$ so 
  $\tA(\mu)<\tA(\nu)=A(\nu)+3\le3$.
  Now $\tA(\mu)=1+\tbeta_1+\sum_1^{j-1}(\tbeta_{l+1}-n_l\tbeta_l)$,
  so we may apply Lemma~\ref{L103} to the SKP
  $[(\tU_l)_0^{j+1};(\tbeta_l)_0^{j+1}]$
  and conclude that $D_l=d_l$ for $1\le l\le j+1$.
  But then 
  \begin{equation*}
    \nu(U_{j+1})
    =\nu(\tU_{j+1})-D_{j+1}
    =\tbeta_{j+1}-d_{j+1}
    >0,
  \end{equation*}
  contradicting the assumption $\nu(P)\le 0$ 
  for all polynomials $P$.
  Thus $\tbeta_j\le d_j$ for all $j$. 
  By~\cite[Lemma~3.32]{valtree},
  the relative skewness of $\nu$ is given by 
  $\talpha(\nu)=\sup_j\tbeta_j/d_j$.
  This shows that $\talpha(\nu)\le1$, so that
  $\a(\nu)=1-\talpha(\nu)\ge0$,
  completing the proof of~(i).
  
  Applying Lemma~\ref{L103}
  again we conclude that $d_j=D_j$ for all $j$.
  Now $d_j=D_j$ is equivalent to
  the curve $U_j^{-1}(0)$ having one place at infinity.
  Thus~(ii) holds as well.
  Finally, if $\nu$ is quasimonomial, then its
  SKP is finite: see~\cite[Definition~2.23]{valtree}.
  Moreover, the curve valuation $\nu_{\tU_k}\in\tcV_p$
  dominates $\nu$ and $\tm(\nu)=\tm(\tU_k)=d_k=D_k$:
  see~\cite[Lemma~3.42]{valtree}.
  By the definition of the multiplicity on $\cV_0$, this
  yields~(iii).
\end{proof}

\subsubsection{Proof of Theorem~\ref{thm-struct-affine}}
That $\cV_1$ is a complete subtree of $\cV_0$ is clear by definition.
By continuity, 
$A(\nu)\le0$ and $\nu(P)\le0$ for $\nu\in\cV_1$ and $P\in R$.

For~(i), pick $\nu\in\cV_1$ quasimonomial.
If $\nu=\affroot$, then $m(\nu)=1$ and $\nu$
is dominated by the pencil valuation $\nu_{|X|}$ for any 
affine function $X$, so suppose $\nu\ne\affroot$.
The key fact is that there exists a polynomial
$P$ with one place at infinity such that the 
curve valuation $\nu_P$ dominates $\nu$ and $m(\nu)=\deg(P)$.
The existence of $P$ is guaranteed by Lemma~\ref{L104}~(iii).
We claim $\nu$ is also dominated by the pencil valuation $\nu_{|P|}$. 
Indeed, both $\nu$ and $\nu_{|P|}$ are dominated by $\nu_P$,
hence they are comparable.
If $\nu\in\cV_1^o$, then $\nu<0$ on all nonconstant polynomials,
but $\nu_{|P|}(P)=0$ so $\nu<\nu_{|P|}$. 
If $\nu\in\cV_1\setminus\cV_1^o$ we obtain $\nu\le\nu_{|P|}$
by approximation.

Thus $\nu$ is dominated by the valuation associated to the 
pencil $|P|$, whose degree is $m(\nu)$. 
On the other hand, if $C$ is any curve with 
one place at infinity and $\nu\le\nu_{|C|}$, then 
$\nu<\nu_C$. Now the multiplicity of the curve valuation 
$\nu_C\in\cV_0$ equals the intersection multiplicity of $C$ with the
line at infinity, hence equals the degree of $C$ by Bezout's theorem.
Thus $m(\nu)\le m(\nu_C)=\deg(C)$, 
completing the proof of~(i).

\smallskip
Next we consider~(ii). By definition, $\a(\affroot)=1$.
For $\nu\ne\affroot$,
Lemma~\ref{L104} applies and shows that $\a(\nu)\ge0$. 
By Theorem~\ref{T101} it is then clear that
skewness restricts to decreasing parameterization of $\cV_1$
with values in $[0,1]$.

If $P$ is any nonconstant polynomial,
then~\eqref{e109} shows that
\begin{equation*}
  \nu(P)
  =-\sum_1^d m_i\,\a(\nu\wedge\nu_i)
  \le-\sum_1^d m_i\,\a(\nu)
  =-\a(\nu)\deg(P)
\end{equation*}
for any $\nu\in\cV_1$.
Equality holds iff all the curve valuations $\nu_i$
associated to irreducible local branches of $P$ dominate $\nu$.
In particular, this holds when $P$ has one place at infinity
and $\nu_P>\nu$. 
By~(i), such a polynomial $P$ exists when $\nu\in\cV_1$
is quasimonomial. 
When $\nu$ is not quasimonomial,
we may find $\nu_n\in\cV_1$ quasimonomial increasing to $\nu$.
For each $n$ there exists a polynomial $P_n$ with
$\nu_n(P_n)=-\a(\nu_n)\deg(P_n)$.
Thus
\begin{equation*}
  \frac{\nu(P_n)}{\deg(P_n)}
  \ge\frac{\nu_n(P_n)}{\deg(P_n)}
  =-\a(\nu_n)
  \to-\a(\nu),
\end{equation*}
which completes the proof of~(ii).

All the statements in~(iii) are consequences of Theorem~\ref{T101} so
we now turn to~(iv).
Consider an end $\nu$ of $\cV_1$.
We have seen that $A(\nu)\le 0$, which implies that
$\nu$ cannot be a curve valuation.
Hence $\nu$ is either infinitely singular or quasimonomial.
When $\nu$ is infinitely singular, we are in case~(c).  Otherwise,
$\nu$ is quasimonomial, hence $\a(\nu)=-\sup\nu(P)/\deg(P)$ is
attained for some polynomial $P$ having one place at infinity, 
by~(ii). 
It thus defines a pencil valuation $\nu_{|P|}$, and we have seen above
that $\nu_{|P|}\ge\nu$.  
If $\nu<\nu_{|P|}$, then $\mu(Q)\le 0$ for all polynomials $Q$ 
and all $\mu\in[\nu,\nu_{|P|}]$. By assumption, $\nu$ is an end of
$\cV_1$ thus $A(\nu)=0$. We are in case~(b). 
The remaining case is
when $\nu=\nu_{|P|}$, \ie $\nu$ is a pencil valuation. 
But $A(\nu)\le0$, so Proposition~\ref{L102} implies that $\nu$ is a rational
pencil valuation. We are then in case~(a).
Finally, one checks that any valuation
of the type~(a),~(b), or~(c) is an end of $\cV_1$ and
does not belong to $\cV_1^o$.
This completes the proof of~(iv) and hence the
proof of Theorem~\ref{thm-struct-affine}.
%
%
%
%

%
%
\end{document}